\RequirePackage{silence}
\documentclass[EJP,preprint,NODS,nobera]{ejpecp}

\usepackage{lmodern}

\usepackage[T1]{fontenc}
\usepackage{enumitem}
\usepackage{placeins}
\usepackage{needspace}
\usepackage{etoolbox}
\usepackage{graphicx}
\usepackage{xcolor}
\usepackage{tikz}
\usetikzlibrary{arrows.meta,calc,patterns.meta}
\usepackage{yhmath}

\BeforeBeginEnvironment{theorem}{\Needspace{5\baselineskip}}
\BeforeBeginEnvironment{proposition}{\Needspace{5\baselineskip}}
\BeforeBeginEnvironment{lemma}{\Needspace{5\baselineskip}}
\BeforeBeginEnvironment{corollary}{\Needspace{5\baselineskip}}
\BeforeBeginEnvironment{definition}{\Needspace{4\baselineskip}}
\BeforeBeginEnvironment{conjecture}{\Needspace{4\baselineskip}}
\BeforeBeginEnvironment{proof}{\Needspace{4\baselineskip}}

\makeatletter
\let\@EJPLOGO\@empty
\makeatother

\newcommand{\R}{\mathbb{R}}
\newcommand{\Z}{\mathbb{Z}}
\newcommand{\N}{\mathbb{N}}

\newcommand{\pP}{\mathbb{P}}
\newcommand{\E}{\mathbb{E}}
\newcommand{\CE}[2]{\E\Cb{#1\,\middle|\,#2}}
\newcommand{\CP}[2]{\pP\bb{#1\,\middle|\,#2}}

\newcommand{\bb}[1]{\left(#1\right)}
\newcommand{\Cb}[1]{\left[#1\right]}
\newcommand{\cb}[1]{\left\{#1\right\}}
\newcommand{\norm}[1]{\left\Vert#1\right\Vert}
\newcommand{\abs}[1]{\left|#1\right|}

\newcommand{\ind}{\mathbf{1}}
\newcommand{\vn}{\varnothing}
\newcommand{\mb}[1]{\mathbf{#1}}
\newcommand{\mbs}[1]{\boldsymbol{#1}}
\newcommand{\Lcal}{\mathcal{L}}
\newcommand{\Dcal}{\mathcal{D}}
\newcommand{\Ncal}{\mathcal{N}}

\newcommand{\pw}[4]{\left\{\begin{array}{cl}
  #1 &\text{ if } #2 \\
  #3 &\text{ if } #4
\end{array}\right.}
\newcommand{\pwthree}[6]{\left\{\begin{array}{cl}
  #1 &\text{ if } #2 \\
  #3 &\text{ if } #4 \\
  #5 &\text{ if } #6
\end{array}\right.}

\SHORTTITLE{Ergodicity without monotonicity}
\TITLE{Ergodicity without monotonicity:\\patch-positive spin systems}
\AUTHORS{%
Maciej~G{\l}uchowski\footnote{Department of Mathematics, University of
California, Los Angeles, CA 90095, USA. Corresponding author: Maciej
G{\l}uchowski. \mbox{\EMAIL{mgluchowski@math.ucla.edu}};
\mbox{\EMAIL{gmenz@math.ucla.edu}}}\orcid{0009-0009-4239-9674}
\and
Georg~Menz\footnotemark[1]\orcid{0000-0001-9295-5063}}
\hypersetup{
  pdftitle={Ergodicity without monotonicity: patch-positive spin systems},
  pdfauthor={Maciej G{\l}uchowski, Georg Menz}
}
\AMSSUBJ{60K35}
\AMSSUBJSECONDARY{60E15; 82C22}

\ABSTRACT{%
Monotonicity is the standard tool for proving that an interacting particle
system is ergodic. What can replace it when the system is not monotone and its
invariant measure is unknown? The simplest such systems are kinetically
constrained models such as FA-1f and East exposed to noise. In these models
vacancies enable updates, so adding vacancies should help the system forget
its initial state. We prove
this for a broad class of spin systems, which we call patch-positive. If
vacancies are created spontaneously at a uniformly positive rate, the system
converges to the same equilibrium from a large class of initial states, and
from every initial state if the threshold density is at most~$\tfrac12$ or
the interactions are oriented. In particular, the East model with noise is
ergodic at every density, and we obtain new cases of the positive-rates
conjecture. The key idea is to compare the dual processes of the system with
and without the added deaths.}

\begin{document}

\makeatletter
\patchcmd{\l@section}{1.0em}{0.4em}{}{\PackageWarning{main}{TOC spacing patch failed}}
\makeatother
\setcounter{tocdepth}{1}
\begingroup\small
\tableofcontents
\endgroup
\newpage

\section{Introduction}
\label{sec:introduction}

Interacting particle systems are stochastic models in which local random
updates at many sites combine to produce collective behaviour. We study spin
systems, continuous-time Markov processes on~$\cb{0,1}^{\Lambda}$, where
$\Lambda$ is a countable set of sites, in which one spin is updated at a time
and the update rule may depend on nearby spins; see~\cite{Liggett1985,Swart2026} for the
general theory. A basic question is whether a given spin system is ergodic,
that is, whether it has a unique invariant measure~$\pi$ and the law of the
process converges to~$\pi$ from every initial configuration. This question is
difficult in general, and the answer depends on the local update rule in a
subtle way.

The classical tool for this question is monotonicity. A spin system is
monotone, traditionally called attractive, if processes started from
configurations ordered coordinatewise,~$\eta\leq\xi$, can be coupled so that
the order is preserved at all later times~\cite{Harris1977},
\cite[Chapter~III, Theorem~2.2]{Liggett1985}; equivalently, stochastically
ordered initial laws have stochastically ordered
evolutions~\cite[Chapter~II, Theorem~2.2 and Definition~2.3]{Liggett1985}.
Examples are independent Bernoulli flips, ferromagnetic stochastic Ising
models and the contact process. For monotone spin systems, evolutions started
from the all-$0$ and all-$1$ configurations are monotone in time and converge
to the lower and upper invariant measures~$\pi^-$ and~$\pi^+$, every other
evolution is sandwiched between them, and ergodicity reduces to showing
that~$\pi^-=\pi^+$~\cite[Chapter~III, Theorem~2.3 and
Corollary~2.4]{Liggett1985}. The question that motivates this paper is what
can be done without monotonicity, in particular when no explicit formula for
the candidate invariant measure is available.

The simplest way to break monotonicity is to constrain a monotone system. We
start from independent Bernoulli refreshes: every site is resampled
from~$\operatorname{Bern}(p)$ at rate one, regardless of the other spins.
These dynamics are monotone, since two ordered initial configurations driven
by the same clocks and the same coins stay ordered. We then constrain the
dynamics: an update at a site is carried out only if a local condition on the
neighbouring spins is satisfied. This gives the kinetically constrained spin
models (KCSMs)~\cite{FredricksonAndersen1984,CancriniEtAl2008,HartarskyToninelli2025}.
We use the convention that state~$0$ is vacant and facilitating. In the
Fredrickson--Andersen one-spin facilitated model (FA-1f), an update at~$i$ is
allowed when at least one neighbouring site is vacant. The East model imposes
the same constraint on an oriented neighbourhood graph; on~$\Z$, a site can
update only if its right neighbour is vacant. KCSMs are toy models for the
freezing and melting of glasses: fully occupied regions are frozen and can
only be unfrozen by vacancies arriving from outside.

The constraint destroys monotonicity. For example, in FA-1f or East
let~$j\in N(i)$, $\eta(x)=\ind(x\neq i,j)$ and~$\xi(x)=\ind(x\neq i)$.
Then~$\eta\leq\xi$, but an update at~$i$ that refreshes to~$1$ is legal
in~$\eta$ and blocked in~$\xi$, so afterwards~$\eta(i)=1>\xi(i)=0$. The
reason is general. Monotonicity requires the~$0\to1$ rate to increase and
the~$1\to0$ rate to decrease as spins are added~\cite[Chapter~III,
Definition~2.1 and Theorem~2.2]{Liggett1985}, but in a KCSM both rates are
proportional to the same constraint. Hence a KCSM is monotone only if its
constraint is trivial. On the other hand, the constraint does not change the
equilibrium: the Bernoulli product measure~$\mu_p$ is reversible for every
KCSM with refresh density~$p$. In Section~\ref{sec:setup} we give precise definitions of the models we use.

Since vacancies facilitate updates, adding vacancies should unfreeze
blocked regions and remove the bottlenecks that slow down the relaxation of
KCSMs. It is therefore natural to ask whether adding noise, for example
independent~$1\to0$ transitions or independent Bernoulli refreshes at a small
rate~$\varepsilon>0$, makes a KCSM ergodic. More generally one may ask
whether adding noise to an interacting particle system causes, facilitates,
or at least preserves ergodicity.
{For example, for finite-range spin systems on~$\Z^d$ with bounded flip rates, adding sufficiently
fast independent refreshes gives ergodicity~\cite[Chapter~I, Theorem~4.1]{Liggett1985}.}

Even for KCSMs this question is not easy. The constraint has already removed
monotonicity. The other structure on which the analysis of KCSMs rests is the
explicit reversible measure~$\mu_p$, which makes spectral-gap and
Dirichlet-form methods available~\cite{CancriniEtAl2008}. Noise that is not
itself reversible with respect to~$\mu_p$ destroys this structure, and the
invariant measures of the perturbed system are not known explicitly. Thus
KCSMs with added noise are a natural test case for the problem of proving
ergodicity for a non-monotone system without an explicit formula for its
invariant measure.

Even without noise, the out-of-equilibrium behaviour of KCSMs is an active
area of research; see the monograph~\cite{HartarskyToninelli2025} and the
survey~\cite{MartinelliShapiraToninelli2025}. For FA-1f, Cancrini, Martinelli, Roberto and Toninelli proved
that the spectral gap is positive for every~$p<1$~\cite{CancriniEtAl2008}.
In infinite volume this does not suffice to control the evolution from
initial laws far from~$\mu_p$, such as deterministic configurations, and no
uniform logarithmic Sobolev inequality is available for these
models~\cite[Corollary~3.13]{HartarskyToninelli2025}.
Earlier work establishes convergence to equilibrium from out-of-equilibrium
initial laws in particular regimes, for instance for FA-1f in the high
facilitating regime~\cite{BlondelEtAl2013}, for every KCSM when the
equilibrium density of vacancies is close to one~\cite{HartarskyToninelli2024},
and for the East
model~\cite{CancriniEtAl2010,Blondel2013,Mareche2019}. See
also~\cite{MountfordValle2019,BlondelDeshayesToninelli2019,MartinelliShapiraToninelli2025}.
The general question remains open. Since the all-$1$ configuration is frozen,
convergence can only be expected from initial states containing a vacancy.

\begin{conjecture}[{\cite[Conjecture~1]{MartinelliShapiraToninelli2025}}]
\label{conj:fa1f}
{Consider FA-1f on~$\Z$ with refresh density~$p\in(0,1)$.
For every initial law~$\nu$ that gives probability one to configurations
with at least one vacancy, we have~$(\nu P_t)(f)\to\mu_p(f)$ for every
local function~$f$.}
\end{conjecture}

Adding noise also connects KCSMs with the positive-rates problem. For
finite-range spin systems in one dimension with positive rates, the natural
conjecture that they are ergodic is false in general, as shown by
G\'acs~\cite{Gacs1986,Gacs2001}; see also Gray's guide~\cite{Gray2001}.
For simple spin systems, binary systems on~$\Z$ in which the update at~$i$
depends only on the spins at~$i$ and~$i+1$, the problem remains open, and
earlier work covers large parts of the parameter
space~\cite{GluchowskiMenz2025,GluchowskiMenz2026}. The part that remains
open consists of East models with added noise.

\begin{conjecture}
\label{conj:east-deaths}
Let~$p\in(0,1)$,~$\varepsilon>0$ and~$s\in[0,1)$. The East model on~$\Z$ with
refresh density~$p$, in which in addition every site is resampled
from~$\operatorname{Bern}(s)$ at rate~$\varepsilon$, is ergodic.
\end{conjecture}

For~$s\in(0,1)$ these systems have positive rates, and the open part of the
positive-rates conjecture for simple spin systems consists of such models
with~$s<p$; see Remark~\ref{rem:simple-remaining-region}. The case~$s=0$, in
which the noise consists of independent~$1\to0$ transitions only, lies at
the boundary of the positive-rates region and is the most constrained case.

We introduce a class of spin systems that we call \emph{patch-positive}.
In Definition~\ref{def:patch-positivity} we specify the class by an explicit
sign condition on the local update rates. The class contains FA-1f and East,
but not FA-$j$f with~$j\geq2$. It also contains models outside the KCSM
class, such as the biased annihilating branching process (BABP). Each
patch-positive system has a threshold profile~$\mb p^\star\in[0,1]^\Lambda$;
for FA-1f and East it equals the refresh density~$p$.

In our main result, Theorem~\ref{thm:common-invariant-limit}, we show that
noise makes patch-positive systems ergodic on a large class of initial laws. Suppose that the neighbourhood graph has polynomial growth and that the
system contains independent~$1\to0$ transitions, which we call pure deaths,
at a uniformly positive rate~$\varepsilon$. Then there is an invariant
measure~$\pi$ such that the law of the process converges to~$\pi$ from every
initial law in an explicit class~$\mathcal M_-$, and the convergence is
exponential, uniformly over the initial laws in suitable subclasses. We define
the class~$\mathcal M_-$ through the joint moments
\begin{equation*}
\mu\bb{\chi_A^\star},
\qquad
\chi_A^\star(\eta)
=
\prod_{i\in A}\bb{\eta(i)-p_i^\star},
\end{equation*}
of the spins centered at the threshold profile, see
Definition~\ref{def:m-minus}. It contains every configuration with finitely
many vacancies, every product measure with densities at
least~$\max\cb{2p_i^\star-1,0}$, and, if~$\mb p^\star\leq\mbs{\tfrac12}$,
every probability measure.

In two corollaries we obtain convergence uniformly over all initial
configurations.
If~$\mb p^\star\leq\mbs{\tfrac12}$, the system is uniformly exponentially
ergodic (Corollary~\ref{cor:uniform-exponential-ergodicity}). If the
neighbourhood graph is orientable, as for the East model, the system is
uniformly ergodic without any restriction on~$\mb p^\star$
(Corollary~\ref{cor:oriented-ergodicity}).

For KCSMs we add noise through independent Bernoulli refreshes at
rate~$\varepsilon$ with density~$s_i$ at site~$i$. The perturbed FA-1f is
patch-positive if and only if~$s_i\geq p$ for every~$i$, and it then contains
pure deaths at rate~$\varepsilon\inf_i(1-s_i)$. Hence FA-1f with such noise
is uniformly exponentially ergodic if~$p\leq\tfrac12$, and the East model
with such noise is uniformly ergodic for every~$p\in(0,1)$
(Corollary~\ref{cor:perturbed-fa-babp-ergodicity}). The same results hold for BABP. For simple spin systems we use the orientable case to identify a new region
of the parameter space where the positive-rates conjecture holds
(Corollary~\ref{cor:simple-positive-rates}).

Thus Conjecture~\ref{conj:east-deaths} holds for~$s\geq p$. The
restriction~$s_i\geq p$ is essential for our method, so the conjecture
remains open for~$s<p$, including the case~$s=0$: adding only pure deaths to
FA-1f or East leaves the patch-positive class for every~$p\in(0,1)$. Conversely, removing pure deaths from a patch-positive
system keeps it patch-positive (Proposition~\ref{prop:death-removal}).
We leave open Conjecture~\ref{conj:fa1f}, which concerns the dynamics without
noise.

The main construction is a representation of the action of the Markov
semigroup on monomials. For~$\chi_A(\eta)=\prod_{i\in A}\eta(i)$, the
starting point is the set-valued Feynman--Kac representation
\begin{equation}
\label{eq:intro-signed-duality}
P_t\chi_A(\eta)
=
\E_A\Cb{
\sigma_t
\exp\bb{\int_0^t V(A_s)\,ds}
\chi_{A_t}(\eta)
}.
\end{equation}
Here~$A_s$ is the finite set-valued dual process,~$\sigma_s\in\cb{+,-}$, and~$V$
is the Feynman--Kac potential. Set-valued dual constructions for spin systems
are standard; see, for example,
\cite{Harris1974,Harris1978,Griffeath1979,HolleyStroock1979,Gray1986,JansenKurt2014,SturmSwart2018,SturmSwartVollering2020}.
For non-monotone systems the right-hand side
of~\eqref{eq:intro-signed-duality} is signed, which is an obstruction for
convergence arguments. 

The new step is to reveal only part of the graphical construction. We define a~$\sigma$-algebra
$\mathcal G_t$ that records interactions between sites that occur in the
dual up to time~$t$, including their times and the sites they connect, while
forgetting their type and the purely local events between them. Once this
information is fixed, the recorded interactions divide each site line~$\cb{i}\times[0,t]$ into
consecutive time intervals. We call these intervals \emph{patches}. Thus a
patch is a piece of one site line between two consecutive recorded
interactions, or between such an interaction and the observation time.
Conditioning on~$\mathcal G_t$, the forgotten randomness has a conditional
law that factorizes over the patches, and the signed random variable
in~\eqref{eq:intro-signed-duality} is itself a product of local factors.
Consequently, schematically, writing~$\mathcal B_t$ for the patches
ending before~$t$ and~$\mathcal E_t$ for those reaching~$t$,
\begin{equation}
\label{eq:intro-patch-factorization}
\E_A\Cb{
\sigma_t
\exp\bb{\int_0^t V(A_s)\,ds}
\chi_{A_t}(\eta)
\,\middle|\,
\mathcal G_t
}
=
\prod_{P\in\mathcal B_t}C(P)
\prod_{P\in\mathcal E_t}C\bb{\eta\bb{i(P)},P}.
\end{equation}
In Section~\ref{sec:patches} we define patches and their contributions
precisely. The conditioning allows for local averaging, which
gives explicit contributions and captures the cancellations of the sign
process in~\eqref{eq:intro-signed-duality}.

Patch positivity means exactly that every contribution~$C(P)$ of a patch
ending before time~$t$ is nonnegative. A patch reaching time~$t$
contributes~$C(z,P)$, an affine function of the remaining spin argument~$z$,
and the threshold~$p_i^\star$ is the smallest~$p\in[0,1]$ for which
$C(p,P)\geq0$ for every such patch at~$i$. For initial laws with nonnegative
centered moments~$\mu(\chi_A^\star)$, the representation then involves only
nonnegative quantities. These initial laws form the subclass~$\mathcal M_\star$ of~$\mathcal M_-$.

The main technical idea is a natural coupling of two dual processes.
Removing the pure deaths from the system does not change the dual
process, but only increases the
Feynman--Kac potential. This holds for every spin system, but in general the
dual weights carry signs and the coupling gives no inequality. Patch
positivity makes the comparison termwise:
every patch contribution of the original system is dominated by the
corresponding contribution of the system without the added deaths, and on
patches that end with an outgoing interaction the original contribution is
smaller by the factor~$e^{-\varepsilon u}$, where~$u$ is the length of the
patch. The ancestry of a successful interaction at time~$u$ consists of such
patches, with lengths adding up to~$u$. Since the system without the added
deaths still gives an integrable bound, the event with a successful interaction after a late time~$T$ contributes to~$P_t\chi_A$
only~$O(e^{-\varepsilon T})$. On the remaining event the dependence on the
initial law decays because end patches relax, and polynomial growth together
with finite propagation controls the number of end patches. Thus the two
coupled duals allow to push through an analogue of the classic argument for ergodicity via duality: use the duality to split~$P_t\chi_A$ into one term that decays and another that doesn't depend on initial configuration.

We can also use the patch representation for systems that are not
patch-positive. Patch
positivity makes the patch contributions nonnegative and allows termwise
comparison. A different possibility is quantitative contraction even when the
contributions have mixed signs: if the relevant patch factors satisfy a
uniform bound of the form~$\abs{C(z,P)}\leq1-\varepsilon$, then realizations
with many patches have small contribution, while realizations with few
patches have long patches whose terminal contributions depend only weakly on
the initial configuration. Models such as FA-2f, which are not
patch-positive, are a natural test case for this possibility.

In Corollaries~\ref{cor:uniform-exponential-ergodicity}
and~\ref{cor:oriented-ergodicity} we reach all initial configurations in two
ways. In the first we use the fact that~$\mathcal M_-$ contains every
probability measure when~$\mb p^\star\leq\mbs{\tfrac12}$. In the second we
use a burn-in argument, which we explain in Remark~\ref{rem:threshold-one-half}.

The duality in~\eqref{eq:intro-signed-duality} belongs to the standard
graphical and set-valued duality framework for spin systems
\cite{Harris1974,Harris1978,Griffeath1979,HolleyStroock1979,Gray1986,JansenKurt2014,SturmSwart2018,SturmSwartVollering2020}.
The new ingredient is the choice of the smaller~$\sigma$-algebra~$\mathcal G_t$.
It is this conditional factorization, rather than the duality itself, that
exposes the local cancellations used in our arguments. Information
percolation also records updates that carry information
\cite{LubetzkySly2016}, but uses them for a different decomposition.
Sly and Yang~\cite{SlyYang2026} develop a "sausage" decomposition of backward
information-percolation for the non-monotone Potts model, using
signed influence weights to obtain cutoff.
Classical comparison results for monotone spin systems are based on couplings
that preserve the coordinatewise order on configurations and on related
association inequalities \cite{Harris1977,Liggett1985,Liggett2006}. Duality
can also yield correlation inequalities in particular interacting particle
systems \cite{GiardinaRedigVafayi2010}. Patch-positive systems need not admit such couplings. We compare their dual
processes.

We state the results before introducing the patches so readers can read
the results independently of the construction. In Section~\ref{sec:setup}
we introduce the spin-system framework and the examples. In
Section~\ref{sec:main-results} we define patch positivity, state the main
theorem and its consequences, and apply them to FA-1f, the East model, BABP
and simple spin systems. In Section~\ref{sec:probabilistic-interpretation}
we introduce the dual process and the patches. In Section~\ref{sec:proofs}
we prove the results. 

\section{Setup and examples}
\label{sec:setup}

Here we fix the spin-system setup and notation for the rest of the paper.
We describe the dynamics through a neighbourhood-dependent update rate and
refresh parameter, which directly extends the usual KCSM formulation. We also give examples of spin models that are particularly relevant to this article. For a comprehensive introduction to spin systems, with proofs of existence and basic properties, we point to~\cite{Liggett1985}.

Let~$\Lambda$ be a countable set of sites and let
$\Omega=\cb{0,1}^{\Lambda}$ be the configuration space. We write~$\cdot \Subset \Lambda$ for a finite set of sites. For each~$i\in\Lambda$,
fix a  neighbourhood~$N(i)\Subset\Lambda\setminus\cb{i}$ and assume
\begin{equation*}
\sup_{i\in\Lambda}\abs{N(i)}<\infty.
\end{equation*}
For~$\eta\in\Omega$, $i\in\Lambda$, and~$x\in\cb{0,1}$,
let~$\eta^{i,x}$ be the configuration obtained from~$\eta$ by setting the spin
at~$i$ equal to~$x$, and let~$\eta^i$ be the configuration obtained by flipping
the spin at~$i$. We call a function~$f:\Omega\to\R$ local if it depends on only finitely many
coordinates. For~$p\in[0,1]$, define the Bernoulli refresh operator at~$i$ by
\begin{equation*}
E_{i,p}f(\eta)
:=
pf(\eta^{i,1})+(1-p)f(\eta^{i,0}).
\end{equation*}

\begin{definition}[Spin system]
We specify a spin system by functions
\begin{equation*}
\lambda_i:\Omega\to[0,\infty),
\qquad
\rho_i:\Omega\to[0,1],
\qquad i\in\Lambda,
\end{equation*}
which depend only on the spins in~$N(i)$ and satisfy
$\sup_{i,\eta}\lambda_i(\eta)<\infty$. Its generator acts on local functions
by
\begin{equation}
\label{eq:spin-generator}
\Lcal f(\eta)
=
\sum_{i\in\Lambda}
\lambda_i(\eta)
\bb{E_{i,\rho_i(\eta)}-I}f(\eta).
\end{equation}
Thus, at rate~$\lambda_i(\eta)$ the state of site~$i$ is updated, and upon update its new state is a~$\operatorname{Bern}(\rho_i(\eta))$ random variable. We call~$\lambda_i$ the update rate and~$\rho_i$ the refresh parameter.
\end{definition}

\begin{remark}
We can also write the generator~\eqref{eq:spin-generator} in the
usual flip-rate form
\begin{equation*}
\Lcal f(\eta)
=
\sum_{i\in\Lambda}
c_i(\eta)\bb{f(\eta^i)-f(\eta)},
\end{equation*}
where
\begin{equation*}
c_i(\eta)
=
\lambda_i(\eta)
\bb{
\bb{1-\eta(i)}\rho_i(\eta)
+
\eta(i)\bb{1-\rho_i(\eta)}
}.
\end{equation*}
Conversely, we can describe any spin system with flip rates~$c$ by an
update rate and refresh parameter. Since the two rates
$c_i(\eta^{i,0})$ and~$c_i(\eta^{i,1})$ depend only on the spins in~$N(i)$,
we obtain the refresh representation
\begin{equation*}
\lambda_i(\eta)
=
c_i(\eta^{i,0})+c_i(\eta^{i,1}),
\qquad
\rho_i(\eta)
=
\frac{c_i(\eta^{i,0})}{\lambda_i(\eta)}
\end{equation*}
when~$\lambda_i(\eta)>0$; if~$\lambda_i(\eta)=0$, the value of
$\rho_i(\eta)$ is irrelevant. 
\end{remark}

The uniform bound on the update rates gives a Feller process
$\cb{\eta_t}_{t\geq0}$ by the standard Harris construction
\cite{Harris1978,Liggett1985}. We define its semigroup~$(P_t)_{t\geq0}$ by
\begin{equation*}
P_tf(\eta)=\E_\eta\Cb{f(\eta_t)}.
\end{equation*}
For an initial law~$\mu$, we write~$\mu P_t$ for the law at time~$t$, so that
$(\mu P_t)(f)=\mu(P_tf)$. We also use the shorthand~$\mu(f) = \int_\Omega f d\mu$.

\begin{definition}
For probability measures~$\mu_t$ and~$\mu$ on~$\Omega$, we write
$\mu_t\Rightarrow\mu$ if
\begin{equation*}
\mu_t(f)\longrightarrow\mu(f)
\end{equation*}
for every local function~$f$. A probability measure~$\pi$ is invariant if
$\pi P_t=\pi$ for every~$t\geq0$. The spin system is \emph{ergodic} if there
is an invariant measure~$\pi$ such that
\begin{equation*}
\mu P_t\Rightarrow\pi
\end{equation*}
for every initial law~$\mu$. It is \emph{uniformly ergodic} if, for every local
function~$f$,
\begin{equation*}
\sup_{\eta\in\Omega}
\abs{P_tf(\eta)-\pi(f)}
\longrightarrow0.
\end{equation*}
It is \emph{uniformly exponentially ergodic} if there is~$\gamma>0$ such that,
for every local function~$f$, some~$C_f<\infty$ satisfies
\begin{equation*}
\sup_{\eta\in\Omega}
\abs{P_tf(\eta)-\pi(f)}
\leq C_fe^{-\gamma t}.
\end{equation*}
\end{definition}

\begin{remark}
Since~$\Omega$ is compact in the product topology, convergence against local
functions is equivalent to weak convergence of probability measures.
\end{remark}

We now give a few examples of familiar spin systems and how to represent them using the refresh convention.

\begin{example}[Independent refreshes]
The simplest spin system has~$\lambda_i\equiv1$ and~$\rho_i\equiv p_i$:
every site is refreshed from~$\operatorname{Bern}(p_i)$ at rate one,
independently of all other spins. Its unique invariant measure is the
Bernoulli product measure with profile~$(p_i)_{i\in\Lambda}$. These dynamics
are monotone and serve as the starting point for the kinetically constrained
models below.
\end{example}

\begin{example}[Kinetically constrained spin models]
In a kinetically constrained spin model (KCSM), the refresh parameter is
constant,
\begin{equation*}
\rho_i(\eta)=p,
\qquad p\in(0,1),
\end{equation*}
while the update rate~$\lambda_i$ is a constraint depending on the neighbouring
spins. We use the convention that state~$0$ is vacant and state~$1$ is occupied. The constraint is typically that at least one of a few prescribed subsets of~$N(i)$ is filled with vacancies. 
Whenever the constraint at~$i$ is satisfied, the spin at~$i$ is refreshed
from the Bernoulli law with parameter~$p$. This construction ensures that~$\otimes_{i\in \Lambda}\operatorname{Bern}(p)$ is a reversible invariant measure.
\end{example}

\begin{example}[Fredrickson-Andersen one spin facilitated model (FA-1f)]
For FA-1f~\cite{FredricksonAndersen1984,CancriniEtAl2008,HartarskyToninelli2025},
assume that~$N(i)\neq\vn$ for every~$i\in\Lambda$ and set
\begin{equation*}
\lambda_i(\eta)=1-\prod_{j\in N(i)}\eta(j),
\qquad
\rho_i(\eta)=p.
\end{equation*}
Thus a site can update whenever at least one site in its neighbourhood is
vacant. This model is a standard example of a KCSM.
\end{example}

\begin{example}[East model]
Let~$\Lambda=\Z^d$ and let~$\cb{e_k}_{k\leq d}$ be the standard basis of
$\Z^d$. The East model is a particular case of FA-1f, with
\begin{equation*}
N(i)=\cb{i+e_k:k\leq d}.
\end{equation*}
The neighbourhoods create an oriented lattice where sites only interact one-way. 
\end{example}

\begin{example}[Biased annihilating branching process]
The biased annihilating branching process (BABP)~\cite{NeuhauserSudbury1993}
has
\begin{equation*}
\lambda_i(\eta)
=
\sum_{j\in N(i)}\bb{1-\eta(j)},
\qquad
\rho_i(\eta)=p.
\end{equation*}
Thus the update rate is the number of vacant neighbours. In the complementary
variables~$\xi(i)=1-\eta(i)$, vacancies are particles.
\end{example}

We also consider spin systems obtained by adding noise in the form of independent refreshes to
the dynamics. 

\begin{definition}[Independent refresh noise]
For~$\mb s=(s_i)_{i\in\Lambda}\in[0,1]^\Lambda$ and~$\varepsilon\geq0$, define
\begin{equation*}
\Ncal^{\mb s}f
:=
\sum_{i\in\Lambda}\bb{E_{i,s_i}-I}f,
\qquad
\Lcal^{\varepsilon,\mb s}
:=
\Lcal+\varepsilon\Ncal^{\mb s}.
\end{equation*}
The generator~$\Lcal^{\varepsilon,\mb s}$ is again of the form
\eqref{eq:spin-generator}, with update rate and refresh parameter
\begin{equation*}
\lambda_i^{\varepsilon,\mb s}
=
\lambda_i+\varepsilon,
\qquad
\rho_i^{\varepsilon,\mb s}
=
\frac{\lambda_i\rho_i+\varepsilon s_i}
{\lambda_i+\varepsilon}
\end{equation*}
when~$\lambda_i+\varepsilon>0$. At zero update rate the refresh parameter is irrelevant.
\end{definition}

\begin{definition}[Pure-death component]
We say that~$\Lcal$ contains a pure-death component of rate~$\varepsilon>0$ if
there is another spin-system generator~$\widetilde\Lcal$ such that
\begin{equation}
\label{eq:pure-death-component}
\Lcal
=
\widetilde\Lcal
+
\varepsilon\Ncal^{\mb 0}.
\end{equation}
\end{definition}

\begin{remark}
Equivalently, the~$1\to0$ transition rate is at least~$\varepsilon$ at every
site and in every configuration:
\begin{equation*}
c_i(\eta^{i,1})
=
\lambda_i(\eta)\bb{1-\rho_i(\eta)}
\geq\varepsilon.
\end{equation*}
\end{remark}

Some of our results require assumptions on geometry of~$\Lambda$. For~$i\in\Lambda$, we set
\begin{equation*}
N_*(i):=N(i)\cup\cb{i},
\end{equation*}
and for~$A\subseteq\Lambda$ set
\begin{equation*}
N_*(A)=\bigcup_{i\in A}N_*(i),
\qquad
B(A,0)=A,
\qquad
B(A,n+1)=N_*\bb{B(A,n)}.
\end{equation*}
For~$i\in\Lambda$, write~$B(i,n)=B(\cb{i},n)$. 

\begin{definition}[Polynomial growth]
    We say that the
neighbourhood graph has polynomial growth if there are~$K,D<\infty$ such that
\begin{equation}
\label{eq:polynomial-growth}
\abs{B(i,n)}
\leq
K(1+n)^D
\qquad
\text{for every }i\in\Lambda\text{ and }n\in\N.
\end{equation}
\end{definition}

\begin{definition}[Orientability]
\label{def:orientability}
A set~$V\subseteq\Lambda$ is \emph{autonomous} if~$N_*(V)=V$. The
neighbourhood graph is \emph{orientable} if every finite set~$A\subseteq\Lambda$
is contained in a finite set~$\bar A$ for which there is
$U\subseteq\Lambda\setminus\bar A$ such that
\begin{equation*}
N_*(U)=U,
\qquad
N_*(\bar A\cup U)=\bar A\cup U.
\end{equation*}
For a given choice of~$\bar A$ and~$U$ for~$A$, we write~$U(A)$
for the chosen set~$U$.
\end{definition}

\begin{example}
\label{ex:east-orientable}
The neighbourhood graph of the East model is orientable. Indeed, for any
$A\Subset\Z^d$, choose a finite box
\begin{equation*}
\bar A=\prod_{k=1}^d \cb{a_k, a_k +1,\dots,b_k}
\end{equation*}
containing~$A$, and set
\begin{equation*}
U=\prod_{k=1}^d \cb{a_k, a_k+1\dots} \setminus \bar{A}.
\end{equation*}
Then~$N_*(U)=U$ and~$N_*(\bar A\cup U)=\bar A\cup U$.
\end{example}

For~$A\Subset\Lambda$, define the monomial
\begin{equation*}
\chi_A(\eta)
:=
\prod_{i\in A}\eta(i),
\qquad
\chi_\vn:=1.
\end{equation*}
The family~$\cb{\chi_A:A\Subset \Lambda}$ is a basis for the local
functions. Thus, every local function~$g$ has a unique
expansion
\begin{equation}
\label{eq:rate-expansion}
g(\eta)
=
\sum_{S\Subset \Lambda}
\widehat g(S)\chi_S(\eta).
\end{equation}
We call~$\widehat g(S)$ the coefficient of~$\chi_S$ in this expansion.

 Finally, for a profile
$\mb p=(p_i)_{i\in\Lambda}\in[0,1]^\Lambda$, let~$\mu_{\mb p}$ be the Bernoulli product measure~$\otimes_{i\in\Lambda}\operatorname{Bern}(p_i)$.
For a constant profile, we write~$\mu_p$.
\section{Main results}
\label{sec:main-results}

We now state the main results. We first introduce a condition on the update rates that we call patch
positivity, together with a threshold profile and the class of initial laws
for which our convergence result holds. The name comes from the representation we develop in Sections~\ref{sec:dual}
and~\ref{sec:patches}. In Proposition~\ref{prop:patch-contribution-positivity}
we give the condition a probabilistic interpretation. We
check the condition for FA-1f, the East model, BABP and their noisy
versions, state the main theorem with its applications to these models, and
finally apply it to the positive-rates problem for simple spin systems.

\Needspace{9\baselineskip}
\subsection{Patch positivity}
\label{subsec:patch-positivity}

\begin{definition}[Patch positivity]
\label{def:patch-positivity}
A spin system is \emph{patch-positive} if, for every~$i\in\Lambda$ and every
nonempty~$S\subseteq N(i)$,
\begin{equation}
\label{eq:patch-positivity}
\widehat{\lambda_i}(S)\leq0,
\qquad
\widehat{\lambda_i\bb{\rho_i-\rho_i(\mb 0)}}(S)\geq0.
\end{equation}
\end{definition}

\begin{definition}[Threshold profile]
For a patch-positive spin system, define~$\mb p^\star=(p_i^\star)_{i\in\Lambda}$ by
\begin{equation}
\label{eq:threshold-profile}
p_i^\star
=
\max\bb{
\cb{0}
\cup
\cb{
\frac{\widehat{\lambda_i\rho_i}(S)}
{\widehat{\lambda_i}(S)}
:
\vn\neq S\subseteq N(i),\ 
\widehat{\lambda_i}(S)<0
}
}.
\end{equation}
We call~$\mb p^\star$ the \emph{threshold profile}.
\end{definition}

Patch positivity is a sign condition on the coefficients of the expansions of
$\lambda_i$ and~$\lambda_i\rho_i$ in the monomial basis. We check it for the
examples of Section~\ref{sec:setup}.

\begin{example}[FA-1f and East]
\label{ex:pp-fa-east}
For FA-1f, the only nonzero nonconstant coefficient of the update rate~$\lambda_i=1-\chi_{N(i)}$ is
\begin{equation*}
\widehat{\lambda_i}\bb{N(i)}=-1.
\end{equation*}
Since~$\rho_i=p$ is constant,
\begin{equation*}
\lambda_i\bb{\rho_i-\rho_i(\mb 0)}=0.
\end{equation*}
Thus FA-1f is patch-positive. Moreover,
\begin{equation*}
\frac{\widehat{\lambda_i\rho_i}\bb{N(i)}}
{\widehat{\lambda_i}\bb{N(i)}}
=p,
\end{equation*}
so its threshold profile is~$\mb p^\star=\mb p$.
The East model is a particular case of FA-1f, hence it is also patch-positive with the same threshold profile.
\end{example}

\begin{example}[FA-$j$f, $j\geq2$]
\label{ex:pp-fajf}
Patch positivity does not hold for all FA models. For FA-$j$f with~$j\geq2$,
the update rate is the indicator that at least~$j$ sites in~$N(i)$ are vacant.
Its expansion in the monomials~$\chi_S$ contains positive nonconstant
coefficients, so the first condition in~\eqref{eq:patch-positivity} fails.
\end{example}

\begin{example}[BABP]
\label{ex:pp-babp}
For BABP, the only nonzero nonconstant coefficients of
\begin{equation*}
\lambda_i(\eta)
=
\sum_{j\in N(i)}\bb{1-\eta(j)}
\end{equation*}
are
\begin{equation*}
\widehat{\lambda_i}\bb{\cb{j}}=-1,
\qquad j\in N(i).
\end{equation*}
Again~$\rho_i=p$ is constant, so
$\lambda_i\bb{\rho_i-\rho_i(\mb0)}=0$. Hence BABP is patch-positive.
Moreover, for every~$j\in N(i)$,
\begin{equation*}
\frac{\widehat{\lambda_i\rho_i}\bb{\cb{j}}}
{\widehat{\lambda_i}\bb{\cb{j}}}
=
p,
\end{equation*}
so~$p_i^\star=p$ if~$N(i)\neq\varnothing$ and~$p_i^\star=0$ if~$N(i)=\varnothing$.
Without noise, Sudbury~\cite{Sudbury1997} proved qualitative convergence of
BABP from translation-invariant initial laws, and on~$\Z^d$ every
homogeneous product law~$\mu_u$ with~$u<1$ converges exponentially
to~$\mu_p$ \cite[Application~2 and Section~5.1.1]{MartinelliShapiraToninelli2025}.
\end{example}

\begin{example}[FA-1f, East and BABP with noise]
\label{ex:pp-noisy}
Let~$\Lcal$ denote either FA-1f or BABP with refresh parameter~$p\in(0,1)$,
and assume that~$N(i)\neq\vn$ for every~$i\in\Lambda$. For
$\varepsilon>0$ and~$\mb s\in[0,1]^\Lambda$, consider
\begin{equation*}
\Lcal^{\varepsilon,\mb s}
=
\Lcal+\varepsilon\Ncal^{\mb s}.
\end{equation*}
Its update rate and refresh parameter are
\begin{equation*}
\lambda_i^{\varepsilon,\mb s}
=
\lambda_i+\varepsilon,
\qquad
\rho_i^{\varepsilon,\mb s}
=
\frac{p\lambda_i+\varepsilon s_i}
{\lambda_i+\varepsilon}.
\end{equation*}
For every nonempty~$S\subseteq N(i)$,
\begin{equation*}
\widehat{\lambda_i^{\varepsilon,\mb s}}(S)
=
\widehat{\lambda_i}(S),
\qquad
\widehat{\lambda_i^{\varepsilon,\mb s}
\rho_i^{\varepsilon,\mb s}}(S)
=
p\widehat{\lambda_i}(S),
\end{equation*}
and therefore
\begin{equation*}
\widehat{\lambda_i^{\varepsilon,\mb s}
\bb{\rho_i^{\varepsilon,\mb s}
-\rho_i^{\varepsilon,\mb s}(\mb0)}}(S)
=
\bb{p-\rho_i^{\varepsilon,\mb s}(\mb0)}
\widehat{\lambda_i}(S).
\end{equation*}
The nonconstant coefficients of~$\lambda_i$ are nonpositive in both models,
and at least one is strictly negative. Hence~$\Lcal^{\varepsilon,\mb s}$ is
patch-positive if and only if
\begin{equation}
\label{eq:perturbed-fa-babp-patch-positive}
\mb s\geq\mb p.
\end{equation}
Whenever this holds, its threshold profile remains~$\mb p^\star=\mb p$.
\end{example}

\begin{remark}
\label{rem:pure-deaths-not-patch-positive}
Condition~\eqref{eq:perturbed-fa-babp-patch-positive} excludes noise that
consists only of pure deaths,~$\mb s=\mb 0$. Thus FA-1f and East with added
pure deaths are not patch-positive for any~$p\in(0,1)$, and the
case~$s<p$ of Conjecture~\ref{conj:east-deaths} lies outside the scope of our
results. In the opposite
direction, removing pure deaths from a patch-positive system preserves patch
positivity and the threshold profile; see
Proposition~\ref{prop:death-removal}. We use this direction in the proof of Theorem~\ref{thm:common-invariant-limit}.
\end{remark}

The threshold profile determines the centering of the spins. For~$A\Subset\Lambda$, define the centered monomial
\begin{equation*}
\chi_A^\star(\eta)
=
\prod_{i\in A}\bb{\eta(i)-p_i^\star},
\qquad
\chi_\vn^\star=1,
\end{equation*}
and let
\begin{equation}
\label{eq:m-star-definition}
\mathcal M_\star
=
\cb{
\mu:
\mu\bb{\chi_A^\star}\geq0
\text{ for every }A\Subset\Lambda
}.
\end{equation}
Thus~$\mathcal M_\star$ consists of the laws under which, on every finite
set, the spins exceed the threshold profile jointly in the sense of
nonnegative centered moments. The Bernoulli product measure~$\mu_{\mb p^\star}$
belongs to~$\mathcal M_\star$, since all of its nonconstant centered
monomials have expectation zero, and so do~$\mu_{\mb 1}$ and every product
measure~$\mu_{\mb p}$ with~$\mb p\geq\mb p^\star$. For the convergence results
we use a larger class of initial measures.

\begin{definition}
\label{def:m-minus}
For~$K\geq0$, let~$\mathcal M_{-,K}$ be the set of probability measures~$\mu$
whose centered moments satisfy the lower bounds
\begin{equation}
\label{eq:m-minus-moment-bound}
\mu\bb{\chi_A^\star}
\geq
-K\prod_{i\in A}\bb{1-p_i^\star}
\qquad
\text{for every }A\Subset\Lambda,
\end{equation}
and set
\begin{equation*}
\mathcal M_-
=
\bigcup_{K\geq0}\mathcal M_{-,K}.
\end{equation*}
Since~$\mu_{\mb 1}\bb{\chi_A^\star}=\prod_{i\in A}\bb{1-p_i^\star}$,
condition~\eqref{eq:m-minus-moment-bound} is equivalent to
\begin{equation}
\label{eq:m-minus-affine-characterization}
\frac{\mu+K\mu_{\mb 1}}{1+K}\in\mathcal M_\star .
\end{equation}
\end{definition}

\begin{remark}
\label{rem:m-minus-examples}
We have~$\mathcal M_{-,0}=\mathcal M_\star$, while~$\mathcal M_-$ also
contains measures whose centered moments may be negative, but not more
negative than a fixed multiple of those of~$\mu_{\mb 1}$. If
$p_i^-:=\max\cb{2p_i^\star-1,0}$, then a Bernoulli product
measure~$\mu_{\mb p}$ belongs to~$\mathcal M_{-,1}$ exactly
when~$\mb p\geq\mb p^-$. In particular, if~$\mb p^\star\leq\mbs{\tfrac12}$,
then~$\abs{\chi_A^\star}\leq\prod_{i\in A}\bb{1-p_i^\star}$, so every
probability measure on~$\Omega$ belongs to~$\mathcal M_{-,1}$. Moreover,
if the system contains a pure-death component, then every point mass on a
configuration with finitely many~$0$s belongs to~$\mathcal M_-$; see
Section~\ref{subsec:proof-common-invariant-limit}.
\end{remark}

\subsection{Ergodicity of patch-positive spin systems}
\label{subsec:main-theorem}

We now state our main result. For patch-positive spin systems, a uniformly
positive rate of pure deaths causes convergence to a common invariant
measure from every initial law in~$\mathcal M_-$.

\begin{theorem}[Common invariant limit]
\label{thm:common-invariant-limit}
Suppose that the neighbourhood graph has polynomial growth with exponent~$D$,
the spin system is patch-positive, and it contains a pure-death component of
rate~$\varepsilon>0$. Then there is an invariant probability measure~$\pi$
such that, for every local function~$f$, there is~$K_f<\infty$ satisfying
\begin{equation}
\label{eq:uniform-m-minus-convergence}
\sup_{\mu\in\mathcal M_{-,K}}
\abs{(\mu P_t)(f)-\pi(f)}
\leq
(1+2K)K_f(1+t)^D e^{-\varepsilon t/2}
\end{equation}
for every~$K\geq0$ and~$t\geq0$. Consequently,
\begin{equation*}
\mu P_t\Rightarrow\pi
\qquad
\text{for every }\mu\in\mathcal M_-.
\end{equation*}
\end{theorem}

In Section~\ref{subsec:proof-common-invariant-limit} we prove
Theorem~\ref{thm:common-invariant-limit}. In the patch representation we
compare the system with the one we obtain by removing the pure deaths
(Proposition~\ref{prop:death-removal}). We show that dual interactions
after time~$T$ contribute only~$O(e^{-\varepsilon T})$, while the end patches
relax at rate~$\varepsilon$.

A direct consequence of Theorem~\ref{thm:common-invariant-limit} and the fact
that~$\mathcal M_{-,1}$ contains every probability measure
if~$\mb p^\star\leq\mbs{\tfrac12}$ (Remark~\ref{rem:m-minus-examples}) is the
following corollary.

\begin{corollary}[Uniform exponential ergodicity]
\label{cor:uniform-exponential-ergodicity}
Under the hypotheses of Theorem~\ref{thm:common-invariant-limit}, suppose also
that
\begin{equation*}
\mb p^\star\leq\mbs{\tfrac12}.
\end{equation*}
Then~$\pi$ is the unique invariant measure and, for every local function~$f$,
\begin{equation}
\label{eq:uniform-pointwise-convergence}
\sup_{\eta\in\Omega}
\abs{P_tf(\eta)-\pi(f)}
\leq
K_f(1+t)^D e^{-\varepsilon t/2}.
\end{equation}
\end{corollary}

\begin{corollary}[Uniform ergodicity under orientability]
\label{cor:oriented-ergodicity}
Under the hypotheses of Theorem~\ref{thm:common-invariant-limit}, suppose also
that the neighbourhood graph is orientable. Then~$\pi$ is the unique invariant
measure and, for every local function~$f$,
\begin{equation*}
\lim_{t\to\infty}
\sup_{\eta\in\Omega}
\abs{P_tf(\eta)-\pi(f)}
=
0.
\end{equation*}
\end{corollary}

In Section~\ref{subsec:proof-orientability} we prove
Corollary~\ref{cor:oriented-ergodicity} by a burn-in argument. We use
orientability to reduce the dynamics near a finite set to a finite
patch-positive system, which drives every evolution towards~$\mathcal M_\star$.
We can then apply Theorem~\ref{thm:common-invariant-limit}.

\begin{remark}
\label{rem:threshold-one-half}
In the proof of Corollary~\ref{cor:oriented-ergodicity} we use a burn-in
argument to show that every evolution approaches~$\mathcal M_\star$ on
finite sets,
\begin{equation*}
\liminf_{t\to\infty}\inf_{\eta\in\Omega}P_t\chi_A^\star(\eta)\geq0
\qquad\text{for every }A\Subset\Lambda,
\end{equation*}
and then we apply Theorem~\ref{thm:common-invariant-limit}. Orientability makes
this possible because, after conditioning on the dynamics on an autonomous
set outside a finite box around~$A$, the dynamics inside the box is a finite
patch-positive system with a time-dependent boundary condition. Such a system
preserves nonnegative centered moments and, thanks to the pure deaths, enters
$\mathcal M_\star$ with a uniformly positive probability in every unit of
time; see Section~\ref{subsec:proof-orientability}.
\end{remark}

\paragraph{Applications.}
\label{subsec:fa-east-applications}
Without noise, FA-1f and East have the explicit reversible measure~$\mu_p$, but
convergence to it is known only from restricted classes of initial laws; see
Conjecture~\ref{conj:fa1f}. Once independent refreshes with
density~$\mb s\neq\mb p$ are added, the invariant measure is no longer known
explicitly. We apply Theorem~\ref{thm:common-invariant-limit} and its corollaries
directly to the noisy models of Example~\ref{ex:pp-noisy}.

\begin{corollary}[FA-1f and BABP with noise]
\label{cor:perturbed-fa-babp-ergodicity}
\label{prop:perturbed-fa-babp-convergence}
Let~$\Lcal$ be FA-1f or BABP with refresh density~$p\in(0,1)$ on a
neighbourhood graph of polynomial growth with exponent~$D$, and assume
that~$N(i)\neq\vn$ for every~$i\in\Lambda$. Let~$\varepsilon>0$
and~$\mb s\in[0,1]^\Lambda$ satisfy
\begin{equation*}
\mb s\geq\mb p,
\qquad
\delta
:=
\varepsilon\inf_{i\in\Lambda}\bb{1-s_i}
>0,
\end{equation*}
and let~$P_t^{\varepsilon,\mb s}$ be the semigroup of
$\Lcal^{\varepsilon,\mb s}=\Lcal+\varepsilon\Ncal^{\mb s}$.
\begin{enumerate}[label=\textup{(\roman*)}]
\item There is an invariant probability measure~$\pi$ such that, for every
local function~$f$, there is~$K_f<\infty$ satisfying
\begin{equation}
\label{eq:perturbed-fa-babp-convergence}
\sup_{\mu\in\mathcal M_{-,K}}
\abs{(\mu P_t^{\varepsilon,\mb s})(f)-\pi(f)}
\leq
(1+2K)K_f(1+t)^D e^{-\delta t/2}
\end{equation}
for every~$K\geq0$ and~$t\geq0$. In particular,
$\mu P_t^{\varepsilon,\mb s}\Rightarrow\pi$ for every homogeneous product
law~$\mu_u$ with~$u\geq\max\cb{2p-1,0}$ and for every point mass on a
configuration with finitely many~$0$s. If~$\mb s=\mb p$, then~$\pi=\mu_p$.

\item If~$p\leq\frac12$, then~$\pi$ is the unique invariant measure and,
for every local function~$f$,
\begin{equation*}
\sup_{\eta\in\Omega}
\abs{P_t^{\varepsilon,\mb s}f(\eta)-\pi(f)}
\leq
K_f(1+t)^D e^{-\delta t/2}.
\end{equation*}

\item If the neighbourhood graph is orientable, then~$\pi$ is the unique
invariant measure and, for every local function~$f$,
\begin{equation*}
\lim_{t\to\infty}
\sup_{\eta\in\Omega}
\abs{P_t^{\varepsilon,\mb s}f(\eta)-\pi(f)}
=
0.
\end{equation*}
In particular, since the East model on~$\Z^d$ is FA-1f on the orientable
neighbourhood graph~$N(i)=\cb{i+e_k:k\leq d}$
(Example~\ref{ex:east-orientable}), the East model with such noise is
uniformly ergodic for every~$p\in(0,1)$.
\end{enumerate}
\end{corollary}

\begin{proof}
By Example~\ref{ex:pp-noisy},~$\Lcal^{\varepsilon,\mb s}$ is patch-positive
with threshold profile~$\mb p$. The added refreshes contain~$1\to0$
transitions at site~$i$ at rate~$\varepsilon(1-s_i)\geq\delta$, so
$\Lcal^{\varepsilon,\mb s}$ contains a pure-death component of rate~$\delta$.
Part~\textup{(i)} therefore follows from
Theorem~\ref{thm:common-invariant-limit} and
Remark~\ref{rem:m-minus-examples}, part~\textup{(ii)} from
Corollary~\ref{cor:uniform-exponential-ergodicity}, and part~\textup{(iii)}
from Corollary~\ref{cor:oriented-ergodicity}. If~$\mb s=\mb p$, then~$\mu_p$ is
invariant both for~$\Lcal$ and for~$\Ncal^{\mb p}$, hence
for~$\Lcal^{\varepsilon,\mb p}$. Since~$\mu_p\in\mathcal M_\star$, convergence
from~$\mu_p$ gives~$\pi=\mu_p$.
\end{proof}

\subsection{Positive rates for simple spin systems}
\label{subsec:simple-positive-rates}

We now turn to the positive-rates problem for a class of one-dimensional spin
systems.

\begin{definition}[Simple spin systems]
\label{def:simple-spin-systems}
A simple spin system~\cite{GluchowskiMenz2025, GluchowskiMenz2026} is a translation-invariant spin system on~$\Z$ with
$N(i)=\cb{i+1}$. It is parametrized by four numbers
$p_{1\mid ab}\in[0,1]$, $a,b\in\cb{0,1}$. When the spins at~$i$ and~$i+1$
are~$a$ and~$b$, respectively, a rate-one update at~$i$ sets the spin to~$1$
with probability~$p_{1\mid ab}$, and to~$0$ with probability~$p_{0\mid ab} = 1 -p_{1\mid ab}$. Equivalently, its refresh representation is
\begin{equation*}
\lambda_i(\eta)
=
p_{1\mid0\eta(i+1)}+1-p_{1\mid1\eta(i+1)},
\qquad
\rho_i(\eta)
=
\frac{p_{1\mid0\eta(i+1)}}{\lambda_i(\eta)}
\end{equation*}
when~$\lambda_i(\eta)>0$; again, if~$\lambda_i(\eta)=0$ then~$\rho_i(\eta)$ is irrelevant.
\end{definition}

Both flips have
positive rates for either value of the right neighbour precisely when
\begin{equation}
\label{eq:simple-positive-rates}
p_{1\mid11},p_{1\mid10}<1,
\qquad
p_{1\mid01},p_{1\mid00}>0.
\end{equation}
The broad positive-rates conjecture for one-dimensional finite-range systems
is false in general: G\'acs constructed counterexamples
\cite{Gacs1986,Gacs2001}; see also Gray's guide~\cite{Gray2001}. Classical
positive-rate results for monotone nearest-neighbour spin systems include
\cite{Gray1982}. The binary, one-sided nearest-neighbour problem considered
here remains open.

\begin{conjecture}[Positive rates for simple spin systems]
\label{conj:simple-positive-rates}
Every simple spin system satisfying~\eqref{eq:simple-positive-rates} is ergodic.
\end{conjecture}

Earlier sufficient conditions cover parts of the parameter space
\cite{GluchowskiMenz2025,GluchowskiMenz2026}. We use the orientability result to cover an additional region. The remaining region consists of noisy East
models that are not patch-positive; see
Remark~\ref{rem:simple-remaining-region}.

\begin{corollary}[Positive rates for simple spin systems]
\label{cor:simple-positive-rates}
Suppose that a simple spin system satisfies~\eqref{eq:simple-positive-rates}
and
\begin{equation}
\label{eq:simple-positive-rates-criterion}
p_{1\mid01}(1-p_{1\mid10})
\geq
p_{1\mid00}(1-p_{1\mid11}).
\end{equation}
Then the system has a unique invariant measure~$\pi$ and, for every local
function~$f$,
\begin{equation*}
\lim_{t\to\infty}
\sup_{\eta\in\Omega}
\abs{P_tf(\eta)-\pi(f)}
=
0.
\end{equation*}
\end{corollary}

\begin{proof}
In the refresh representation,
\begin{equation*}
\lambda_i(\eta)
=
p_{1\mid0\eta(i+1)}+1-p_{1\mid1\eta(i+1)},
\qquad
\rho_i(\eta)
=
\frac{p_{1\mid0\eta(i+1)}}{\lambda_i(\eta)}.
\end{equation*}
By translation invariance, we write~$\lambda(b)=\lambda_i(\eta^{i+1,b})$ for~$b\in\cb{0,1}$.
Exchanging the roles of~$0$ and~$1$ swaps the two values of~$\lambda$ and
leaves~\eqref{eq:simple-positive-rates-criterion} unchanged. We may therefore
assume that~$\lambda(1)\leq\lambda(0)$. Since
$N(i)=\cb{i+1}$,
\begin{align*}
\widehat{\lambda_i}\bb{\cb{i+1}}
&=
p_{1\mid01}-p_{1\mid11}-p_{1\mid00}+p_{1\mid10}\\
&= \lambda(1)-\lambda(0)
\leq0,\\
\widehat{\lambda_i\bb{\rho_i-\rho_i(\mb0)}}\bb{\cb{i+1}}
&=
\frac{
p_{1\mid01}(1-p_{1\mid10})
-
p_{1\mid00}(1-p_{1\mid11})
}{
p_{1\mid00}+1-p_{1\mid10}
}
\geq0.
\end{align*}
Thus the system is patch-positive. Its~$1\to0$ rates are bounded below by
\begin{equation*}
\min\cb{1-p_{1\mid11},1-p_{1\mid10}}>0,
\end{equation*}
so it contains a pure-death component. The one-sided neighbourhood graph has
polynomial growth and is orientable. We can therefore apply Corollary~\ref{cor:oriented-ergodicity}.
\end{proof}

\begin{figure}[tbp]
\centering
\begingroup
\color{black}
\begin{tikzpicture}[scale=4, line cap=round, line join=round]
  \pgfmathsetmacro{\criticalz}{3-2*sqrt(2)}
  \foreach \z/\xshift/\yshift/\denom in
    {0.01/0/0/100,0.125/1.55/0/8,0.25/0/-1.55/4} {
    \begin{scope}[xshift=\xshift cm,yshift=\yshift cm]
      \pgfmathsetmacro{\oldcross}{(sqrt(2)+\z)/(1+sqrt(2))}
      \pgfmathsetmacro{\redend}{1-\z/sqrt(2)}

      \fill[lightgray!25] (0,0) -- (1,0) -- (1,\z) -- (0.5,\z)
        -- (0.5,{0.5-\z}) -- (1,{1-\z})
        -- (1,1) -- (0,1) -- cycle;
      \fill[darkgray!25] (0.5,\z) -- (\redend,\z)
        -- (\oldcross,{\oldcross-\z})
        -- (0.5,{0.5-\z}) -- cycle;

      \ifdim\z pt<\criticalz pt
        \pgfmathsetmacro{\newcross}{1-sqrt(\z/sqrt(2))}
        \pgfmathsetmacro{\newend}{(1+\z+sqrt(1-6*\z+\z*\z))/2}
        \fill[red!25] (\oldcross,{\oldcross-\z})
          -- (\newend,{\newend-\z})
          -- plot[domain=\newend:\newcross,samples=100,variable=\x]
            ({\x},{\z/(1-\x)}) -- cycle;
      \fi

      \draw[thick] (0.5,\z) -- (0.5,{0.5-\z}) -- (1,{1-\z});
      \draw[thick] (0.5,\z) -- (1,\z);
      \draw[thick,dashed,darkgray!70!black]
        (\oldcross,{\oldcross-\z}) -- (\redend,\z);
      \ifdim\z pt<\criticalz pt
        \draw[thick,red!70!black]
          plot[domain=\newcross:\newend,samples=100,variable=\x]
            ({\x},{\z/(1-\x)});
      \fi

      \draw (0,0) rectangle (1,1);
      \draw[->] (-0.02,0) -- (1.08,0)
        node[below right] {$p_{1\mid10}$};
      \draw[->] (0,-0.02) -- (0,1.08)
        node[above left] {$p_{1\mid01}$};
      \foreach \t in {0,1} {
        \draw (\t,0) -- (\t,-0.02) node[below] {\small $\t$};
        \draw (0,\t) -- (-0.02,\t) node[left] {\small $\t$};
      }
      \draw (0.5,0) -- (0.5,-0.02) node[below] {\small $\tfrac12$};
      \draw (0,0.5) -- (-0.02,0.5) node[left] {\small $\tfrac12$};
      \node at (0.5,-0.20) {$p_{1\mid00}=\tfrac1{\denom}$};
    \end{scope}
  }
  \begin{scope}[shift={(1.43,-0.75)},every node/.style={anchor=west,font=\small}]
    \fill[lightgray!25] (0,0) rectangle (0.065,0.055);
    \node at (0.09,0.0275) {Corollaries 7.1 and 7.2 in~\cite{GluchowskiMenz2025}};
    \fill[darkgray!25] (0,-0.15) rectangle (0.065,-0.095);
    \node at (0.09,-0.1225) {Long-lived-state criterion~\cite{GluchowskiMenz2026}};
    \fill[red!25] (0,-0.30) rectangle (0.065,-0.245);
    \node at (0.09,-0.2725) {Corollary~\ref{cor:simple-positive-rates}};
    \draw (0,-0.45) rectangle (0.065,-0.395);
    \node at (0.09,-0.4225) {Not proven};
  \end{scope}
\end{tikzpicture}
\endgroup
\caption{Sections of the cube~$\cb{p_{1\mid11}=0}$. The horizontal coordinate is
$p_{1\mid10}$ and the vertical coordinate is~$p_{1\mid01}$; the three panels
have~$p_{1\mid00}=1/100$,~$1/8$, and~$1/4$. Light gray is the region covered
by~\cite{GluchowskiMenz2025}, where
$p_{1\mid10}<\tfrac12$, or
$p_{1\mid10}<p_{1\mid01}+p_{1\mid00}$, or
$p_{1\mid01}\leq p_{1\mid00}$. Dark gray is the additional region from~\cite{GluchowskiMenz2026}, where
$\max\cb{p_{1\mid01},p_{1\mid00}}<\sqrt2(1-p_{1\mid10})$, and red is the
additional region given by Corollary~\ref{cor:simple-positive-rates}, where
$p_{1\mid00}\leq(1-p_{1\mid10})p_{1\mid01}$. White denotes parameters not
covered by these criteria. The improvement is greatest when spontaneous births are sparse, meaning~$p_{1\mid 00}\approx 0$.}
\label{fig:simple-ips-ergodicity-regions}
\end{figure}
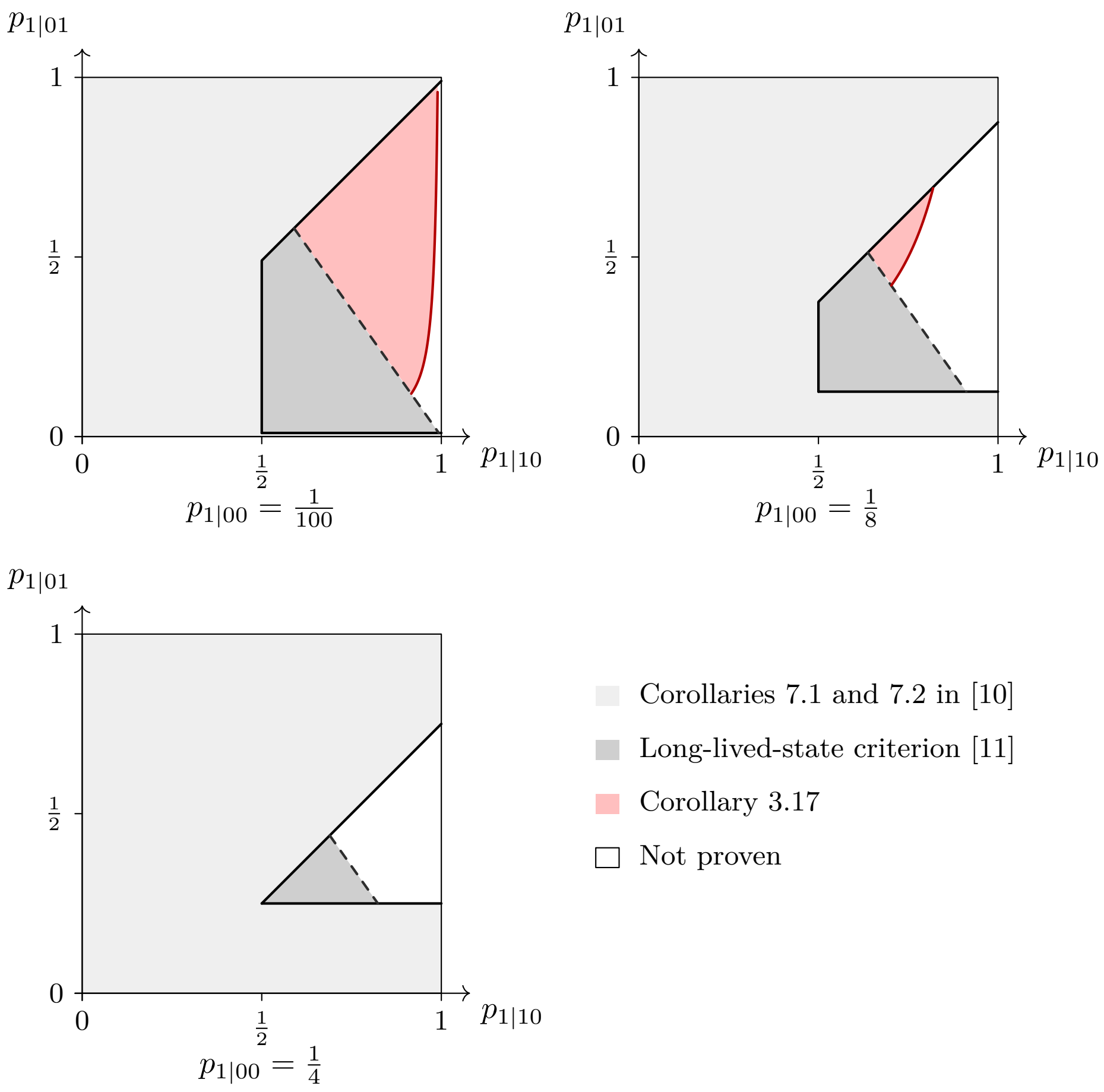

To compare this criterion with the earlier sufficient conditions, we use the
same normalization as in~\cite{GluchowskiMenz2025}. Symmetry of exchanging the two states and a time-scaling trick~\cite{GluchowskiMenz2025} reduce the problem to two cubes in the parameter space
\begin{equation*}
\cb{p_{1\mid11}=0, p_{1\mid 10}<1, p_{1\mid 01}>0, p_{1\mid 00}>0}
\quad\text{or}\quad
\cb{p_{1\mid11}<1, p_{1\mid 10}=0, p_{1\mid 01}>0, p_{1\mid 00}>0}.
\end{equation*}
The spin systems with parameters in the second cube are already known to be ergodic
\cite[Section~7]{GluchowskiMenz2025}, so we focus on~$p_{1\mid11}=0$. On this
face,~\eqref{eq:simple-positive-rates-criterion} becomes
\begin{equation}
\label{eq:simple-positive-rates-face}
p_{1\mid00}
\leq
(1-p_{1\mid10})p_{1\mid01}.
\end{equation}
In Figure~\ref{fig:simple-ips-ergodicity-regions} we compare this condition
with the earlier sufficient conditions on three sections of the parameter
space.

\begin{remark}[The remaining region]
\label{rem:simple-remaining-region}
In Figure~\ref{fig:simple-ips-ergodicity-regions} the parameters we leave
uncovered correspond to noisy East models with~$\mb s<\mb p$. Indeed, on the face~$p_{1\mid11}=0$ these parameters
satisfy~$p_{1\mid01}>p_{1\mid00}$. After exchanging the states~$0$ and~$1$
and changing time by the factor
$c=p_{1\mid10}+p_{1\mid01}-p_{1\mid00}>0$, such a system becomes the East
model on~$\Z$ with refresh density~$p$ and added noise~$\varepsilon\Ncal^{\mb s}$
with constant~$\mb s\equiv s$, where
\begin{equation*}
p=\frac{p_{1\mid10}}{c},
\qquad
\varepsilon=\frac{1-p_{1\mid10}+p_{1\mid00}}{c},
\qquad
s=\frac{1-p_{1\mid10}}{1-p_{1\mid10}+p_{1\mid00}},
\end{equation*}
and~$p,s\in(0,1)$,~$\varepsilon>0$. In these coordinates
\begin{equation*}
p_{1\mid01}\bb{1-p_{1\mid10}}-p_{1\mid00}\bb{1-p_{1\mid11}}
=
c^2\varepsilon\,(s-p),
\end{equation*}
so criterion~\eqref{eq:simple-positive-rates-criterion} is exactly the
patch-positivity condition~$s\geq p$ of Example~\ref{ex:pp-noisy}. Consequently,
removing the restriction~$\mb s\geq\mb p$ from
Corollary~\ref{cor:perturbed-fa-babp-ergodicity}\textup{(iii)} for the East
model, that is, proving uniform ergodicity of the East model with added
noise~$\varepsilon\Ncal^{\mb s}$ for all~$\varepsilon>0$ and~$s\in(0,1)$, would,
together with the reductions and results
of~\cite{GluchowskiMenz2025,GluchowskiMenz2026}, prove
Conjecture~\ref{conj:simple-positive-rates}; this is the case~$s\in(0,p)$ of
Conjecture~\ref{conj:east-deaths}. Its endpoint~$s=0$, pure deaths, lies
outside the positive-rates region, since there~\eqref{eq:simple-positive-rates}
fails.
\end{remark}
\FloatBarrier
\section{The patch representation}
\label{sec:probabilistic-interpretation}
\label{subsec:proofs-section-4}

We now construct the representation that we use to prove the main
convergence theorem. We begin with the classic dual process~$\bb{\sigma,A}$, with~$\sigma_t\in\cb{+,-}$,~$A_t\Subset \Lambda$, which traces dependence of the spin system on the initial configuration. 
By duality we obtain the representation
\begin{equation*}
P_t\chi_A(\eta) = \E\Cb{\sigma_t \exp\bb{\int_0^t V(A_s)ds}\chi_{A_t}(\eta)},
\end{equation*}
for a potential~$V$. For an ergodicity argument one typically needs to argue that the event $\cb{A_t\neq\vn}$ has negligible contribution, which boils down to showing~$e^{\int V}$ is small or capturing cancellations from~$\sigma$.

Ours is the second route. We aim to average with respect to a~$\sigma$-algebra that is small enough to yield sufficient averaging but rich enough for the calculation to be tractable. Our choice is a~$\sigma$-algebra that reveals when and where the interactions between sites in the dual occur, and nothing else. In particular, we leave out the types of those interactions. These interactions divide the timeline at each site into time intervals, which we call patches. Conditional on this~$\sigma$-algebra the remaining randomness factorizes over patches, with laws simple enough to permit explicit computation.

We obtain the patch representation and can now explain the conditions we
introduced in Section~\ref{sec:main-results}. Patch positivity makes the
contributions to~$P_t\chi_A$ of patches that don't reach time~$t$ (bulk
patches) nonnegative, and the threshold profile determines which initial
laws give a nonnegative contribution from patches that cross~$t$ (end
patches). 

\subsection{Dual process}
\label{sec:dual}

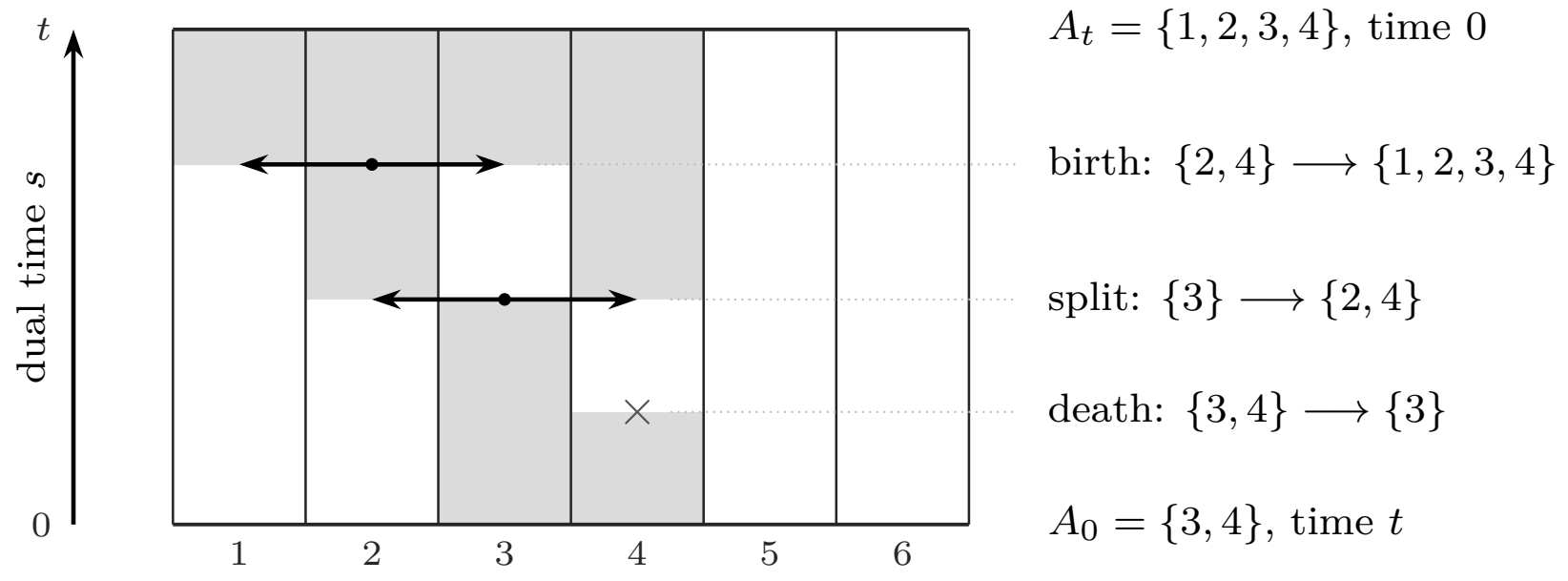
\begin{figure}[tbp]
\centering
\begingroup
\definecolor{dualfigactive}{gray}{0.86}
\definecolor{dualfigused}{gray}{0.32}
\definecolor{dualfigoutline}{gray}{0.15}
\begin{tikzpicture}[x=0.92cm,y=0.78cm]
\tikzset{
  dualfig site/.style={draw=dualfigoutline,line width=0.48pt},
  dualfig active/.style={fill=dualfigactive,draw=none},
  dualfig interaction/.style={draw=black,line width=0.86pt,
    -{Stealth[length=2.0mm,width=1.35mm]}},
  dualfig horizon/.style={draw=dualfigoutline,line width=0.72pt},
  dualfig guide/.style={draw=gray!50,line width=0.40pt,densely dotted},
  dualfig label/.style={font=\footnotesize,text=black,inner sep=1pt},
  dualfig small/.style={font=\scriptsize,text=dualfigoutline,inner sep=1pt},
  dualfig death/.style={font=\normalsize,text=dualfigused,inner sep=0pt}
}
\path[dualfig active] (2,0) rectangle (3,2);
\path[dualfig active] (3,0) rectangle (4,1);
\path[dualfig active] (3,2) rectangle (4,4.4);
\path[dualfig active] (1,2) rectangle (2,4.4);
\path[dualfig active] (0,3.2) rectangle (1,4.4);
\path[dualfig active] (2,3.2) rectangle (3,4.4);
\foreach \x in {0,...,6}{\draw[dualfig site] (\x,0)--(\x,4.4);}
\draw[dualfig horizon] (0,0)--(6,0);
\draw[dualfig horizon] (0,4.4)--(6,4.4);
\foreach \x/\n in {0.5/1,1.5/2,2.5/3,3.5/4,4.5/5,5.5/6}{
  \node[dualfig small,anchor=north] at (\x,-0.10) {$\n$};
}
\node[dualfig death] at (3.5,1) {$\times$};
\draw[dualfig guide] (3.75,1)--(6.35,1);
\node[dualfig label,anchor=west] at (6.55,1)
  {death: $\cb{3,4}\longrightarrow\cb{3}$};
\draw[dualfig interaction] (2.5,2)--(1.5,2);
\draw[dualfig interaction] (2.5,2)--(3.5,2);
\fill[black] (2.5,2) circle[radius=1.25pt];
\draw[dualfig guide] (3.75,2)--(6.35,2);
\node[dualfig label,anchor=west] at (6.55,2)
  {split: $\cb{3}\longrightarrow\cb{2,4}$};
\draw[dualfig interaction] (1.5,3.2)--(0.5,3.2);
\draw[dualfig interaction] (1.5,3.2)--(2.5,3.2);
\fill[black] (1.5,3.2) circle[radius=1.25pt];
\draw[dualfig guide] (2.75,3.2)--(6.35,3.2);
\node[dualfig label,anchor=west] at (6.55,3.2)
  {birth: $\cb{2,4}\longrightarrow\cb{1,2,3,4}$};
\draw[dualfig interaction] (-0.75,0)--(-0.75,4.4);
\node[dualfig label,rotate=90,anchor=south] at (-0.92,2.2)
  {dual time~$s$};
\node[dualfig small,anchor=east] at (-0.87,0) {$0$};
\node[dualfig small,anchor=east] at (-0.87,4.4) {$t$};
\node[dualfig label,anchor=west] at (6.55,0)
  {$A_0=\cb{3,4}$, time~$t$};
\node[dualfig label,anchor=west] at (6.55,4.4)
  {$A_t=\cb{1,2,3,4}$, time~$0$};
\end{tikzpicture}
\endgroup
\caption{A trajectory of the dual set, with dual time increasing upward
(backward in the original spin-system time). Gray shading marks the sites
in~$A_s$. Starting from~$A_0=\cb{3,4}$, a death removes site~$4$.
A split replaces site~$3$ by~$\cb{2,4}$, and a birth adds sites~$1$ and~$3$
while keeping sites~$2$ and~$4$. Thus~$A_t=\cb{1,2,3,4}$. We omit the sign coordinate~$\sigma_s$.}
\label{fig:dual-idea}
\end{figure}

The dual traces the dependence on the initial configuration backward in
time: the dual set~$A_s$ consists of the sites whose spins at time~$t-s$ the
value~$\chi_A(\eta_t)$ may still depend on. Roughly speaking, when a site
of~$A_s$ is updated, it is removed if its new spin is drawn independently of
all other spins (a \emph{death}); it is replaced by a set~$S$ of neighbours
if the new spin is determined by the spins in~$S$ (a \emph{split}); and~$S$
is added if the outcome depends both on its earlier spin and on~$S$, for
example because the update may be blocked (a \emph{birth}). If the dual set becomes empty,~$\chi_A(\eta_t)$ no longer depends on the
initial configuration. In Figure~\ref{fig:dual-idea} we illustrate the dual
process. Such set-valued duals go back
to Harris~\cite{Harris1978} and
Griffeath~\cite[Chapter~II, Section~1 and Chapter~III, Section~1]{Griffeath1979}
for additive and cancellative systems, see also
\cite[Chapter~III, Sections~4 and~6]{Liggett1985}
and~\cite{Gray1986,JansenKurt2014,SturmSwart2018},
and~\cite{LubetzkySly2016} for a related backward exploration. In general the dual transitions carry signs, which we record in the sign
process~$\sigma_s\in\cb{+,-}$. These come from the inclusion-exclusion in
the expansion of the rates. 

We use a monomial duality function. For~$A\Subset\Lambda$ and~$\sigma\in\cb{+,-}$  define
\begin{equation}
\label{eq:monomial-duality-function}
H((A,\sigma),\eta)
=
\sigma\chi_A(\eta).
\end{equation}
We seek a signed set-valued Markov process~$(A_t,\sigma_t)$ and a potential~$V$
such that
\begin{equation}
\label{eq:monomial-generator-duality}
\Lcal_\eta H((A,\sigma),\eta)
=
\Dcal H((A,\sigma),\eta)
+
V(A)H((A,\sigma),\eta),
\end{equation}
where~$\Dcal$ is the generator of~$(A_t,\sigma_t)$. We read off the transitions and the potential from the expansion
of~$\Lcal\chi_A$ in the monomial basis,
\begin{equation}
\label{eq:signed-dual-generator-expansion}
\Lcal\chi_A
=
\sum_{i\in A}\sum_{S\subseteq N(i)}
a_i^\delta(S)\chi_{\bb{A\setminus\cb{i}}\cup S}
+
\sum_{i\in A}\sum_{S\subseteq N(i)}
a_i^\beta(S)\chi_{A\cup S},
\end{equation}
where
\begin{equation}
\label{eq:signed-dual-coefficients}
a_i^\delta(S)
=
\widehat{\lambda_i\rho_i}(S),
\qquad
a_i^\beta(S)
=
-\widehat{\lambda_i}(S).
\end{equation}
(see the proof of Proposition~\ref{thm:monomial-duality}). In the second
sum, the term with~$S=\varnothing$ leaves the index unchanged; it is a diagonal
term with coefficient~$a_i^\beta(\varnothing)=-\widehat{\lambda_i}(\varnothing)
=-\lambda_i(\mb 0)$. The absolute values of the remaining coefficients give the
rates of the marked updates, and their signs are carried by the second coordinate.

\begin{definition}[Signed dual]
\label{def:signed-dual}
With~$\operatorname{sgn}_{\pm}(r)=+$ for~$r\geq0$
and~$\operatorname{sgn}_{\pm}(r)=-$ for~$r<0$, set, for~$S\subseteq N(i)$,
\begin{equation}
\label{eq:death-split-data}
\delta_i(S)
=
\abs{a_i^\delta(S)},
\qquad
\sigma_i^\delta(S)
=
\operatorname{sgn}_{\pm}\bb{a_i^\delta(S)},
\end{equation}
and, for~$\varnothing\neq S\subseteq N(i)$,
\begin{equation}
\label{eq:birth-data}
\beta_i(S)
=
\abs{a_i^\beta(S)},
\qquad
\sigma_i^\beta(S)
=
\operatorname{sgn}_{\pm}\bb{a_i^\beta(S)}.
\end{equation}
For~$i\in A$, define
\begin{align*}
D_{i,S}(A,\sigma)
&=
\bb{
\bb{A\setminus\cb{i}}\cup S,
\sigma\sigma_i^\delta(S)
},
\\
B_{i,S}(A,\sigma)
&=
\bb{
A\cup S,
\sigma\sigma_i^\beta(S)
}.
\end{align*}
The \emph{signed dual} is the Markov process~$(A_t,\sigma_t)$ with generator
\begin{align}
\label{eq:signed-dual-generator}
\Dcal g(A,\sigma)
={}&
\sum_{i\in A}\sum_{S\subseteq N(i)}
\delta_i(S)
\bb{
g\bb{D_{i,S}(A,\sigma)}-g(A,\sigma)
}
\notag\\
&+
\sum_{i\in A}
\sum_{\substack{S\subseteq N(i)\\S\neq\varnothing}}
\beta_i(S)
\bb{
g\bb{B_{i,S}(A,\sigma)}-g(A,\sigma)
},
\end{align}
and \emph{potential}
\begin{equation}
\label{eq:dual-potential}
V(A)
=
\sum_{i\in A}V_i,
\qquad
V_i
=
\sum_{S\subseteq N(i)}\delta_i(S)
+
\sum_{\substack{S\subseteq N(i)\\S\neq\varnothing}}\beta_i(S)
+
a_i^\beta(\varnothing).
\end{equation}
We write~$\pP_A$ and~$\E_A$ for the law and expectation of the signed dual
started from~$(A,+)$.
\end{definition}

\begin{remark}
The map~$D_{i,\varnothing}$ removes~$i$ (death),~$D_{i,S}$
with~$S\neq\varnothing$ replaces~$i$ by~$S$ (split), and~$B_{i,S}$ adds~$S$
(birth). By the uniform bounds on~$\abs{N(i)}$ and the rates, the dual jump
rate per active site and~$\abs{V_i}$ are uniformly bounded.
\end{remark}

The duality holds for infinite spin systems given an integrability condition, but it always holds in finite volume. Thus, to make our results not depend on a technical condition we introduce the zero-boundary finite approximation of the spin system. For~$R\Subset\Lambda$
and~$\eta\in\Omega$, let~$\eta^{R,0}$ agree with~$\eta$ on~$R$ and vanish
outside~$R$.

\begin{definition}[Zero-boundary system]
\label{def:finite-volume}
Let~$R\Subset\Lambda$. The \emph{zero-boundary system in~$R$} is the
Markov chain on~$\cb{0,1}^R$ with generator
\begin{equation}
\label{eq:zero-boundary-generator}
\Lcal^{R,0}f(\eta)
=
\sum_{i\in R}
\lambda_i\bb{\eta^{R,0}}
\bb{E_{i,\rho_i(\eta^{R,0})}-I}f(\eta),
\end{equation}
acting on functions~$f$ that depend only on the spins in~$R$. We write~$P_t^{R,0}$
for its semigroup. For the signed dual started from~$A\subseteq R$, the
\emph{exit time} from~$R$ and the \emph{confinement event} are
\begin{equation}
\label{eq:confinement-event}
\tau_R=\inf\cb{s\geq0:A_s\nsubseteq R},
\qquad
E_t^R=\cb{\tau_R>t}.
\end{equation}
\end{definition}

The finite-volume zero-boundary system approximates the infinite-volume spin system in a sense made precise by the following lemma. We prove it in Appendix~\ref{app:passage-to-infinite-volume}.

\begin{lemma}[Finite propagation]
\label{lem:finite-propagation}
Let~$f$ be supported on~$A\Subset\Lambda$.
\begin{enumerate}[label=\textup{(\alph*)}]
\item For every~$a>0$, there are~$v<\infty$ and~$C_A<\infty$ such that
\begin{equation}
\label{eq:finite-propagation}
\norm{\bb{P_t-P_t^{R,0}}f}_\infty
\leq
C_Ae^{-at}\norm f_\infty
\end{equation}
for every~$t\geq0$ and every finite~$R$ containing~$B\bb{A,\lceil vt\rceil}$.
\item For every~$t\geq0$,~$\norm{\bb{P_t-P_t^{R,0}}f}_\infty\to0$
as~$R\uparrow\Lambda$ along finite sets.
\end{enumerate}
\end{lemma}

\begin{proposition}[Monomial Feynman--Kac duality in finite volume]
\label{thm:monomial-duality}
For every $R\Subset\Lambda$, $A\subseteq R$, $t\geq0$ and $\eta\in\Omega$,
\begin{equation}
\label{eq:monomial-feynman-kac-duality}
P_t^{R,0}\chi_A(\eta)
=
\E_A\Cb{
\sigma_t
\exp\bb{\int_0^t V(A_s)\,ds}
\chi_{A_t}(\eta)
\ind\bb{E_t^R}
}.
\end{equation}
\end{proposition}

\begin{proof}[Proof of Proposition~\ref{thm:monomial-duality}]
If~$i\notin A$, then~$\bb{E_{i,\rho_i(\eta)}-I}\chi_A=0$. If~$i\in A$, then
\begin{equation*}
\lambda_i\bb{E_{i,\rho_i}-I}\chi_A
=
\lambda_i\rho_i\chi_{A\setminus\cb{i}}
-
\lambda_i\chi_A.
\end{equation*}
Expanding~$\lambda_i\rho_i$ and~$\lambda_i$ in monomials and
using~$\chi_S\chi_B=\chi_{S\cup B}$ gives~\eqref{eq:signed-dual-generator-expansion}
with the coefficients~\eqref{eq:signed-dual-coefficients}.
By~\eqref{eq:death-split-data} and~\eqref{eq:birth-data}, each dual rate
multiplied by its sign is the corresponding coefficient. In~$\Dcal H+V(A)H$,
the rate terms subtracted by~$\Dcal$ cancel those in~$V(A)$, leaving the marked
update terms and~$\sum_{i\in A}a_i^\beta(\varnothing)H$.
This gives the generator identity~\eqref{eq:monomial-generator-duality}.

Now fix~$R\Subset\Lambda$. For~$S\Subset\Lambda$ we
have~$\chi_S(\eta^{R,0})=\chi_S(\eta)\ind\cb{S\subseteq R}$. Hence the same
computation with the rates evaluated at~$\eta^{R,0}$ gives, for~$B\subseteq R$,
\begin{equation}
\label{eq:confined-spin-generator-on-H}
\Lcal^{R,0}\chi_B
=
\sum_{i\in B}
\sum_{S\subseteq N(i)\cap R}
a_i^\delta(S)
\chi_{\bb{B\setminus\cb{i}}\cup S}
+
\sum_{i\in B}
\sum_{S\subseteq N(i)\cap R}
a_i^\beta(S)
\chi_{B\cup S}.
\end{equation}
Let~$\Dcal_R$ be the generator of the signed set process on the finite
state space
\begin{equation*}
\cb{(B,\sigma):B\subseteq R,\ \sigma\in\cb{+,-}}
\end{equation*}
obtained from~$\Dcal$ by removing the jumps with target~$S\nsubseteq R$, and let
\begin{equation*}
\kappa_R(B)
=
\sum_{i\in B}
\sum_{\substack{\varnothing\neq S\subseteq N(i)\\S\nsubseteq R}}
\bb{\delta_i(S)+\beta_i(S)}
\end{equation*}
be the total rate of the removed jumps. As above, we compare~\eqref{eq:confined-spin-generator-on-H}
with~\eqref{eq:signed-dual-generator} and~\eqref{eq:dual-potential} to obtain
\begin{equation}
\label{eq:confined-generator-duality}
\Lcal_\eta^{R,0}H((B,\sigma),\eta)
=
\Dcal_RH((B,\sigma),\eta)
+
\bb{V(B)-\kappa_R(B)}H((B,\sigma),\eta),
\qquad B\subseteq R.
\end{equation}
Both generators are bounded operators on finite spaces, so this identity
integrates to
\begin{equation*}
P_t^{R,0}H((A,+),\cdot)(\eta)
=
e^{t\bb{\Dcal_R+V-\kappa_R}}H(\cdot,\eta)(A,+),
\end{equation*}
see for example~\cite{JansenKurt2014}. By the Feynman--Kac formula for the
finite-state chain with generator~$\Dcal_R$, the right-hand side is the
expectation of~$\exp\bb{\int_0^t(V-\kappa_R)}\,H$ along that chain. Up to
time~$\tau_R$, this chain and the signed dual make the same jumps at the same
rates; the dual additionally makes the removed jumps, at total
rate~$\kappa_R$, and the first of them occurs exactly at~$\tau_R$. Hence the
weight~$\exp\bb{-\int_0^t\kappa_R}$ is the probability that the dual makes no
removed jump along a given path of the chain, and the right-hand side
equals
\begin{equation*}
\E_A\Cb{\sigma_t\exp\bb{\int_0^tV(A_s)\,ds}\chi_{A_t}(\eta)\ind\bb{E_t^R}}.
\end{equation*}
This proves~\eqref{eq:monomial-feynman-kac-duality}.
\end{proof}

\subsection{Patches}
\label{sec:patches}

Without the sign~$\sigma_t$, the duality~\eqref{eq:monomial-feynman-kac-duality}
would express~$P_t^{R,0}\chi_A$ through nonnegative weights, and inequalities
between the potentials of two systems would carry over directly to the
semigroups. To control the cancellations of the signs, we introduce
a~$\sigma$-algebra that reveals when sites interact, but not how. The interactions divide each site line into intervals, which we call
patches. In Section~\ref{subsec:patch-representation} we show that,
conditionally on the interactions, the law of the dual process factorizes
over the patches.
In Figure~\ref{fig:patch-geometry} we show a realization of the dual
(left) and the patches it induces, split at a finite time~$t$ into bulk
patches and end patches (right). We first construct the dual
graphically.

\begin{definition}[Graphical construction]
\label{def:graphical-construction}
For
every~$i\in\Lambda$ and~$S\subseteq N(i)$, let
$I_{i,S}^\delta\subseteq(0,\infty)$ be a Poisson point process of
rate~$\delta_i(S)$. For every nonempty~$S\subseteq N(i)$, let
$I_{i,S}^\beta\subseteq(0,\infty)$ be a Poisson point process of
rate~$\beta_i(S)$. All these processes are independent, and we set
$I_{i,\varnothing}^\beta=\varnothing$. We denote the resulting collection of
Poisson marks by
\begin{equation*}
I^{\mathrm P}
=
\cb{
(i,t,\alpha,S):
i\in\Lambda,\ 
\alpha\in\cb{\delta,\beta},\ 
S\subseteq N(i),\ 
t\in I_{i,S}^\alpha
}.
\end{equation*}
For a mark~$(i,t,\alpha,S)$, the site~$i$ is its \emph{source},~$S$ its
\emph{target}, and~$\alpha$ its \emph{kind}. Starting from~$(A,+)$, a
mark~$(i,t,\delta,S)$ with~$i\in A_{t-}$ applies~$D_{i,S}$, while a
mark~$(i,t,\beta,S)$ with~$i\in A_{t-}$ applies~$B_{i,S}$. A mark whose
source is not active does not change the process.
\end{definition}
This gives the signed dual process with generator~$\Dcal$
from~\eqref{eq:signed-dual-generator}. The marks through which dependence passes between distinct sites are exactly
those with nonempty target that occur at active sources.
\begin{definition}[Successful interaction and its record]
\label{def:successful-interaction}
A mark~$(i,t,\alpha,S)\in I^{\mathrm P}$ is a \emph{successful
interaction} if
\begin{equation*}
S\neq\varnothing
\qquad\text{and}\qquad
i\in A_{t-}.
\end{equation*}
Its \emph{record} is the triple~$(i,t,S)$.
Together with the formal initial record~$(\infty,0,A)$, with
source~$\infty\notin\Lambda$, the \emph{successful-interaction skeleton} up
to time~$T<\infty$ is
\begin{equation}
\label{eq:successful-interaction-skeleton}
\begin{aligned}
\mathcal I_T
&=
\cb{(\infty,0,A)}
\\
&\quad\cup
\cb{
(i,t,S):
t\leq T
\text{ and }(i,t,\alpha,S)\text{ is successful for some }
\alpha\in\cb{\delta,\beta}
}.
\end{aligned}
\end{equation}
We also set~$\mathcal I
=
\bigcup_{T<\infty}\mathcal I_T$.
\end{definition}

\begin{remark}\label{rem:interaction-records}
The set of records~$\mathcal I_T$ gives the times, sources and targets of
successful interactions up to time~$T$. It omits deaths, marks at inactive
sources and the kinds of the successful interactions.
\end{remark}

\begin{lemma}[Local finiteness]
\label{lem:dual-local-finiteness}
For every~$T<\infty$, the successful-interaction skeleton~$\mathcal I_T$ is
almost surely finite.
\end{lemma}

\begin{proof}[Proof of Lemma~\ref{lem:dual-local-finiteness}]
The total rate of marks whose source is active is bounded by~$C\abs{A_t}$ for some finite~$C$, uniformly in the state. Moreover, a single
mark can add at most
\begin{equation*}
\sup_{i\in\Lambda}\abs{N(i)}<\infty
\end{equation*}
sites to the active set. Hence~$\abs{A_t}$ is dominated by a nonexplosive
continuous-time branching process. Almost surely only finitely many marks occur
at active sources on~$[0,T]$, and therefore only finitely many of them are
successful interactions.
\end{proof}

Each record~$(i,t,S)\in\mathcal I$ is outgoing at~$i$ and incoming at every
site of~$S$. On each site line these times cut time into consecutive
intervals, the patches.

\begin{definition}[Patch]
\label{def:patch}
For~$i\in\Lambda$, the times of incoming and outgoing successful interactions
at~$i$ are
\begin{align*}
\mathcal T^{\mathsf I}(i)
&=
\cb{t\geq0:(j,t,S)\in\mathcal I
  \text{ for some }j,S\text{ with }i\in S},
\\
\mathcal T^{\mathsf O}(i)
&=
\cb{t>0:(i,t,S)\in\mathcal I\text{ for some }S},
\\
\mathcal T(i)
&=
\mathcal T^{\mathsf I}(i)\cup\mathcal T^{\mathsf O}(i).
\end{align*}
By the formal record,~$0\in\mathcal T^{\mathsf I}(i)$ for~$i\in A$: the
initial set counts as an incoming interaction. For~$s\in\mathcal T(i)$, the
next successful interaction at~$i$ occurs at
\begin{equation*}
e_i(s)
=
\inf\cb{u>s:u\in\mathcal T(i)},
\qquad
\inf\varnothing=\infty.
\end{equation*}
Let~$(j,s,S)\in\mathcal I$ be a record for which~$i=j$ or~$i\in S$. The patch
beginning at~$(i,s)$ is
\begin{equation*}
P
=
\bb{i(P),[s(P),e(P)),\mathsf X(P)\mathsf Y(P),S(P)},
\end{equation*}
where
\begin{equation*}
i(P)=i,
\qquad
s(P)=s,
\qquad
e(P)=e_i(s),
\qquad
S(P)=S.
\end{equation*}
The initial label~$\mathsf X(P)$ records whether the interaction at~$s(P)$ is
incoming or outgoing:
\begin{equation*}
\mathsf X(P)
=
\pw{\mathsf I}{i(P)\in S(P)}{\mathsf O}{i(P)\notin S(P)}.
\end{equation*}
The terminal label~$\mathsf Y(P)$ records the next successful interaction on
the same site line, with~$\mathsf E$ meaning that there is no such interaction:
\begin{equation*}
\mathsf Y(P)
=
\pwthree
{\mathsf I}{e(P)<\infty\text{ and }e(P)\in\mathcal T^{\mathsf I}\bb{i(P)}}
{\mathsf O}{e(P)<\infty\text{ and }e(P)\in\mathcal T^{\mathsf O}\bb{i(P)}}
{\mathsf E}{e(P)=\infty}.
\end{equation*}
Let~$\mathcal P$ be the family of all patches.
\end{definition}

\begin{definition}[Bulk and end patches]
\label{def:bulk-end-patches}
Fix a finite time~$t\geq0$. For~$P\in\mathcal P$ with~$s(P)\leq t$, the patch
\emph{cut at time~$t$} is
\begin{equation*}
\begin{aligned}
P^{\downarrow t}
&=
\bb{
 i(P),
 [s(P),\min\cb{e(P),t}),
 \mathsf X(P)\mathsf Y^{\downarrow t}(P),
 S(P)
},
\\
\mathsf Y^{\downarrow t}(P)
&=
\pw{\mathsf Y(P)}{e(P)\leq t}{\mathsf E}{t<e(P)}.
\end{aligned}
\end{equation*}
The families of \emph{bulk patches}, \emph{end patches}, and
\emph{finite-horizon patches} are
\begin{align*}
\mathcal B_t
&=
\cb{P\in\mathcal P:e(P)\leq t},\\
\mathcal E_t
&=
\cb{P^{\downarrow t}:P\in\mathcal P,\ s(P)\leq t<e(P)},\\
\mathcal P_t
&=
\mathcal B_t\cup\mathcal E_t.
\end{align*}
For~$P\in\mathcal E_t$ and finite~$u\geq t$, the \emph{extension} of~$P$
through time~$u$ is
\begin{equation}
\label{eq:geometric-patch-extension}
P^{\uparrow u}
=
\bb{i(P),[s(P),u),\mathsf X(P)\mathsf E,S(P)}.
\end{equation}
\end{definition}

\begin{remark}
In the infinite-horizon family~$\mathcal P$, a patch with terminal
label~$\mathsf E$ is infinite, while a patch with terminal label~$\mathsf I$
or~$\mathsf O$ is finite. Bulk patches end by time~$t$; end patches cross~$t$ and are cut there, and
$P^{\uparrow t}=P$ for~$P\in\mathcal E_t$. By
Lemma~\ref{lem:dual-local-finiteness},~$\mathcal P_t$ is almost surely
finite, and distinct end patches have distinct sites.
In Figure~\ref{fig:patch-geometry} we show a realization and its patch
family. The thick red line separates~$\mathcal B_t$ from~$\mathcal E_t$.
\end{remark}
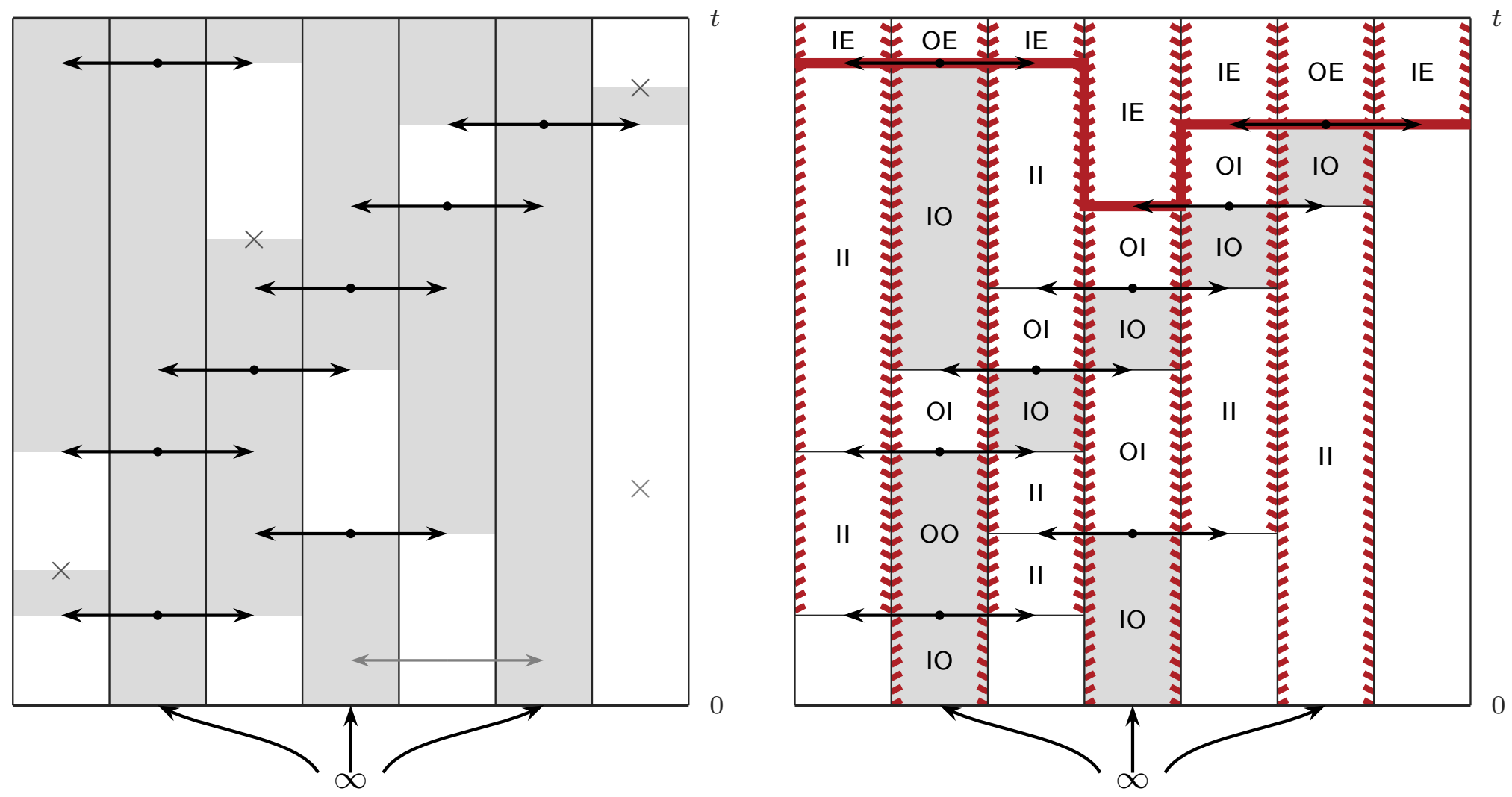
\begin{figure}[tbp]
\centering
\resizebox{\linewidth}{!}{\begingroup
\definecolor{patchfigred}{RGB}{175,32,38}
\definecolor{patchfigactive}{gray}{0.86}
\definecolor{patchfigused}{gray}{0.32}
\definecolor{patchfigunused}{gray}{0.50}
\definecolor{patchfigoutline}{gray}{0.15}
\begin{tikzpicture}[x=0.92cm,y=0.78cm]
\tikzset{
  patchfig site/.style={draw=patchfigoutline,line width=0.48pt},
  patchfig active/.style={fill=patchfigactive,draw=none},
  patchfig success/.style={draw=black,line width=0.86pt,-{Stealth[length=2.0mm,width=1.35mm]}},
  patchfig unused/.style={draw=patchfigunused,line width=0.68pt,-{Stealth[length=1.75mm,width=1.18mm]}},
  patchfig horizon/.style={draw=patchfigoutline,line width=0.72pt},
  patchfig boundary/.style={draw=patchfigoutline,line width=0.40pt},
  patchfig label/.style={font=\scriptsize,text=black,inner sep=0pt},
  patchfig separator/.style={draw=patchfigred,line width=2.70pt,line cap=butt,line join=miter},
  patchfig small/.style={font=\scriptsize,text=patchfigoutline},
  patchfig used death/.style={font=\normalsize,text=patchfigused,inner sep=0pt},
  patchfig unused death/.style={font=\normalsize,text=patchfigunused,inner sep=0pt}
}
\newcommand{\PatchFigGrid}{%
  \foreach \x in {0,...,7}{\draw[patchfig site](\x,0)--(\x,8.4);}
  \draw[patchfig horizon](0,0)--(7,0);
  \draw[patchfig horizon](0,8.4)--(7,8.4);}
\newcommand{\PatchFigBarrier}[3]{%
  \begin{scope}
    \clip (#1,#2) rectangle ({#1+1},#3);
    \foreach \k in {0,...,54}{%
      \pgfmathsetmacro{\yy}{0.165*\k}%
      \draw[draw=patchfigred,line width=1.65pt,line cap=butt]
        (#1,\yy) -- ({#1+0.105},{\yy+0.064});
      \draw[draw=patchfigred,line width=1.65pt,line cap=butt]
        ({#1+0.895},{\yy+0.064}) -- ({#1+1},\yy);
    }%
  \end{scope}%
}
\newcommand{\PatchFigInitial}{%
  \node[font=\small,inner sep=0pt](patchfiginf) at (3.5,-0.92){$\infty$};
  \draw[patchfig success]([xshift=-1.4mm]patchfiginf.north west)..controls(2.95,-0.40)and(2.05,-0.31)..(1.5,0);
  \draw[patchfig success](patchfiginf.north)--(3.5,0);
  \draw[patchfig success]([xshift=1.4mm]patchfiginf.north east)..controls(4.05,-0.40)and(4.95,-0.31)..(5.5,0);}
\begin{scope}
  \path[patchfig active] (0,1.1) rectangle (1,1.65);
  \path[patchfig active] (0,3.1) rectangle (1,8.4);
  \path[patchfig active] (1,0) rectangle (2,8.4);
  \path[patchfig active] (2,1.1) rectangle (3,5.7);
  \path[patchfig active] (2,7.85) rectangle (3,8.4);
  \path[patchfig active] (3,0) rectangle (4,2.1);
  \path[patchfig active] (3,4.1) rectangle (4,8.4);
  \path[patchfig active] (4,2.1) rectangle (5,6.1);
  \path[patchfig active] (4,7.1) rectangle (5,8.4);
  \path[patchfig active] (5,0) rectangle (6,8.4);
  \path[patchfig active] (6,7.1) rectangle (7,7.55);
  \PatchFigGrid
  \draw[patchfig success] (1.5,1.1)--(0.5,1.1);
  \draw[patchfig success] (1.5,1.1)--(2.5,1.1);
  \fill[black] (1.5,1.1) circle[radius=1.25pt];
  \draw[patchfig success] (3.5,2.1)--(2.5,2.1);
  \draw[patchfig success] (3.5,2.1)--(4.5,2.1);
  \fill[black] (3.5,2.1) circle[radius=1.25pt];
  \draw[patchfig success] (1.5,3.1)--(0.5,3.1);
  \draw[patchfig success] (1.5,3.1)--(2.5,3.1);
  \fill[black] (1.5,3.1) circle[radius=1.25pt];
  \draw[patchfig success] (2.5,4.1)--(1.5,4.1);
  \draw[patchfig success] (2.5,4.1)--(3.5,4.1);
  \fill[black] (2.5,4.1) circle[radius=1.25pt];
  \draw[patchfig success] (3.5,5.1)--(2.5,5.1);
  \draw[patchfig success] (3.5,5.1)--(4.5,5.1);
  \fill[black] (3.5,5.1) circle[radius=1.25pt];
  \draw[patchfig success] (4.5,6.1)--(3.5,6.1);
  \draw[patchfig success] (4.5,6.1)--(5.5,6.1);
  \fill[black] (4.5,6.1) circle[radius=1.25pt];
  \draw[patchfig success] (5.5,7.1)--(4.5,7.1);
  \draw[patchfig success] (5.5,7.1)--(6.5,7.1);
  \fill[black] (5.5,7.1) circle[radius=1.25pt];
  \draw[patchfig success] (1.5,7.85)--(0.5,7.85);
  \draw[patchfig success] (1.5,7.85)--(2.5,7.85);
  \fill[black] (1.5,7.85) circle[radius=1.25pt];
  \draw[patchfig unused] (4.5,0.55)--(3.5,0.55);
  \draw[patchfig unused] (4.5,0.55)--(5.5,0.55);
  \node[patchfig used death] at (0.5,1.65) {$\times$};
  \node[patchfig used death] at (2.5,5.7) {$\times$};
  \node[patchfig used death] at (6.5,7.55) {$\times$};
  \node[patchfig unused death] at (6.5,2.65) {$\times$};
  \PatchFigInitial
  \node[patchfig small,anchor=west] at (7.08,8.4){$t$};
  \node[patchfig small,anchor=west] at (7.08,0){$0$};
\end{scope}
\begin{scope}[shift={(8.10,0)}]
  \path[patchfig active] (1,0) rectangle (2,1.1);
  \path[patchfig active] (1,1.1) rectangle (2,3.1);
  \path[patchfig active] (1,4.1) rectangle (2,7.85);
  \path[patchfig active] (2,3.1) rectangle (3,4.1);
  \path[patchfig active] (3,0) rectangle (4,2.1);
  \path[patchfig active] (3,4.1) rectangle (4,5.1);
  \path[patchfig active] (4,5.1) rectangle (5,6.1);
  \path[patchfig active] (5,6.1) rectangle (6,7.1);
  \PatchFigBarrier{0}{1.1}{8.4}
  \PatchFigBarrier{1}{0}{8.4}
  \PatchFigBarrier{2}{1.1}{8.4}
  \PatchFigBarrier{3}{0}{8.4}
  \PatchFigBarrier{4}{2.1}{8.4}
  \PatchFigBarrier{5}{0}{8.4}
  \PatchFigBarrier{6}{7.1}{8.4}
  \PatchFigGrid
  \draw[patchfig boundary] (0,1.1)--(1,1.1);
  \draw[patchfig boundary] (0,3.1)--(1,3.1);
  \draw[patchfig boundary] (0,7.85)--(1,7.85);
  \draw[patchfig boundary] (0,8.4)--(1,8.4);
  \draw[patchfig boundary] (1,0)--(2,0);
  \draw[patchfig boundary] (1,1.1)--(2,1.1);
  \draw[patchfig boundary] (1,3.1)--(2,3.1);
  \draw[patchfig boundary] (1,4.1)--(2,4.1);
  \draw[patchfig boundary] (1,7.85)--(2,7.85);
  \draw[patchfig boundary] (1,8.4)--(2,8.4);
  \draw[patchfig boundary] (2,1.1)--(3,1.1);
  \draw[patchfig boundary] (2,2.1)--(3,2.1);
  \draw[patchfig boundary] (2,3.1)--(3,3.1);
  \draw[patchfig boundary] (2,4.1)--(3,4.1);
  \draw[patchfig boundary] (2,5.1)--(3,5.1);
  \draw[patchfig boundary] (2,7.85)--(3,7.85);
  \draw[patchfig boundary] (2,8.4)--(3,8.4);
  \draw[patchfig boundary] (3,0)--(4,0);
  \draw[patchfig boundary] (3,2.1)--(4,2.1);
  \draw[patchfig boundary] (3,4.1)--(4,4.1);
  \draw[patchfig boundary] (3,5.1)--(4,5.1);
  \draw[patchfig boundary] (3,6.1)--(4,6.1);
  \draw[patchfig boundary] (3,8.4)--(4,8.4);
  \draw[patchfig boundary] (4,2.1)--(5,2.1);
  \draw[patchfig boundary] (4,5.1)--(5,5.1);
  \draw[patchfig boundary] (4,6.1)--(5,6.1);
  \draw[patchfig boundary] (4,7.1)--(5,7.1);
  \draw[patchfig boundary] (4,8.4)--(5,8.4);
  \draw[patchfig boundary] (5,0)--(6,0);
  \draw[patchfig boundary] (5,6.1)--(6,6.1);
  \draw[patchfig boundary] (5,7.1)--(6,7.1);
  \draw[patchfig boundary] (5,8.4)--(6,8.4);
  \draw[patchfig boundary] (6,7.1)--(7,7.1);
  \draw[patchfig boundary] (6,8.4)--(7,8.4);
  \node[patchfig label] at (0.5,2.1) {$\mathsf{II}$};
  \node[patchfig label] at (0.5,5.475) {$\mathsf{II}$};
  \node[patchfig label] at (0.5,8.125) {$\mathsf{IE}$};
  \node[patchfig label] at (1.5,0.55) {$\mathsf{IO}$};
  \node[patchfig label] at (1.5,2.1) {$\mathsf{OO}$};
  \node[patchfig label] at (1.5,3.6) {$\mathsf{OI}$};
  \node[patchfig label] at (1.5,5.975) {$\mathsf{IO}$};
  \node[patchfig label] at (1.5,8.125) {$\mathsf{OE}$};
  \node[patchfig label] at (2.5,1.6) {$\mathsf{II}$};
  \node[patchfig label] at (2.5,2.6) {$\mathsf{II}$};
  \node[patchfig label] at (2.5,3.6) {$\mathsf{IO}$};
  \node[patchfig label] at (2.5,4.6) {$\mathsf{OI}$};
  \node[patchfig label] at (2.5,6.475) {$\mathsf{II}$};
  \node[patchfig label] at (2.5,8.125) {$\mathsf{IE}$};
  \node[patchfig label] at (3.5,1.05) {$\mathsf{IO}$};
  \node[patchfig label] at (3.5,3.1) {$\mathsf{OI}$};
  \node[patchfig label] at (3.5,4.6) {$\mathsf{IO}$};
  \node[patchfig label] at (3.5,5.6) {$\mathsf{OI}$};
  \node[patchfig label] at (3.5,7.25) {$\mathsf{IE}$};
  \node[patchfig label] at (4.5,3.6) {$\mathsf{II}$};
  \node[patchfig label] at (4.5,5.6) {$\mathsf{IO}$};
  \node[patchfig label] at (4.5,6.6) {$\mathsf{OI}$};
  \node[patchfig label] at (4.5,7.75) {$\mathsf{IE}$};
  \node[patchfig label] at (5.5,3.05) {$\mathsf{II}$};
  \node[patchfig label] at (5.5,6.6) {$\mathsf{IO}$};
  \node[patchfig label] at (5.5,7.75) {$\mathsf{OE}$};
  \node[patchfig label] at (6.5,7.75) {$\mathsf{IE}$};
  \draw[patchfig separator] (0,7.85) -- (1,7.85) -- (2,7.85) -- (3,7.85) -- (3,6.1) -- (4,6.1) -- (4,7.1) -- (5,7.1) -- (6,7.1) -- (7,7.1);
  \draw[patchfig success] (1.5,1.1)--(0.5,1.1);
  \draw[patchfig success] (1.5,1.1)--(2.5,1.1);
  \fill[black] (1.5,1.1) circle[radius=1.25pt];
  \draw[patchfig success] (3.5,2.1)--(2.5,2.1);
  \draw[patchfig success] (3.5,2.1)--(4.5,2.1);
  \fill[black] (3.5,2.1) circle[radius=1.25pt];
  \draw[patchfig success] (1.5,3.1)--(0.5,3.1);
  \draw[patchfig success] (1.5,3.1)--(2.5,3.1);
  \fill[black] (1.5,3.1) circle[radius=1.25pt];
  \draw[patchfig success] (2.5,4.1)--(1.5,4.1);
  \draw[patchfig success] (2.5,4.1)--(3.5,4.1);
  \fill[black] (2.5,4.1) circle[radius=1.25pt];
  \draw[patchfig success] (3.5,5.1)--(2.5,5.1);
  \draw[patchfig success] (3.5,5.1)--(4.5,5.1);
  \fill[black] (3.5,5.1) circle[radius=1.25pt];
  \draw[patchfig success] (4.5,6.1)--(3.5,6.1);
  \draw[patchfig success] (4.5,6.1)--(5.5,6.1);
  \fill[black] (4.5,6.1) circle[radius=1.25pt];
  \draw[patchfig success] (5.5,7.1)--(4.5,7.1);
  \draw[patchfig success] (5.5,7.1)--(6.5,7.1);
  \fill[black] (5.5,7.1) circle[radius=1.25pt];
  \draw[patchfig success] (1.5,7.85)--(0.5,7.85);
  \draw[patchfig success] (1.5,7.85)--(2.5,7.85);
  \fill[black] (1.5,7.85) circle[radius=1.25pt];
  \PatchFigInitial
  \node[patchfig small,anchor=west] at (7.08,8.4){$t$};
  \node[patchfig small,anchor=west] at (7.08,0){$0$};
\end{scope}
\end{tikzpicture}
\endgroup}
\caption{A realization of the graphical construction (left) and the induced
finite-horizon patch family at time~$t$ (right). On the left, gray shading marks
times when a site is dual-active, black arrows are successful interactions,
gray arrows are rings at inactive sources, and crosses are deaths. On the right,
gray shading highlights patches with outgoing terminal boundary, which remain
active throughout. Red striped sides delimit patches, the labels give their
initial and terminal boundary types, and the thick red line separates bulk
patches from end patches.}
\label{fig:patch-geometry}
\end{figure}
\FloatBarrier

\subsection{Factorization into patch laws}
\label{subsec:patch-laws}

 Now we can define the~$\sigma$-algebra that facilitates the averaging we need to control the cancellations from the sign process. For~$T<\infty$, define
\begin{equation*}
\mathcal G_T:=\sigma(\mathcal I_T),
\qquad
\mathcal G_\infty:=\sigma(\mathcal I).
\end{equation*}
Conditioning on~$\mathcal{G}_T$ fixes the interaction times, sources and targets, and hence the patch family~$\mathcal P_T$. The kinds of the interactions and the marks inside patches remain random, see Remark~\ref{rem:interaction-records}.

Conditional on a realization of~$\mathcal{I}_T$ the events inside each patch must be consistent with~$\mathcal{I}_T$: no successful interaction can originate inside a patch, and patches of type~$\mathsf{IO}$ or~$\mathsf{OO}$ must be active throughout to facilitate the terminal outgoing interaction. We first describe a law for the randomness in one patch and the consistency conditions. We show that, conditional on~$\mathcal G_T$, the law of~$\cb{\bb{\sigma_t,A_t}}_{t\leq T}$ factorizes over patches, with each patch law obtained by imposing its own consistency conditions.

\begin{definition}[Patch law]
\label{def:patch-law}
Fix a patch~$P$, write~$i=i(P)$,~$s=s(P)$,~$e=e(P)$,~$S=S(P)$, and let
\begin{equation*}
\Sigma_P^\circ
=
\cb{(i,u,\gamma,R)\in I^{\mathrm P}:s<u<e}
\end{equation*}
be the marks from~$i$ inside the patch. The~$(\delta,R)$ marks have rate~$\delta_i(R)$,
and for nonempty~$R$ the~$(\beta,R)$ marks have rate~$\beta_i(R)$. These Poisson
processes are independent. If~$\mathsf X(P)=\mathsf I$, the initial
boundary is incoming and we set~$\Sigma_P=\Sigma_P^\circ$.
If~$\mathsf X(P)=\mathsf O$, the boundary record specifies the target~$S$ but
not whether the interaction is a split or a birth; let
$\alpha(P)\in\cb{\delta,\beta}$ denote its kind, independent
of~$\Sigma_P^\circ$, with
\begin{equation}
\label{eq:patch-initial-kind-law}
\pP_P\bb{\alpha(P)=\delta}
=
\frac{\delta_i(S)}{\delta_i(S)+\beta_i(S)},
\qquad
\pP_P\bb{\alpha(P)=\beta}
=
\frac{\beta_i(S)}{\delta_i(S)+\beta_i(S)},
\end{equation}
and set
\begin{equation*}
\Sigma_P
=
\Sigma_P^\circ\cup\cb{(i,s,\alpha(P),S)}.
\end{equation*}
The \emph{patch law}~$\pP_P$ is the law of~$\Sigma_P$, and~$\E_P$ denotes
expectation under~$\pP_P$.
The marks in~$\Sigma_P$ determine a one-site process~$X^P\in\cb{0,1}$, where
$X_u^P=1$ means that the patch site is active at time~$u$. Its initial value is
\begin{equation*}
X_s^P
=
\pw{1}{\mathsf X(P)=\mathsf I}
{\ind\cb{\alpha(P)=\beta}}{\mathsf X(P)=\mathsf O}.
\end{equation*}
We then read the interior marks in chronological order. When~$X_{u-}^P=1$, a
$\delta$-mark sets~$X_u^P=0$ and a~$\beta$-mark leaves~$X_u^P=1$; when
$X_{u-}^P=0$, either mark leaves~$X_u^P=0$. Between marks,~$X^P$ is constant.
For a finite patch with terminal label~$\mathsf E$ and~$e>s$, set~$X_e^P=X_{e-}^P$.
\end{definition}

Almost surely~$\delta_i(S)+\beta_i(S)>0$ in~\eqref{eq:patch-initial-kind-law},
since the interaction starting the patch has positive rate.

\begin{definition}[Consistent patch law]
\label{def:consistent-patch-law}
The \emph{consistency event} of~$P$ is
\begin{equation*}
\operatorname{Con}(P)
=
\cb{
\begin{aligned}
X_{u-}^P=0
&\quad\text{for every }(i,u,\gamma,R)\in\Sigma_P^\circ\\[-0.2em]
&\quad\text{with }R\neq\varnothing
\end{aligned}
}
\cap
\pw{\cb{X_{e-}^P=1}}{\mathsf Y(P)=\mathsf O}
{\Omega_P}{\mathsf Y(P)\in\cb{\mathsf I,\mathsf E}},
\end{equation*}
where~$\Omega_P$ is the sample space of the patch law.
If~$\pP_P\bb{\operatorname{Con}(P)}>0$, the \emph{consistent patch law} is
\begin{equation*}
\pP_P^{\mathrm{con}}(\cdot)
=
\pP_P\bb{\,\cdot\,\middle|\,\operatorname{Con}(P)},
\end{equation*}
and~$\E_P^{\mathrm{con}}$ denotes expectation under this law.
\end{definition}

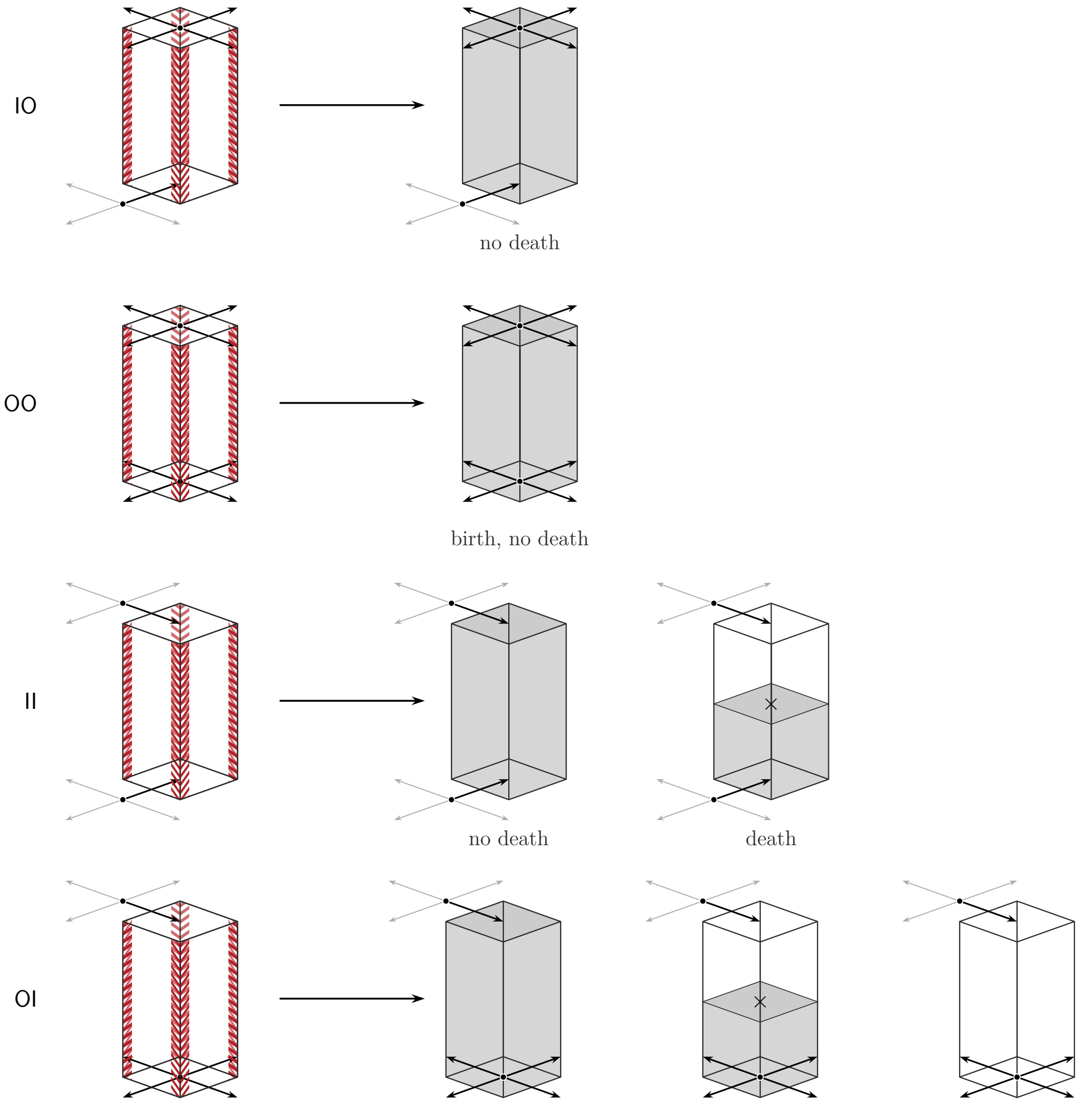
\begin{figure}[tbp]
\centering
\begingroup
\definecolor{patchred}{RGB}{180,34,40}
\definecolor{activefront}{gray}{0.84}
\definecolor{activeback}{gray}{0.89}
\definecolor{activetop}{gray}{0.80}
\definecolor{secondary}{gray}{0.68}
\definecolor{outline}{gray}{0.17}
\tikzset{
  patch outline/.style={draw=outline,line width=0.62pt,line join=round},
  slice outline/.style={draw=outline,line width=0.43pt},
  active front/.style={fill=activefront,draw=none},
  active back/.style={fill=activeback,draw=none},
  active top/.style={fill=activetop,draw=none},
  inactive face/.style={fill=white,draw=none},
  inactive side/.style={fill=white,draw=none},
  patch stripe/.style={fill=patchred,draw=none},
  patch stripe back/.style={fill=patchred,fill opacity=0.66,draw=none},
  boundary interaction/.style={draw=black,line width=0.86pt,-{Stealth[length=1.8mm,width=1.25mm]}},
  external secondary/.style={draw=secondary,line width=0.52pt,-{Stealth[length=1.35mm,width=0.95mm]}},
  used death/.style={font=\large,text=black,inner sep=0pt},
  patch type label/.style={font=\large,inner sep=0pt},
  patch case label/.style={font=\fontsize{11}{12}\selectfont,text=outline,align=center,inner sep=0pt},
  patch decomposition arrow/.style={draw=black,line width=0.95pt,-{Stealth[length=2.4mm,width=1.55mm]}},
}

\newcommand{\PatchCoordinates}{%
  \coordinate (W) at (-1.18,0);
  \coordinate (N) at (0,0.42);
  \coordinate (E) at (1.18,0);
  \coordinate (S) at (0,-0.42);
  \coordinate (Wt) at (-1.18,3.2);
  \coordinate (Nt) at (0,3.62);
  \coordinate (Et) at (1.18,3.2);
  \coordinate (St) at (0,2.78);
}

\newcommand{\PatchOutline}{%
  \draw[patch outline] (W)--(N)--(E)--(S)--cycle;
  \draw[patch outline] (Wt)--(Nt)--(Et)--(St)--cycle;
  \draw[patch outline] (W)--(Wt) (N)--(Nt) (E)--(Et) (S)--(St);
}

\newcommand{\PatchStripeFamily}[3]{%
  \begin{scope}
    \coordinate (root) at (#1);
    \coordinate (inner) at ($(#1)!0.16!(#2)$);
    \clip (root)--(inner)--($(inner)+(0,3.2)$)--($(root)+(0,3.2)$)--cycle;
    \foreach \zz in {-0.4,-0.235,...,3.56}{%
      \pgfmathsetmacro{\zza}{\zz+0.075}
      \pgfmathsetmacro{\zzb}{\zz+0.2003}
      \pgfmathsetmacro{\zzc}{\zz+0.2753}
      \path[#3] ($(root)+(0,\zz)$)--($(root)+(0,\zza)$)--($(inner)+(0,\zzc)$)--($(inner)+(0,\zzb)$)--cycle;
    }
  \end{scope}
  \begin{scope}
    \coordinate (root) at (#2);
    \coordinate (inner) at ($(#2)!0.16!(#1)$);
    \clip (root)--(inner)--($(inner)+(0,3.2)$)--($(root)+(0,3.2)$)--cycle;
    \foreach \zz in {-0.4,-0.235,...,3.56}{%
      \pgfmathsetmacro{\zza}{\zz+0.075}
      \pgfmathsetmacro{\zzb}{\zz+0.2003}
      \pgfmathsetmacro{\zzc}{\zz+0.2753}
      \path[#3] ($(root)+(0,\zz)$)--($(root)+(0,\zza)$)--($(inner)+(0,\zzc)$)--($(inner)+(0,\zzb)$)--cycle;
    }
  \end{scope}
}

\newcommand{\RearPatchBarriers}{%
  \PatchStripeFamily{W}{N}{patch stripe back}%
  \PatchStripeFamily{N}{E}{patch stripe back}%
}
\newcommand{\FrontPatchBarriers}{%
  \PatchStripeFamily{E}{S}{patch stripe}%
  \PatchStripeFamily{S}{W}{patch stripe}%
}

\newcommand{\SourceDot}[2]{%
  \fill[white] (#1,#2) circle[radius=2.15pt];
  \fill[black] (#1,#2) circle[radius=1.55pt];
}
\newcommand{\OutgoingAt}[1]{%
  \draw[boundary interaction] (0,#1)--(1.18,{#1+0.42});
  \draw[boundary interaction] (0,#1)--(-1.18,{#1-0.42});
  \draw[boundary interaction] (0,#1)--(1.18,{#1-0.42});
  \draw[boundary interaction] (0,#1)--(-1.18,{#1+0.42});
  \SourceDot{0}{#1}%
}
\newcommand{\IncomingInitial}{%
  \draw[boundary interaction] (-1.18,-0.42)--(0,0);
  \draw[external secondary] (-1.18,-0.42)--(-2.36,-0.84);
  \draw[external secondary] (-1.18,-0.42)--(0,-0.84);
  \draw[external secondary] (-1.18,-0.42)--(-2.36,0);
  \SourceDot{-1.18}{-0.42}%
}
\newcommand{\IncomingTerminal}{%
  \draw[boundary interaction] (-1.18,3.62)--(0,3.2);
  \draw[external secondary] (-1.18,3.62)--(-2.36,4.04);
  \draw[external secondary] (-1.18,3.62)--(0,4.04);
  \draw[external secondary] (-1.18,3.62)--(-2.36,3.2);
  \SourceDot{-1.18}{3.62}%
}

\newcommand{\PatchBlock}[2]{%
  \PatchCoordinates
  \RearPatchBarriers
  \PatchOutline
  #1#2
  \FrontPatchBarriers
  \draw[patch outline] (W)--(N)--(E)--(S)--cycle;
  \draw[patch outline] (Wt)--(Nt)--(Et)--(St)--cycle;
  \draw[patch outline] (E)--(Et) (S)--(St) (W)--(Wt);
}

\newcommand{\PrismBase}[2]{%
  \PatchCoordinates
  \path[#1] (W)--(N)--(E)--(S)--cycle;
  \path[inactive side] (W)--(N)--(Nt)--(Wt)--cycle;
  \path[inactive side] (N)--(E)--(Et)--(Nt)--cycle;
  \path[inactive side] (E)--(S)--(St)--(Et)--cycle;
  \path[inactive side] (S)--(W)--(Wt)--(St)--cycle;
  #2
}
\newcommand{\FullPrism}{%
  \PrismBase{active top}{%
    \path[active back] (W)--(N)--(Nt)--(Wt)--cycle;
    \path[active back] (N)--(E)--(Et)--(Nt)--cycle;
    \path[active front] (E)--(S)--(St)--(Et)--cycle;
    \path[active front] (S)--(W)--(Wt)--(St)--cycle;
    \path[active top] (Wt)--(Nt)--(Et)--(St)--cycle;
  }%
  \PatchOutline
}
\newcommand{\EmptyPrism}{%
  \PrismBase{inactive face}{}%
  \path[inactive face] (Wt)--(Nt)--(Et)--(St)--cycle;
  \PatchOutline
}
\newcommand{\DeathPrism}{%
  \PatchCoordinates
  \path[active top] (W)--(N)--(E)--(S)--cycle;
  \path[inactive side] (W)--(N)--(Nt)--(Wt)--cycle;
  \path[inactive side] (N)--(E)--(Et)--(Nt)--cycle;
  \path[inactive side] (E)--(S)--(St)--(Et)--cycle;
  \path[inactive side] (S)--(W)--(Wt)--(St)--cycle;
  \coordinate (Wd) at (-1.18,1.55);
  \coordinate (Nd) at (0,1.97);
  \coordinate (Ed) at (1.18,1.55);
  \coordinate (Sd) at (0,1.13);
  \path[active back] (W)--(N)--(Nd)--(Wd)--cycle;
  \path[active back] (N)--(E)--(Ed)--(Nd)--cycle;
  \path[active front] (E)--(S)--(Sd)--(Ed)--cycle;
  \path[active front] (S)--(W)--(Wd)--(Sd)--cycle;
  \path[active top] (Wd)--(Nd)--(Ed)--(Sd)--cycle;
  \draw[slice outline] (Wd)--(Nd)--(Ed)--(Sd)--cycle;
  \path[inactive face] (Wt)--(Nt)--(Et)--(St)--cycle;
  \PatchOutline
  \node[used death] at (0,1.55) {$\times$};
}

\newcommand{\PatchRowIO}{%
\begin{tikzpicture}[x=1cm,y=1cm]
\path[use as bounding box] (-3.25,-1.50) rectangle (8.6,4.18);
\node[patch type label,anchor=east] at (-2.58,1.42) {$\mathsf{IO}$};
\begin{scope}[scale=0.88]\PatchBlock{\IncomingInitial}{\OutgoingAt{3.2}}\end{scope}
\draw[patch decomposition arrow] (1.80,1.42)--(4.43,1.42);
\begin{scope}[shift={(6.15,0)},scale=0.88]\FullPrism\IncomingInitial\OutgoingAt{3.2}\node[patch case label] at (0,-1.20) {no death};\end{scope}
\end{tikzpicture}}

\newcommand{\PatchRowOO}{%
\begin{tikzpicture}[x=1cm,y=1cm]
\path[use as bounding box] (-3.25,-1.50) rectangle (8.6,4.18);
\node[patch type label,anchor=east] at (-2.58,1.42) {$\mathsf{OO}$};
\begin{scope}[scale=0.88]\PatchBlock{\OutgoingAt{0}}{\OutgoingAt{3.2}}\end{scope}
\draw[patch decomposition arrow] (1.80,1.42)--(4.43,1.42);
\begin{scope}[shift={(6.15,0)},scale=0.88]\FullPrism\OutgoingAt{0}\OutgoingAt{3.2}\node[patch case label] at (0,-1.20) {birth, no death};\end{scope}
\end{tikzpicture}}

\newcommand{\PatchRowII}{%
\begin{tikzpicture}[x=1cm,y=1cm]
\path[use as bounding box] (-3.25,-1.50) rectangle (12.95,4.18);
\node[patch type label,anchor=east] at (-2.58,1.42) {$\mathsf{II}$};
\begin{scope}[scale=0.88]\PatchBlock{\IncomingInitial}{\IncomingTerminal}\end{scope}
\draw[patch decomposition arrow] (1.80,1.42)--(4.43,1.42);
\begin{scope}[shift={(5.95,0)},scale=0.88]\FullPrism\IncomingInitial\IncomingTerminal\node[patch case label] at (0,-1.20) {no death};\end{scope}
\begin{scope}[shift={(10.70,0)},scale=0.88]\DeathPrism\IncomingInitial\IncomingTerminal\node[patch case label] at (0,-1.20) {death};\end{scope}
\end{tikzpicture}}

\newcommand{\PatchRowOI}{%
\begin{tikzpicture}[x=1cm,y=1cm]
\path[use as bounding box] (-3.25,-1.50) rectangle (17.3,4.18);
\node[patch type label,anchor=east] at (-2.58,1.42) {$\mathsf{OI}$};
\begin{scope}[scale=0.88]\PatchBlock{\OutgoingAt{0}}{\IncomingTerminal}\end{scope}
\draw[patch decomposition arrow] (1.80,1.42)--(4.43,1.42);
\begin{scope}[shift={(5.85,0)},scale=0.88]\FullPrism\OutgoingAt{0}\IncomingTerminal\node[patch case label] at (0,-1.20) {birth, no death};\end{scope}
\begin{scope}[shift={(10.5,0)},scale=0.88]\DeathPrism\OutgoingAt{0}\IncomingTerminal\node[patch case label] at (0,-1.20) {birth, then death};\end{scope}
\begin{scope}[shift={(15.15,0)},scale=0.88]\EmptyPrism\OutgoingAt{0}\IncomingTerminal\node[patch case label] at (0,-1.20) {split};\end{scope}
\end{tikzpicture}}

\noindent\makebox[\textwidth][l]{\scalebox{0.67}{\PatchRowIO}}
\par\vspace{-0.65em}
\noindent\makebox[\textwidth][l]{\scalebox{0.67}{\PatchRowOO}}
\par\vspace{-0.65em}
\noindent\makebox[\textwidth][l]{\scalebox{0.67}{\PatchRowII}}
\par\vspace{-0.65em}
\noindent\makebox[\textwidth][l]{\scalebox{0.67}{\PatchRowOI}}
\endgroup
\caption{One-site processes compatible with the four patch types. For each type, the diagram on the left specifies the incoming or
outgoing boundary interactions, and the diagrams on the right show the
possible activity trajectories of the patch site that are consistent with
these boundary interactions. Gray shading indicates times at which the patch
site is active.}
\label{fig:patch-interiors}
\end{figure}

\begin{remark}[Meaning of consistency]
\label{rem:consistency}
The skeleton fixes the boundaries of~$P$, and $\operatorname{Con}(P)$ is the
event that the hidden randomness of~$P$ does not contradict them. First, the
patch contains no further successful interaction: every mark from~$i$ with
nonempty target inside the patch occurs while the site is inactive. Second, if
the patch ends with an outgoing interaction, this interaction is successful,
so the site is active just before~$e$. For the patch types this means:
\begin{itemize}
\item $\mathsf{IO}$: the site stays active, so there is no mark from~$i$
inside the patch;
\item $\mathsf{OO}$: in addition, the initial interaction is a birth;
\item $\mathsf{II}$,~$\mathsf{IE}$: no mark with nonempty target occurs before
the first death;
\item $\mathsf{OI}$,~$\mathsf{OE}$: either the initial interaction is a split,
and there is no further constraint, or it is a birth, and the constraint is as
for~$\mathsf{II}$.
\end{itemize}
In Figure~\ref{fig:patch-interiors} we show the resulting trajectories
of~$X^P$.
\end{remark}

The next lemma shows that, conditionally on~$\mathcal G_T$, the patches are
independent with their consistent laws.

\begin{lemma}[Patch factorization]
\label{thm:patch-factorization}
Fix~$T<\infty$. Almost surely,~$\pP_P\bb{\operatorname{Con}(P)}>0$
for every~$P\in\mathcal P_T$. Conditional on~$\mathcal G_T$, the variables~$\Sigma_P$ are
independent, and the conditional law of~$\Sigma_P$ is
$\pP_P^{\mathrm{con}}$ for each~$P\in\mathcal P_T$. Equivalently, for
bounded measurable functions~$f_P$ such that either all~$f_P$ are nonnegative
or the product on the left below is integrable,
\begin{equation}
\label{eq:patch-factorization}
\CE{
\prod_{P\in\mathcal P_T}f_P\bb{\Sigma_P}
}{\mathcal G_T}
=
\prod_{P\in\mathcal P_T}
\E_P^{\mathrm{con}}\Cb{f_P\bb{\Sigma_P}}.
\end{equation}
\end{lemma}

We give the proof in Appendix~\ref{app:conditional-factorization} because
it is purely technical. We verify that, since~$\mathcal{G}_T$ already
reveals all the ways different sites interact, the process in each patch
evolves independently and according to the conditional
law~$\pP_P^{\operatorname{con}}$.

\subsection{Patch contributions and representation}
\label{subsec:patch-representation}

We now combine the patch factorization Lemma~\ref{thm:patch-factorization} and the Feynman Kac formula~\eqref{eq:monomial-feynman-kac-duality} to obtain the patch representation Theorem~\ref{thm:patch-representation}. The step is conceptually simple: the random variable 
\begin{equation*}
Z_t^\eta
=
\sigma_t
\exp\bb{\int_0^tV(A_u)\,du}
\chi_{A_t}(\eta),
\end{equation*}
from the FK formula, factorizes into contributions from each patch (Lemma~\ref{lem:pathwise-factorization}). The process law conditional on~$\mathcal{G}_T$ factorizes over patches
(Lemma~\ref{thm:patch-factorization}), so we can average each factor
independently with the corresponding law~$\pP_P^{\mathrm{con}}$. In this
subsection we mainly work out the patch contributions to the
representation. In Lemma~\ref{lem:contribution-formulas} we give explicit
formulas for them.

\begin{definition}[Patch weights]
\label{def:patch-weights}
The \emph{sign} of a patch~$P$ is the sign of the interaction that begins it, given by
\begin{equation*}
\sigma(P)
=
\pw{+}{\mathsf X(P)=\mathsf I}
{\sigma_{i(P)}^{\alpha(P)}\bb{S(P)}}{\mathsf X(P)=\mathsf O},
\end{equation*}
For a patch~$P\in\mathcal P$, set
\begin{equation*}
F(P)
=
\sigma(P)
\exp\bb{
V_{i(P)}\int_{s(P)}^{e(P)}X_u^P\,du
}.
\end{equation*}
For an end patch~$P\in\mathcal E_t$
(see Definition~\ref{def:bulk-end-patches}) and~$z\in[0,1]$, set
\begin{equation*}
F(z,P)
=
\sigma(P)
\exp\bb{
V_{i(P)}\int_{s(P)}^tX_u^P\,du
}
z^{X_t^P}.
\end{equation*}
\end{definition}

\begin{lemma}[Pathwise factorization]
\label{lem:pathwise-factorization}
For every~$t\geq0$ and~$\eta\in\Omega$,
\begin{equation}
\label{eq:pathwise-patch-factorization}
Z_t^\eta
=
\prod_{P\in\mathcal B_t}F(P)
\prod_{P\in\mathcal E_t}F\bb{\eta\bb{i(P)},P}.
\end{equation}
\end{lemma}

\begin{proof}[Proof of Lemma~\ref{lem:pathwise-factorization}]
Each successful interaction contributes its sign to the outgoing patch at
its source; deaths have positive sign since
\begin{equation*}
a_i^\delta(\varnothing)
=
\widehat{\lambda_i\rho_i}(\varnothing)
=
\lambda_i(\mb 0)\rho_i(\mb 0)
\geq0.
\end{equation*}
Hence
\begin{equation*}
\sigma_t
=
\prod_{P\in\mathcal B_t}\sigma(P)
\prod_{P\in\mathcal E_t}\sigma(P).
\end{equation*}
Since~$V(A_u)=\sum_{i\in A_u}V_i$, the patch intervals split the potential
integral, while the end patches describe the active set at time~$t$. Hence
\begin{align*}
\int_0^tV(A_u)\,du
&=
\sum_{P\in\mathcal B_t}
V_{i(P)}\int_{s(P)}^{e(P)}X_u^P\,du
+
\sum_{P\in\mathcal E_t}
V_{i(P)}\int_{s(P)}^tX_u^P\,du,
\\
\chi_{A_t}(\eta)
&=
\prod_{P\in\mathcal E_t}
\eta\bb{i(P)}^{X_t^P}.
\end{align*}
Multiplying these identities gives~\eqref{eq:pathwise-patch-factorization}.
\end{proof}

We now compute the conditional expectations of the patch factors. These are the contributions appearing in the patch representation.

\begin{lemma}[Patch contribution formulas]
\label{lem:contribution-formulas}
For~$i\in\Lambda$, let
\begin{equation*}
\alpha_i
=
\sum_{R\subseteq N(i)}\delta_i(R)
+
\sum_{\substack{R\subseteq N(i)\\R\neq\varnothing}}\beta_i(R)
\end{equation*}
be the total update rate for the dual process at site~$i$. Then
$V_i=\alpha_i+a_i^\beta(\varnothing)$. For~$\Delta\in[0,\infty)$
and~$z\in[0,1]$, set
\begin{align}
\varphi_i(\Delta)
&:=
e^{-\alpha_i\Delta}
+
\delta_i(\varnothing)
\int_0^\Delta e^{-\alpha_i u}\,du,
\label{eq:patch-phi}
\\
\psi_i(\Delta,z)
&:=
\delta_i(\varnothing)
\int_0^\Delta e^{a_i^\beta(\varnothing)u}\,du
+
z e^{a_i^\beta(\varnothing)\Delta}.
\label{eq:patch-psi}
\end{align}
Define~$\varphi_i(\infty)$ and~$\psi_i(\infty,z)$ by taking~$\Delta\to\infty$.
For an end patch~$P\in\mathcal E_t$ of type~$\mathsf{IE}$, write~$\Delta(P)=e(P)-s(P)$. Then
\begin{align}
\label{eq:patch-phi-psi-interpretation}
\varphi_i(\Delta(P))
&=
\pP_P\bb{\operatorname{Con}(P)},\notag\\
\psi_i(\Delta(P),\eta(i(P)))
&=
\E_P\Cb{\ind_{\operatorname{Con}(P)}F(\eta(i(P)),P)}.
\end{align}

Let~$P$ be a patch and write
\begin{equation*}
i=i(P),
\qquad
S=S(P),
\qquad
\Delta=e(P)-s(P).
\end{equation*}
\begin{enumerate}[label=\textup{(\alph*)}]
\item The \emph{contribution of a patch}~$P\in\mathcal P$ is
\begin{equation}
\label{eq:full-contribution-formulas}
C(P)
:=
\E_P^{\mathrm{con}}\Cb{F(P)}
=
\begin{cases}
\displaystyle
\frac{\psi_i(\Delta,1)}{\varphi_i(\Delta)},
&\mathsf X(P)\mathsf Y(P)\in\cb{\mathsf{II},\mathsf{IE}},
\\[1.4em]
e^{V_i\Delta},
&\mathsf X(P)\mathsf Y(P)=\mathsf{IO},
\\[0.8em]
\displaystyle
\frac{
\delta_i(S)\sigma_i^\delta(S)
+
\beta_i(S)\sigma_i^\beta(S)\psi_i(\Delta,1)
}{
\delta_i(S)+\beta_i(S)\varphi_i(\Delta)
},
&\mathsf X(P)\mathsf Y(P)\in\cb{\mathsf{OI},\mathsf{OE}},
\\[1.4em]
\sigma_i^\beta(S)e^{V_i\Delta},
&\mathsf X(P)\mathsf Y(P)=\mathsf{OO}.
\end{cases}
\end{equation}
\item The \emph{contribution of an end patch}~$P\in\mathcal E_t$ with site state~$z$ is
\begin{equation}
\label{eq:end-contribution-formulas}
C(z,P)
:=
\E_P^{\mathrm{con}}\Cb{F(z,P)}
=
\pw{
\displaystyle
\frac{\psi_i(\Delta,z)}{\varphi_i(\Delta)}
}{\mathsf X(P)=\mathsf I,}
{
\displaystyle
\frac{
\delta_i(S)\sigma_i^\delta(S)
+
\beta_i(S)\sigma_i^\beta(S)\psi_i(\Delta,z)
}{
\delta_i(S)+\beta_i(S)\varphi_i(\Delta)
}
}{\mathsf X(P)=\mathsf O.}
\end{equation}
\end{enumerate}
\end{lemma}

\begin{proof}[Proof of Lemma~\ref{lem:contribution-formulas}]
By the definition of the consistent patch law, the conditional
expectation of a patch factor is the ratio of~$\E_P\Cb{\ind_{\operatorname{Con}(P)} F(P)}$ or~$\E_P\Cb{\ind_{\operatorname{Con}(P)} F(z,P)}$  to~$\pP_P\bb{\operatorname{Con}(P)}$. We compute the latter first.
Suppose first that the initial interaction is incoming, so~$X_s^P=1$, and the terminal interaction is incoming or absent,
corresponding to~$\mathsf {II}$ or~$\mathsf {IE}$. Consistency occurs when
there is no outgoing interaction before the endpoint or death. Therefore
\begin{equation*}
\pP_P\bb{\operatorname{Con}(P)}
=
e^{-\alpha_i\Delta}
+
\delta_i(\varnothing)
\int_0^\Delta e^{-\alpha_i u}\,du
=
\varphi_i(\Delta).
\end{equation*}
For an incoming end patch, the unnormalized expectation in
\eqref{eq:patch-phi-psi-interpretation} is
\begin{equation*}
\E_P\Cb{F(z,P)\ind_{\operatorname{Con}(P)}}.
\end{equation*}
If there is a death at time~$u$, the site is active
until~$u$ and contributes~$e^{V_i u}$. If there is no death,
the site remains active for the whole patch and contributes~$z e^{V_i\Delta }$. Hence
\begin{equation*}
\E_P\Cb{F(z,P)\ind_{\operatorname{Con}(P)}}
=
\delta_i(\varnothing)
\int_0^\Delta e^{\bb{V_i-\alpha_i}u}\,du
+
z e^{\bb{V_i-\alpha_i}\Delta}
=
\psi_i(\Delta,z).
\end{equation*}
These identities prove~\eqref{eq:patch-phi-psi-interpretation} and the~$\mathsf{IE}$ formula in~\eqref{eq:end-contribution-formulas}.
The same calculation with~$z=1$ gives the bulk patch contributions
for~$\mathsf{II}$ and~$\mathsf{IE}$ in~\eqref{eq:full-contribution-formulas}.

Suppose next that the initial interaction is outgoing and the
terminal interaction is incoming or absent, corresponding to~$\mathsf{OI}$ or~$\mathsf{OE}$. If~$\alpha(P)=\delta$, the source becomes inactive at
the initial interaction, so consistency is automatic and the local factor
is~$\sigma_i^\delta(S)$. If~$\alpha(P)=\beta$, the source remains active
and consistency, again, occurs when there is no outgoing interaction before the endpoint or death. By~\eqref{eq:patch-initial-kind-law},
we have
\begin{align*}
\pP_P\bb{\operatorname{Con}(P)} &= 
\frac{
\delta_i(S)+\beta_i(S)\varphi_i(\Delta)
}{
\delta_i(S)+\beta_i(S)
},\\
\E_P\Cb{F(z,P)\ind_{\operatorname{Con}(P)}} &= \frac{
\delta_i(S)\sigma_i^\delta(S)
+
\beta_i(S)\sigma_i^\beta(S)\psi_i(\Delta,z)
}{
\delta_i(S)+\beta_i(S)
}.
\end{align*}
For an end patch, the ratio gives the~$\mathsf{OE}$ case of
\eqref{eq:end-contribution-formulas}. Setting~$z=1$ gives
\eqref{eq:full-contribution-formulas} for~$\mathsf{OI}$ and~$\mathsf{OE}$. If the terminal interaction is outgoing, consistency event is that the site
remains active throughout the patch. If the initial interaction is
incoming, case~$\mathsf{IO}$, the patch contributes~$e^{V_i\Delta}$. If the initial interaction
is outgoing, case~$\mathsf{OO}$, consistency also forces~$\alpha(P)=\beta$, so the contribution is~$\sigma_i^\beta(S)e^{V_i\Delta}$. These are the remaining two formulas of~\eqref{eq:full-contribution-formulas}.
\end{proof}

\begin{theorem}[Patch representation in finite volume]
\label{thm:patch-representation}
For a probability measure~$\mu$ on~$\Omega$ and~$t\geq0$, let
\begin{equation}
\label{eq:finite-horizon-weight}
W_t^\mu
=
\prod_{P\in\mathcal B_t}C(P)\,
\mu\bb{
\prod_{P\in\mathcal E_t}C\bb{\eta\bb{i(P)},P}
}.
\end{equation}
Then, for every~$R\Subset\Lambda$,~$A\subseteq R$, and~$t\geq0$,
\begin{equation}
\label{eq:patch-representation}
(\mu P_t^{R,0})(\chi_A)
=
\E_A\Cb{W_t^\mu\ind\bb{E_t^R}}.
\end{equation}
\end{theorem}

\begin{proof}[Proof of Theorem~\ref{thm:patch-representation}]
Fix~$\eta\in\Omega$. On~$E_t^R$ we have~$A_s\subseteq R$ for~$s\leq t$, so
\begin{equation*}
\abs{Z_t^\eta\ind(E_t^R)}\leq\exp\bb{t\abs R\sup_i\abs{V_i}}.
\end{equation*}
By
Proposition~\ref{thm:monomial-duality} we obtain
\begin{equation*}
P_t^{R,0}\chi_A(\eta)=\E_A\Cb{Z_t^\eta\ind\bb{E_t^R}}.
\end{equation*}
The~$E_t^R$
is the event that every successful interaction up to time~$t$ has its target
in~$R$, so~$E_t^R\in\mathcal G_t$. Since~$F(P)$ and~$F(z,P)$ are
functions of~$\Sigma_P$, we use the pathwise
factorization~\eqref{eq:pathwise-patch-factorization} and
Lemmas~\ref{thm:patch-factorization} and~\ref{lem:contribution-formulas} to
obtain
\begin{equation*}
\CE{Z_t^\eta\ind\bb{E_t^R}}{\mathcal G_t}
=
\ind\bb{E_t^R}
\prod_{P\in\mathcal B_t}C(P)
\prod_{P\in\mathcal E_t}C\bb{\eta\bb{i(P)},P}
=\ind\bb{E_t^R}\,W_t^{\delta_\eta}.
\end{equation*}
Taking expectations and integrating over~$\eta$ with respect to~$\mu$, which
is valid since the integrand is bounded, proves~\eqref{eq:patch-representation}.
\end{proof}

\subsection{Patch positivity and the threshold}
\label{subsec:proof-patch-positivity-threshold}

We can now interpret patch positivity and the threshold profile that we
introduced in Section~\ref{subsec:patch-positivity}. Patch-positivity ensures that the contribution of bulk patches is non-negative. The profile~$\mb p^\star$ gives the family of measures~$\mathcal M_\star$ for which the contribution of end-patches is non-negative. Thus, under patch-positivity and for initial measures in~$\mathcal M_\star$ the patch representation captures all cancellations from the sign process.

\begin{proposition}[Patch positivity and threshold interpretation]
\label{prop:patch-contribution-positivity}
\leavevmode\par\smallskip
\begin{enumerate}[label=\textup{(\roman*)},leftmargin=*,itemsep=.35em,topsep=0pt,parsep=0pt]
\item The spin system is patch-positive if and only if~$C(P)\geq0$
almost surely for every~$P\in\mathcal P$.
\item If the system is patch-positive, then for each~$i\in\Lambda$, the
threshold~$p_i^\star$ from~\eqref{eq:threshold-profile} is the smallest
$p\in[0,1]$ for which~$C(p,P)\geq0$ for every end patch~$P$ based at~$i$.
Moreover,~$C(z,P)$ is affine and nondecreasing in~$z\in[0,1]$, and
\begin{equation}
\label{eq:threshold-below-zero-neighbour-refresh}
p_i^\star\leq\rho_i(\mb 0).
\end{equation}
\item If the system is patch-positive and~$\mu\in\mathcal M_\star$, then
$W_t^\mu\geq0$ almost surely for every~$t\geq0$.
\end{enumerate}
\end{proposition}

\begin{proof}[Proof of Proposition~\ref{prop:patch-contribution-positivity}]
In terms of the spin rates,
\begin{equation}
\label{eq:dual-to-spin-contribution-data}
\begin{aligned}
\alpha_i
&=
\sum_{R\subseteq N(i)}\abs{\widehat{\lambda_i\rho_i}(R)}
+
\sum_{\substack{R\subseteq N(i)\\R\neq\varnothing}}
\abs{\widehat{\lambda_i}(R)},
\\
a_i^\beta(\varnothing)
&=
-\widehat{\lambda_i}(\varnothing)
=
-\lambda_i(\mb 0),
\\
\delta_i(S)\sigma_i^\delta(S)
&=
\widehat{\lambda_i\rho_i}(S),
&\delta_i(S)
&=
\abs{\widehat{\lambda_i\rho_i}(S)},
\\
\beta_i(S)\sigma_i^\beta(S)
&=
-\widehat{\lambda_i}(S),
&\beta_i(S)
&=
\abs{\widehat{\lambda_i}(S)}.
\end{aligned}
\end{equation}
Since
\begin{equation*}
\widehat{\lambda_i}(\varnothing)=\lambda_i(\mb 0),
\qquad
\widehat{\lambda_i\rho_i}(\varnothing)
=
\lambda_i(\mb 0)\rho_i(\mb 0),
\end{equation*}
we also have
\begin{align*}
\varphi_i(\Delta)
&=
e^{-\alpha_i\Delta}
+
\lambda_i(\mb 0)\rho_i(\mb 0)
\int_0^\Delta e^{-\alpha_i u}\,du,
\\
\psi_i(\Delta,z)
&=
\lambda_i(\mb 0)\rho_i(\mb 0)
\int_0^\Delta e^{-\lambda_i(\mb 0)u}\,du
+
z e^{-\lambda_i(\mb 0)\Delta}.
\end{align*}

If~$\lambda_i(\mb 0)>0$, then
\begin{equation}
\label{eq:all-zero-neighbourhood-relaxation}
\psi_i(u,z)
=
\rho_i(\mb 0)
+
\bb{z-\rho_i(\mb 0)}e^{-\lambda_i(\mb 0)u},
\qquad
\psi_i(u+v,z)
=
\psi_i\bb{u,\psi_i(v,z)}.
\end{equation}
When $\lambda_i(\mb 0)=0$, we have $\psi_i(u,z)=z$, and the second identity
still holds.
Thus, when $\lambda_i(\mb 0)>0$, the value~$z$ is forgotten exponentially
fast and replaced by~$\rho_i(\mb 0)$.

We first prove~\textup{(i)}. Fix a patch~$P$ based at~$i$ and write
$\Delta=e(P)-s(P)\in(0,\infty]$.

Suppose first that the initial boundary is incoming. By
\eqref{eq:full-contribution-formulas},
\begin{equation*}
C(P)
=
\begin{cases}
\displaystyle
\frac{\psi_i(\Delta,1)}{\varphi_i(\Delta)},
&\mathsf X(P)\mathsf Y(P)\in\cb{\mathsf{II},\mathsf{IE}},
\\[1em]
e^{V_i\Delta},
&\mathsf X(P)\mathsf Y(P)=\mathsf{IO}.
\end{cases}
\end{equation*}
For a realized~$\mathsf{II}$ or~$\mathsf{IE}$ patch, the consistency probability~$\varphi_i(\Delta)$ is positive almost surely, and
\begin{equation*}
\delta_i(\varnothing)
=
\widehat{\lambda_i\rho_i}(\varnothing)
=
\lambda_i(\mb 0)\rho_i(\mb 0)
\geq0
\end{equation*}
gives~$\psi_i(\Delta,1)\geq0$. Hence every patch with an incoming
initial boundary has nonnegative contribution.

Now suppose that the initial boundary is outgoing, with nonempty target~$S$.
For an~$\mathsf{OO}$ patch,
\eqref{eq:full-contribution-formulas} and
\eqref{eq:dual-to-spin-contribution-data} give
\begin{equation*}
C(P)
=
\sigma_i^\beta(S)e^{V_i\Delta},
\qquad
\beta_i(S)\sigma_i^\beta(S)
=
-\widehat{\lambda_i}(S)
\geq0.
\end{equation*}
An~$\mathsf{OO}$ patch can occur only when~$\beta_i(S)>0$, so its sign is
nonnegative. For an~$\mathsf{OI}$ or~$\mathsf{OE}$ patch, the denominator in
\eqref{eq:full-contribution-formulas} is positive and its numerator is
\begin{equation}
\widehat{\lambda_i\rho_i}(S)
-
\widehat{\lambda_i}(S)\psi_i(\Delta,1)
=
\widehat{\lambda_i\bb{\rho_i-\rho_i(\mb 0)}}(S)
+
\bb{1-\rho_i(\mb 0)}
\bb{-\widehat{\lambda_i}(S)}
e^{-\lambda_i(\mb 0)\Delta}.
\label{eq:outgoing-full-patch-numerator}
\end{equation}
Here we used~\eqref{eq:all-zero-neighbourhood-relaxation}. The same identity
holds when~$\lambda_i(\mb 0)=0$. Both terms in
\eqref{eq:outgoing-full-patch-numerator} are nonnegative by patch positivity. Thus the coefficient criterion for patch positivity implies~$C(P)\geq 0$.

Conversely if~$C(P)\geq 0$ for every patch, then nonnegativity of the~$\mathsf{OO}$ contributions forces
$-\widehat{\lambda_i}(S)\geq0$ for every nonempty~$S$.
If~$\lambda_i(\mb0)=0$, this and~$\lambda_i\geq0$ force~$\lambda_i\equiv0$.
Otherwise, letting~$\Delta\to\infty$ in
\eqref{eq:outgoing-full-patch-numerator} gives
$\widehat{\lambda_i(\rho_i-\rho_i(\mb0))}(S)\geq0$ for every~$S$
with~$\delta_i(S)+\beta_i(S)>0$.
For the remaining~$S$ both coefficients vanish.
These are exactly the two patch-positivity inequalities.

We next prove~\textup{(ii)} by identifying the threshold for end patches. If
$\lambda_i(\mb 0)=0$, then patch positivity and~$\lambda_i\geq0$ force
$\lambda_i\equiv0$. Indeed, if~$T\subseteq N(i)$ is the set of coordinates
equal to~$1$ in a neighbourhood configuration, then
\begin{equation*}
0
\leq
\lambda_i(\eta)
=
\sum_{S\subseteq T}\widehat{\lambda_i}(S)
\leq0.
\end{equation*}
Thus~$\lambda_i\rho_i\equiv0$, there are no outgoing interactions based
at~$i$, and~\eqref{eq:threshold-profile} gives~$p_i^\star=0$. Going forward assume
that~$\lambda_i(\mb 0)>0$.

An~$\mathsf{IE}$ patch has, by~\eqref{eq:end-contribution-formulas},
\begin{equation*}
C(p,P)
=
\frac{\psi_i(\Delta,p)}{\varphi_i(\Delta)}
\geq0
\qquad
\text{for every }p\in[0,1],
\end{equation*}
because by~\eqref{eq:all-zero-neighbourhood-relaxation} the
$\psi_i(\Delta,p)$ is a convex combination of~$p$ and~$\rho_i(\mb 0)$.
Thus only~$\mathsf{OE}$ patches can have negative contributions and impose a restriction. For such a patch
with target~$S$, the sign of~$C(p,P)$ is the sign of
\begin{equation*}
N_{i,S}(\Delta,p)
=
\widehat{\lambda_i\rho_i}(S)
-
\widehat{\lambda_i}(S)\psi_i(\Delta,p).
\end{equation*}
As a function of~$e^{-\lambda_i(\mb 0)\Delta}$, this is affine, with endpoint
values
\begin{align*}
N_{i,S}(0,p)
&=
\widehat{\lambda_i\rho_i}(S)
-
p\widehat{\lambda_i}(S),
\\
\lim_{\Delta\to\infty}N_{i,S}(\Delta,p)
&=
\widehat{\lambda_i\bb{\rho_i-\rho_i(\mb 0)}}(S)
\geq0.
\end{align*}
If~$\widehat{\lambda_i}(S)=0$, patch positivity gives
$\widehat{\lambda_i\rho_i}(S)\geq0$, there is no restriction on~$p$.
If~$\widehat{\lambda_i}(S)<0$, nonnegativity for every~$\Delta$ is
equivalent to
\begin{equation*}
p
\geq
\frac{\widehat{\lambda_i\rho_i}(S)}
{\widehat{\lambda_i}(S)}.
\end{equation*}
Taking the maximum over all nonempty~$S$ gives
\eqref{eq:threshold-profile}. This threshold is sharp: if~$p<p_i^\star$,
then for some~$S$ with~$\widehat{\lambda_i}(S)<0$ we have
$N_{i,S}(0,p)<0$ and so~$C(p,P)<0$ for short~$\mathsf{OE}$ patches with~$S(P)=S$.

The same formulas also show that every end contribution is nondecreasing in its
terminal argument. For an~$\mathsf{IE}$ patch,
\begin{equation*}
\partial_zC(z,P)
=
\frac{e^{-\lambda_i(\mb 0)\Delta}}{\varphi_i(\Delta)}
\geq0,
\end{equation*}
while for an~$\mathsf{OE}$ patch,
\begin{equation*}
\partial_zC(z,P)
=
\frac{
\bb{-\widehat{\lambda_i}(S)}
e^{-\lambda_i(\mb 0)\Delta}
}{
\delta_i(S)+\beta_i(S)\varphi_i(\Delta)
}
\geq0.
\end{equation*}

Finally, if~$\widehat{\lambda_i}(S)<0$, patch positivity gives
\begin{equation*}
\widehat{\lambda_i\rho_i}(S)
-
\rho_i(\mb 0)\widehat{\lambda_i}(S)
\geq0.
\end{equation*}
Dividing by~$\widehat{\lambda_i}(S)<0$ yields
\begin{equation*}
\frac{\widehat{\lambda_i\rho_i}(S)}
{\widehat{\lambda_i}(S)}
\leq
\rho_i(\mb 0).
\end{equation*}
Taking the maximum in~\eqref{eq:threshold-profile} proves
\eqref{eq:threshold-below-zero-neighbour-refresh}.

Finally, we prove~\textup{(iii)}. By~\textup{(i)}, every bulk patch (see Definition~\ref{def:bulk-end-patches}) has
nonnegative contribution. For an end patch~$P$ based at~$i$, part~\textup{(ii)}
gives that~$C(z,P)$ is affine and nondecreasing, with
$C(p_i^\star,P)\geq0$. Set~$\kappa(P)=\partial_zC(z,P)\geq0$. Then
\begin{equation*}
C\bb{\eta(i),P}
=
C\bb{p_i^\star,P}
+
\kappa(P)\bb{\eta(i)-p_i^\star}.
\end{equation*}
Distinct end patches are based at distinct sites. Therefore
\begin{equation}
\label{eq:centered-end-expansion}
\mu\bb{
\prod_{P\in\mathcal E_t}
C\bb{\eta\bb{i(P)},P}
} =
\sum_{\mathcal Q\subseteq\mathcal E_t}
\mu\bb{
\chi_{\cb{i(P):P\in\mathcal Q}}^\star
}
\prod_{P\in\mathcal Q}\kappa(P)
\prod_{P\in\mathcal E_t\setminus\mathcal Q}
C\bb{p_{i(P)}^\star,P}.
\end{equation}
For~$\mu\in\mathcal M_\star$ every term is nonnegative, so~$W_t^\mu\geq0$.
\end{proof}

\section{Proofs}
\label{sec:proofs}

In this Section we give proofs for the results of Section~\ref{sec:main-results}. We prove the main result,  Theorem~\ref{thm:common-invariant-limit} by adapting the classical
duality argument for ergodicity. After a finite-volume approximation, the patch
representation splits~$P_t^{R,0}\chi_A$ into two terms: one that is weakly dependent on initial law and another that decays.  
We give the decomposition and state the two estimates, then prove the theorem.
In Sections~\ref{subsec:proof-pure-death-comparison}
and~\ref{subsec:terminal-relaxation} we prove the two estimates in detail.
In Section~\ref{subsec:ergodicity-corollaries} we give proofs of the two corollaries, which upgrade results of Theorem~\ref{thm:common-invariant-limit} to ergodicity. We obtain Corollary~\ref{cor:uniform-exponential-ergodicity} directly as a
special case of Theorem~\ref{thm:common-invariant-limit}. We prove
Corollary~\ref{cor:oriented-ergodicity} by a burn-in argument showing that
with orientable site graph all measures converge towards~$\mathcal M_\star$,
and from there towards the unique invariant measure by
Theorem~\ref{thm:common-invariant-limit}.

\subsection{Proof of the common invariant limit, Theorem~\ref{thm:common-invariant-limit}}
\label{subsec:duality-proof-plan}
\label{subsec:proof-common-invariant-limit}
\label{subsec:infinite-volume-proof}
\label{subsec:proofs-main}

The classical duality argument gives ergodicity, assuming~$V(A)\leq -\gamma\ind\bb{A\neq\vn}$, for some~$\gamma>0$. Let
\begin{equation*}
\tau=\inf\cb{s\geq0:A_s=\vn}.
\end{equation*}
The empty dual set is absorbing and~$V(\vn)=0$, so the duality, which holds in infinite volume because the weight is bounded (Remark~\ref{rem:fk-infinite-volume}), gives
\begin{equation}
P_t\chi_A(\eta)
=
\E_A\Cb{\sigma_\tau
 e^{\int_0^\tau V(A_s)\,ds}\ind\bb{\tau\leq t}}
+\E_A\Cb{\sigma_t
 e^{\int_0^t V(A_s)\,ds}\chi_{A_t}(\eta)
 \ind\bb{\tau>t}}.
\label{eq:classical-duality-split}
\end{equation}
The first term is independent of~$\eta$ and converges as~$t\to\infty$. The second has absolute value at most~$e^{-\gamma t}$. Hence for each~$A$ the~$P_t\chi_A(\eta)$ converges to a limit independent of the initial configuration. The monomials are a basis of space of local functions, so uniform exponential ergodicity at rate~$\gamma$ follows.

The argument above fails if the potential~$V$ can be positive, as we lose the decay of the~$\eta$-dependent term in~\eqref{eq:classical-duality-split}. The patch representation gives us a better control over this term, because patch contributions in Theorem~\ref{thm:patch-representation} capture the cancellations from sign~$\sigma_t$, see Proposition~\ref{prop:patch-contribution-positivity}. 

The only technical snag is that we still need to split~$P_t\chi_A$ into three terms: one for finite-volume approximation, one weakly dependent on the initial measure, and one decaying. We cannot use the event~$\cb{\tau >t}$ because it is not in~$\mathcal{G}_t$. Knowing when and where the successful interactions happen before time~$t$ does not determine whether the dual is dead by time~$t$. We solve this problem by using the proxy event~$\cb{\text{there is no successful
interaction in }(T,t]}$, which belongs to~$\mathcal G_t$. The proxy event takes the place of~$\cb{\tau\leq t}$. With no new successful interactions, no new site enters the dual process after time~$T$, and all end patches are already present at time~$T$ and span the entire interval~$[T,t)$ (see Definition~\ref{def:bulk-end-patches}). The dependence on the initial measure sits only in the end patches and decays exponentially with end patch length at rate~$\lambda_i(\mb 0)\geq\varepsilon$. We thus show that the contribution of the proxy event converges exponentially fast to a quantity independent of the initial measure. The complementary event, in which there is a successful interaction after time~$T$, takes the place of~$\cb{\tau>t}$. The pure deaths enter the dual through the potential, making it smaller. Comparing with a system without the pure-death component gives an upper bound of~$e^{-\varepsilon T}$ on the contribution of the complementary event.

Fix~$A\subseteq R\Subset\Lambda$,~$\mu\in\mathcal M_\star$, and
$0\leq T<t$. We use the signed dual started from~$A$ and its patch
family from Section~\ref{sec:probabilistic-interpretation}. Recall
from~\eqref{eq:confinement-event} that~$E_t^R$ is the event that the dual process
stays in~$R$ up to time~$t$. Define
\begin{align}
W&=\prod_{P\in\mathcal P}C(P)\,
\ind\bb{\abs{\mathcal P}<\infty},\notag\\
L_{T,t}
&=\cb{\text{no successful interaction occurs in }(T,t]},
\notag\\
L_T
&=\cb{\text{no successful interaction occurs after }T}.
\label{eq:no-late-interaction-events}
\end{align}
The random variable~$W$ does not depend on~$\mu$ and plays the role
of~$\sigma_\tau e^{\int_0^\tau V(A_s)\,ds}$ in~\eqref{eq:classical-duality-split}.
The event~$L_{T,t}$ belongs to~$\mathcal G_t$ and~$L_T$ to~$\mathcal G_\infty$.
By Theorem~\ref{thm:patch-representation}, we have
\begin{align}
(\mu P_t)(\chi_A)
={}&\mu\bb{\bb{P_t-P_t^{R,0}}\chi_A}
\notag\\
&+\E_A\Cb{W_t^\mu\ind\bb{E_t^R\cap L_{T,t}}}
 +\E_A\Cb{W_t^\mu\ind\bb{E_t^R\cap L_{T,t}^c}}.
\label{eq:two-term-proof-plan}
\end{align}
The first term is the finite-volume error, which is negligible when~$R$ is chosen as a ball around~$A$ whose radius grows
linearly in~$t$, cf. Lemma~\ref{lem:finite-propagation}. The second term corresponds to the convergent term
in~\eqref{eq:classical-duality-split}; we show that it converges to the~$\mu$-independent~$\E_A\Cb{W\ind\bb{E_T^R\cap L_T}}$. The third term corresponds to the decaying term
in~\eqref{eq:classical-duality-split} and we show that it is at most~$e^{-\varepsilon T}$.

\begin{lemma}[Late interactions]
\label{lem:late-interactions}
Under the assumptions of Theorem~\ref{thm:common-invariant-limit}, for any $A\subseteq R\Subset\Lambda$, $\mu\in\mathcal M_\star$ and $0\leq T<t$ we have
\begin{align}
0\leq\E_A\Cb{W_t^\mu\ind\bb{E_t^R\cap L_{T,t}^c}}
&\leq e^{-\varepsilon T},
\label{eq:finite-late-interaction-bound}\\
0\leq\E_A\Cb{W\ind\bb{L_T^c}}
&\leq e^{-\varepsilon T}.
\label{eq:limiting-late-interaction-bound}
\end{align}
\end{lemma}

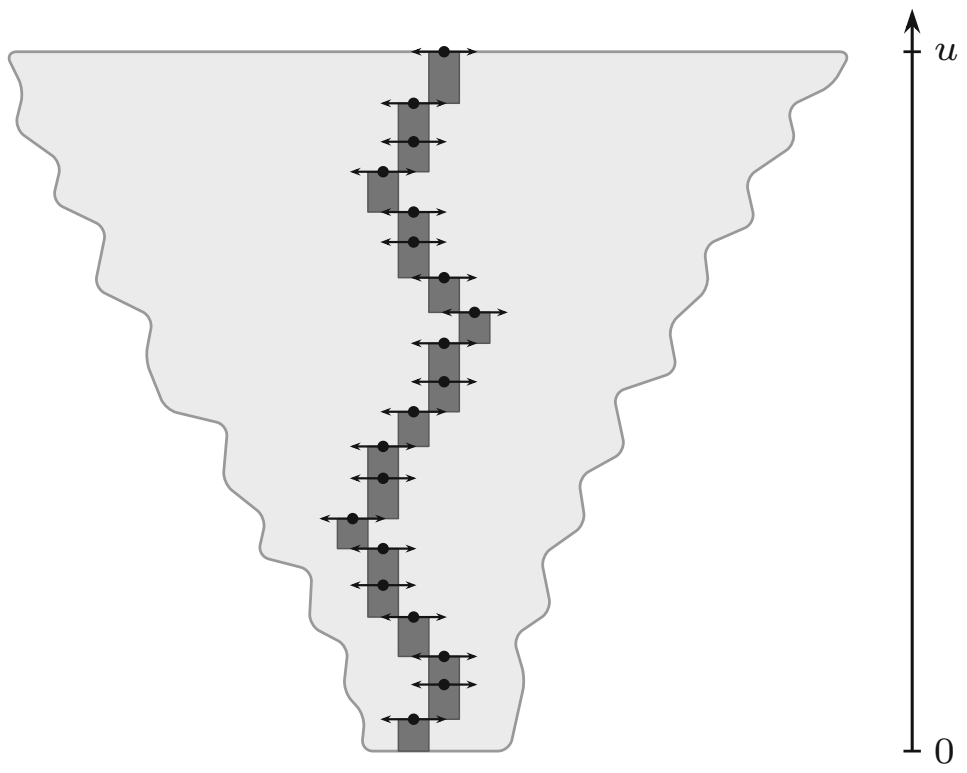
\begin{figure}[tbp]
\centering
\begingroup
\definecolor{convtrailcone}{gray}{0.92}
\definecolor{convtrailconeborder}{gray}{0.60}
\definecolor{convtrailpatch}{gray}{0.46}
\definecolor{convtrailpatchborder}{gray}{0.34}
\definecolor{convtrailink}{gray}{0.08}
\begin{tikzpicture}[x=0.48cm,y=0.62cm]
\path[use as bounding box] (0.20,-0.12) rectangle (13.64,8.02);
\tikzset{
  convtrail cone/.style={fill=convtrailcone,draw=convtrailconeborder,
    line width=0.52pt,line join=round,rounded corners=2.2pt},
  convtrail patch/.style={fill=convtrailpatch,draw=convtrailpatchborder,
    line width=0.28pt,line join=round},
  convtrail arrow/.style={draw=convtrailink,line width=0.42pt,
    -{Stealth[length=0.72mm,width=0.52mm]}},
  convtrail source/.style={fill=convtrailink,draw=none,circle,
    inner sep=0pt,minimum size=2.15pt},
  convtrail axis/.style={draw=convtrailink,line width=0.58pt,
    -{Stealth[length=1.65mm,width=1.08mm]}},
  convtrail tick/.style={draw=convtrailink,line width=0.52pt},
  convtrail label/.style={font=\scriptsize,text=convtrailink,inner sep=1pt}
}

\path[convtrail cone]
  (5.43,0)
  -- (5.48,0.38)
  -- (5.18,0.64)
  -- (5.25,1.12)
  -- (4.71,1.34)
  -- (4.75,1.93)
  -- (4.00,2.08)
  -- (4.12,2.49)
  -- (3.52,2.86)
  -- (3.59,3.48)
  -- (2.73,3.64)
  -- (2.45,4.18)
  -- (2.56,4.63)
  -- (1.75,4.94)
  -- (1.92,5.58)
  -- (1.17,5.86)
  -- (1.30,6.28)
  -- (0.65,6.64)
  -- (0.78,7.05)
  -- (0.52,7.46)
  -- (12.18,7.46)
  -- (11.92,7.10)
  -- (11.30,6.87)
  -- (11.42,6.48)
  -- (10.72,6.11)
  -- (10.86,5.63)
  -- (10.15,5.39)
  -- (10.22,4.94)
  -- (9.66,4.54)
  -- (9.78,4.05)
  -- (8.90,3.82)
  -- (9.07,3.20)
  -- (8.42,2.92)
  -- (8.52,2.42)
  -- (7.90,2.08)
  -- (8.05,1.50)
  -- (7.52,1.22)
  -- (7.69,0.77)
  -- (7.47,0)
  -- cycle;

\foreach \x/\ya/\yb in {
  6.15/0.00/0.34,
  6.57/0.34/0.71,
  6.57/0.71/1.01,
  6.15/1.01/1.43,
  5.73/1.43/1.77,
  5.73/1.77/2.16,
  5.31/2.16/2.48,
  5.73/2.48/2.91,
  5.73/2.91/3.25,
  6.15/3.25/3.62,
  6.57/3.62/3.94,
  6.57/3.94/4.35,
  6.99/4.35/4.68,
  6.57/4.68/5.05,
  6.15/5.05/5.43,
  6.15/5.43/5.75,
  5.73/5.75/6.18,
  6.15/6.18/6.50,
  6.15/6.50/6.91,
  6.57/6.91/7.46
}
{
  \path[convtrail patch]
    ({\x-0.21},\ya) rectangle ({\x+0.21},\yb);
}

\foreach \x/\y in {
  6.15/0.34,6.57/0.71,6.57/1.01,6.15/1.43,5.73/1.77,
  5.73/2.16,5.31/2.48,5.73/2.91,5.73/3.25,6.15/3.62,
  6.57/3.94,6.57/4.35,6.99/4.68,6.57/5.05,6.15/5.43,
  6.15/5.75,5.73/6.18,6.15/6.50,6.15/6.91,6.57/7.46
}
{
  \draw[convtrail arrow] (\x,\y)--({\x-0.46},\y);
  \draw[convtrail arrow] (\x,\y)--({\x+0.46},\y);
  \node[convtrail source] at (\x,\y) {};
}

\draw[convtrail axis] (13.02,0)--(13.02,7.92);
\draw[convtrail tick] (12.90,0)--(13.14,0);
\draw[convtrail tick] (12.90,7.46)--(13.14,7.46);
\node[convtrail label,anchor=west] at (13.23,7.46) {$u$};
\node[convtrail label,anchor=west] at (13.23,0) {$0$};
\end{tikzpicture}
\endgroup
\caption{The ancestry of a successful interaction at time~$u$. Each dark
rectangle is a patch ending at an outgoing successful interaction. The patch
lifetimes sum to~$u$.}
\label{fig:backward-outgoing-chain}
\end{figure}

To see why this bound holds, consider the spin system with the same generator as~$\mathcal{L}$ but without the pure-death component. We show the two systems have the same dual process, but the less noisy one has larger potential~$V$. Thus they have the same random patches and the patch contributions for the noiseless process are larger. For a patch~$P$ ending with an outgoing interaction, the two
contributions differ by the factor~$e^{-\varepsilon(e(P)-s(P))}$.
A successful interaction at time~$u>T$ has an ancestry of such patches
(Figure~\ref{fig:backward-outgoing-chain}) whose lengths sum to~$u$. Their product therefore supplies the factor
\begin{equation*}
\prod_{P\text{ in the ancestry}}
e^{-\varepsilon(e(P)-s(P))}=e^{-\varepsilon u}\leq e^{-\varepsilon T}.
\end{equation*}
The corresponding term for the system without deaths is at most~$1$,
which gives the estimate. The proof of Lemma~\ref{lem:late-interactions} is
given in Section~\ref{subsec:proof-pure-death-comparison} below.

\FloatBarrier

\begin{lemma}[Terminal relaxation]
\label{lem:no-late-relaxation}
Assume the hypotheses of Theorem~\ref{thm:common-invariant-limit}. For
$A\subseteq R\Subset\Lambda$, $\mu\in\mathcal M_\star$ and $0\leq T<t$, we have
\begin{equation}
\label{eq:no-late-relaxation-bound}
\begin{aligned}
\abs{
\E_A\Cb{W_t^\mu\ind\bb{E_t^R\cap L_{T,t}}}
-\E_A\Cb{W\ind\bb{E_T^R\cap L_T}}
}
&\leq(1+\abs R)e^{-\varepsilon(t-T)},\\
\E_A\Cb{W\ind\bb{(E_T^R)^c\cap L_T}}
&\leq\norm{\bb{P_T-P_T^{R,0}}\chi_A}_\infty.
\end{aligned}
\end{equation}
\end{lemma}

\begin{figure}[tbp]
\centering
\begingroup
\definecolor{patchfigred}{RGB}{175,32,38}
\definecolor{patchfigactive}{gray}{0.86}
\definecolor{patchfigoutline}{gray}{0.15}
\begin{tikzpicture}[x=0.63cm,y=0.66cm]
\path[use as bounding box] (-0.08,-0.15) rectangle (20.58,10.66);
\tikzset{
  patchfig site/.style={draw=patchfigoutline,line width=0.48pt},
  patchfig success/.style={draw=black,line width=0.86pt,
    -{Stealth[length=2.0mm,width=1.35mm]}},
  patchfig active/.style={fill=patchfigactive,draw=none},
  patchfig horizon/.style={draw=patchfigoutline,line width=0.72pt},
  patchfig boundary/.style={draw=patchfigoutline,line width=0.40pt},
  patchfig cut/.style={draw=patchfigoutline,line width=1.25pt,
    dash pattern=on 4.2pt off 2.6pt,line cap=round},
  patchfig small/.style={font=\scriptsize,text=patchfigoutline,inner sep=0pt}
}

\newcommand{\ConvEndPatchBarrier}[3]{%
  \path[
    pattern={Lines[angle=39,distance=3.10pt,line width=1.40pt]},
    pattern color=patchfigred
  ] (#1,#2) rectangle ({#1+0.105},#3);
  \path[
    pattern={Lines[angle=-39,distance=3.10pt,line width=1.40pt]},
    pattern color=patchfigred
  ] ({#1+0.895},#2) rectangle ({#1+1},#3);}

\newcommand{\ConvEndInteractionBoundary}[2]{%
  \draw[patchfig boundary] ({#1-1},#2)--({#1+2},#2);}
\newcommand{\ConvEndInteractionMark}[2]{%
  \draw[patchfig success] ({#1+0.5},#2)--({#1-0.5},#2);
  \draw[patchfig success] ({#1+0.5},#2)--({#1+1.5},#2);
  \fill[black] ({#1+0.5},#2) circle[radius=1.25pt];}

\foreach \x/\a/\b in {
  9/0.00/0.48,10/0.48/0.92,11/0.92/1.20,8/0.48/1.38,
  9/1.38/1.58,7/1.38/1.74,12/1.20/2.18,10/1.58/2.40,
  6/1.74/2.72,8/1.74/2.88,13/2.18/2.98,5/2.72/3.22,
  11/2.40/3.44,14/2.98/3.70,6/3.22/3.94,15/3.70/4.02,
  4/3.22/4.12,12/3.44/4.36,9/2.88/4.58,3/4.12/4.68,
  16/4.02/4.86,2/4.68/5.02,14/4.02/5.18,4/4.68/5.28,
  12/4.36/5.32,17/4.86/5.38,7/3.94/5.50,10/4.58/5.66,
  1/5.02/5.74,18/5.38/5.86,15/5.18/5.98,5/5.28/6.10
}{\path[patchfig active] (\x,\a) rectangle ({\x+1},\b);}

\foreach \x/\s in {
  0/5.74,1/5.02,2/4.68,3/4.12,4/3.22,
  5/2.72,6/1.74,7/1.38,8/0.48,9/0.00,
  10/0.00,11/0.92,12/1.20,13/2.18,14/2.98,
  15/3.70,16/4.02,17/4.86,18/5.38,19/5.86
}{\ConvEndPatchBarrier{\x}{\s}{10.50}}

\foreach \x/\s in {
  0/5.74,1/5.02,2/4.68,3/4.12,4/3.22,
  5/2.72,6/1.74,7/1.38,8/0.48,9/0.00,
  10/0.00,11/0.00,12/0.92,13/1.20,14/2.18,
  15/2.98,16/3.70,17/4.02,18/4.86,19/5.38,
  20/5.86
}{\draw[patchfig site](\x,\s)--(\x,10.50);}
\draw[patchfig horizon](0,0)--(20,0);
\draw[patchfig horizon](0,10.50)--(20,10.50);

\foreach \s/\y in {
  9/0.48,10/0.92,11/1.20,8/1.38,9/1.58,7/1.74,
  12/2.18,10/2.40,6/2.72,8/2.88,13/2.98,5/3.22,
  11/3.44,14/3.70,6/3.94,15/4.02,4/4.12,12/4.36,
  9/4.58,3/4.68,16/4.86,2/5.02,14/5.18,4/5.28,
  12/5.32,17/5.38,7/5.50,10/5.66,1/5.74,18/5.86,
  15/5.98,5/6.10
}{\ConvEndInteractionBoundary{\s}{\y}}

\foreach \s/\y in {
  9/0.48,10/0.92,11/1.20,8/1.38,9/1.58,7/1.74,
  12/2.18,10/2.40,6/2.72,8/2.88,13/2.98,5/3.22,
  11/3.44,14/3.70,6/3.94,15/4.02,4/4.12,12/4.36,
  9/4.58,3/4.68,16/4.86,2/5.02,14/5.18,4/5.28,
  12/5.32,17/5.38,7/5.50,10/5.66,1/5.74,18/5.86,
  15/5.98,5/6.10
}{\ConvEndInteractionMark{\s}{\y}}

\draw[patchfig cut] (0,6.35)--(20,6.35);
\node[patchfig small,anchor=west] at (20.18,0) {$0$};
\node[patchfig small,anchor=west] at (20.18,6.35) {$T$};
\node[patchfig small,anchor=west] at (20.18,10.50) {$t$};
\end{tikzpicture}
\endgroup
\caption{On~$L_{T,t}$, no successful interaction occurs in~$(T,t]$. The bulk
patches are unchanged after~$T$, and every end patch has length at least~$t-T$.}
\label{fig:no-late-end-patches}
\end{figure}
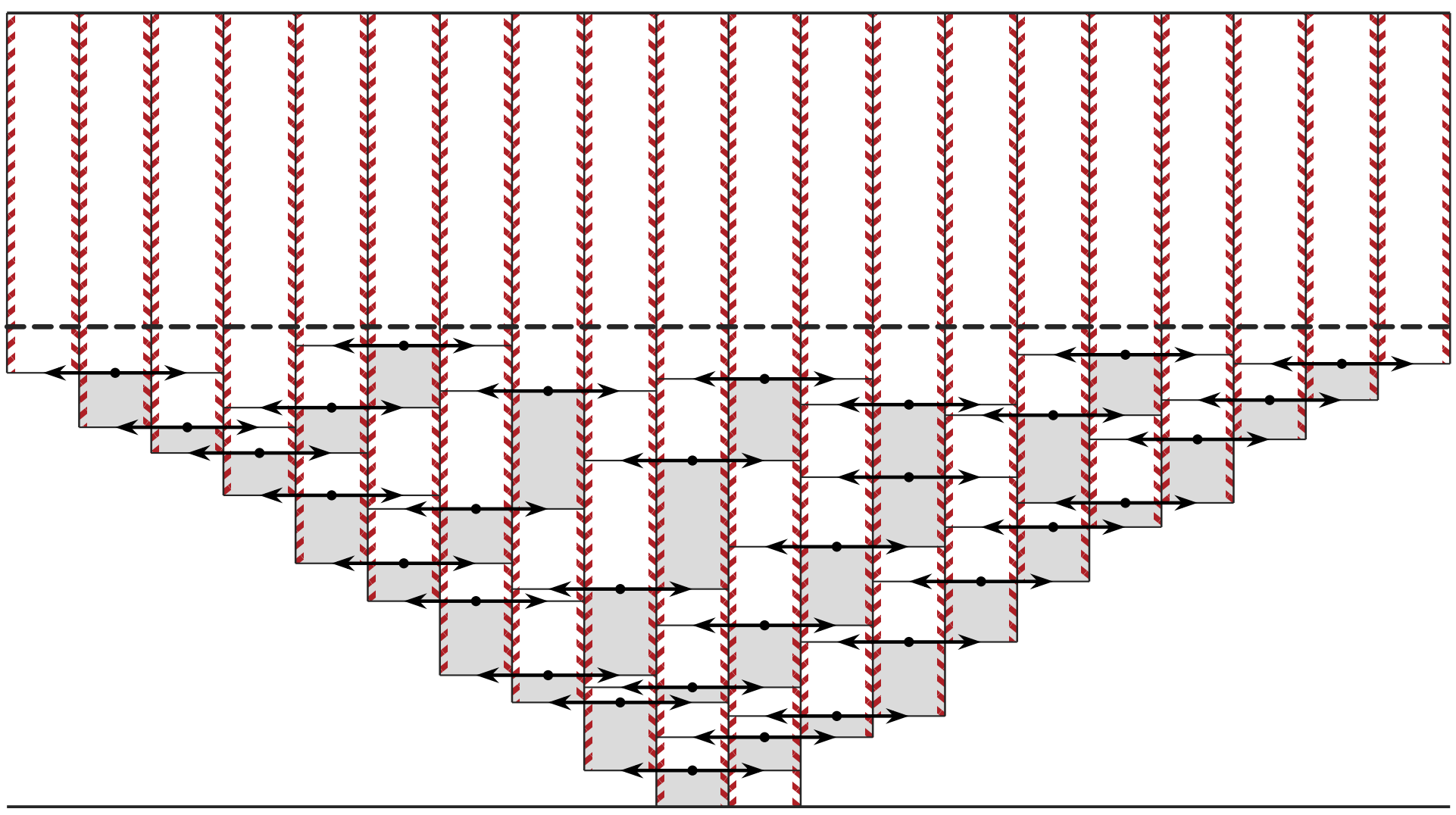

On~$L_{T,t}$, the bulk patches are unchanged after~$T$ and every end patch
present at~$T$ extends to time~$t$ (Figure~\ref{fig:no-late-end-patches}). Only contributions from end patches retain a dependence
on the initial law, and that dependence decays exponentially with end patch length. On~$L_{T,t}$ all end patches have length at least~$t-T$, which gives the decay rate~$\varepsilon(t-T)$. The bound gets worse with the number of end patches, but on~$E_T^R$, there are at most~$\abs R$ of them.
In Section~\ref{subsec:terminal-relaxation} we prove Lemma~\ref{lem:no-late-relaxation}.

\FloatBarrier

\Needspace{12\baselineskip}
\begin{proof}[Proof of Theorem~\ref{thm:common-invariant-limit}]
Fix~$A\Subset\Lambda$,~$\mu\in\mathcal M_\star$, and~$t>0$. Set~$T=t/2$ and
$R=B\bb{A,\lceil vt\rceil}$, with~$v$ from Lemma~\ref{lem:finite-propagation}(a)
for~$a=\varepsilon$. By~\eqref{eq:two-term-proof-plan} and the preceding estimates,
\begin{align*}
(\mu P_t)(\chi_A)
&=
\underbrace{\mu\bb{\bb{P_t-P_t^{R,0}}\chi_A}}_{\abs{\cdot}\lesssim e^{-\varepsilon t}}
+\underbrace{\E_A\Cb{W_t^\mu\ind\bb{E_t^R\cap L_{t/2,t}^c}}}_{0\leq\cdot\leq e^{-\varepsilon t/2}}\\
&\quad+\E_A\Cb{W_t^\mu\ind\bb{E_t^R\cap L_{t/2,t}}}\\
&=\E_A\Cb{W\ind\bb{L_{t/2}}} + O\bb{\abs{R}e^{-\varepsilon t/2}}\\
&=\E_A[W] + O\bb{\abs{R}e^{-\varepsilon t/2}}.
\end{align*}

The first two underbraced bounds are, respectively, the statement of Lemma~\ref{lem:finite-propagation}(a) and the inequality
\eqref{eq:finite-late-interaction-bound} in Lemma~\ref{lem:late-interactions}. The second equality follows from Lemma~\ref{lem:no-late-relaxation} and the same bound for the finite-volume approximation, Lemma~\ref{lem:finite-propagation}(a). The third follows from inequality~\eqref{eq:limiting-late-interaction-bound} of Lemma~\ref{lem:late-interactions}. The polynomial growth assumption gives $\abs{R}\lesssim \bb{1+t}^D$. Thus, for every $A\Subset \Lambda$ there exists $K_A\geq 0$ such that
\begin{equation}
\label{eq:m-star-monomial-convergence}
\abs{(\mu P_t)(\chi_A)-\E_A[W]}
\leq K_A(1+t)^D e^{-\varepsilon t/2}.
\end{equation}
The bound is uniform in~$\mu\in\mathcal M_\star$. The space of probability measures on~$\Omega$ is weakly
compact, and monomial moments determine a law. Thus~$\mu P_t\Rightarrow\pi$
for a law~$\pi$ satisfying~$\pi(\chi_A)=\E_A[W]$ for every finite~$A$.
For every~$s\geq0$, the Feller property and the semigroup identity give
\begin{equation*}
\pi P_s=\lim_{t\to\infty}(\mu P_t)P_s=\lim_{t\to\infty}\mu P_{t+s}=\pi.
\end{equation*}

Every local function is a finite linear combination of monomials. Hence
\eqref{eq:m-star-monomial-convergence} gives,
for every local~$f$,
\begin{equation}
\label{eq:m-star-local-convergence}
\sup_{\mu\in\mathcal M_\star}
\abs{(\mu P_t)(f)-\pi(f)}
\leq
K_f(1+t)^D e^{-\varepsilon t/2}
\end{equation}
for some~$K_f<\infty$.

To extend the bound to~$\mathcal M_-$, let~$K\geq0$ and
$\mu\in\mathcal M_{-,K}$, and set
\begin{equation*}
\overline\mu
=
\frac{\mu+K\mu_{\mb 1}}{1+K}
\in\mathcal M_\star, \qquad \text{so that} \qquad \mu
=
(1+K)\overline\mu-K\mu_{\mb 1}.
\end{equation*}
Both~$\overline{\mu}$ and~$\mu_{\mb 1}$ are in~$\mathcal{M}_*$, so by~\eqref{eq:m-star-local-convergence} we have
\begin{equation}
\label{eq:m-minus-reduction}
\sup_{\mu\in\mathcal M_{-,K}}
\abs{(\mu P_t)(f)-\pi(f)}
\leq(1+2K)K_f(1+t)^D e^{-\varepsilon t/2}.
\end{equation}
This proves~\eqref{eq:uniform-m-minus-convergence}. Since~$\mathcal M_- = \bigcup_{K\geq0}\mathcal M_{-,K}$, it also follows that~$\mu P_t\Rightarrow\pi$ for~$\mu\in\mathcal M_-$.
\end{proof}

We now verify the claim in Remark~\ref{rem:m-minus-examples} for point
masses on configurations with finitely many zeros. The pure-death component gives
\begin{equation*}
\lambda_i(\mb 0)\bb{1-\rho_i(\mb 0)}
\geq
\varepsilon,
\end{equation*}
so~$\rho_i(\mb 0)<1$. By Proposition~\ref{prop:patch-contribution-positivity} we have
$p_i^\star\leq\rho_i(\mb 0)$, and therefore~$p_i^\star<1$ for every~$i$.
Let~$\eta$ have finite zero set
\begin{equation*}
F
=
\cb{i:\eta(i)=0}.
\end{equation*}
For every finite~$A$,
\begin{equation*}
\frac{
\abs{\delta_\eta(\chi_A^\star)}
}{
\mu_{\mb 1}(\chi_A^\star)
}
=
\prod_{i\in A\cap F}
\frac{p_i^\star}{1-p_i^\star}.
\end{equation*}
Consequently, any
\begin{equation*}
K
\geq
\prod_{i\in F}
\max\cb{
1,
\frac{p_i^\star}{1-p_i^\star}
}
\end{equation*}
satisfies
\begin{equation*}
\delta_\eta(\chi_A^\star)
+
K\mu_{\mb 1}(\chi_A^\star)
\geq0
\qquad
\text{for every }A\Subset\Lambda.
\end{equation*}
Hence
\begin{equation*}
\frac{\delta_\eta+K\mu_{\mb 1}}{1+K}
\in
\mathcal M_\star,
\end{equation*}
so~$\delta_\eta\in\mathcal M_{-,K}\subseteq\mathcal M_-$.

\subsection{\texorpdfstring{{Proof of the late-interaction estimate}, Lemma~\ref{lem:late-interactions}}{Proof of the late-interaction estimate}}
\label{subsec:proof-pure-death-comparison}

 To prove Lemma~\ref{lem:late-interactions} we compare the system with the one obtained
by removing its pure deaths. The comparison is simple at the level of the
dual: a pure death at site~$i$ acts on the duality function~$\sigma\chi_A$
through the factor~$-\ind\bb{i\in A}$, so pure deaths change the
potential~$V$ and nothing else. The two systems therefore share the dual
process, the successful interactions and the patch laws, and they differ only
in the patch contributions.
Recall the bulk and end patch families~$\mathcal B_t$
and~$\mathcal E_t$ from Definition~\ref{def:bulk-end-patches}.

\begin{proposition}
\label{prop:death-removal}
Let~$\Lcal$ and~$\Lcal'$ be spin-system generators. Suppose that~$\Lcal$ is
patch-positive and that
\begin{equation}
\label{eq:pure-death-decomposition}
\Lcal f(\eta)
=
\Lcal'f(\eta)
+
\sum_{i\in\Lambda}d_i
\bb{f\bb{\eta^{i,0}}-f(\eta)},
\qquad d_i\geq0.
\end{equation}
Write~$C'$ for the patch contributions of~$\Lcal'$. Then:
\begin{enumerate}[label=\textup{(\roman*)}]
\item The two systems have the same dual process, successful-interaction
skeleton, patch families and patch laws
(cf.~Definition~\ref{def:patch-law}).
\item $\Lcal'$ is patch-positive with the same threshold profile~$\mb p^\star$.
\item Every patch~$P\in\mathcal P$ satisfies~$0\leq C(P)\leq C'(P)$, and every patch~$P$
ending with an outgoing interaction satisfies
$C(P)=e^{-d_{i(P)}(e(P)-s(P))}C'(P)$.
\item For every~$\mu\in\mathcal M_\star$ and every~$t\geq0$, almost surely,
\begin{equation}
\label{eq:death-removal-weight-comparison}
0
\leq
\prod_{P\in\mathcal B_t}C(P)\,
\mu\bb{
\prod_{P\in\mathcal E_t}
C\bb{\eta\bb{i(P)},P}
}
\leq
\prod_{P\in\mathcal B_t}C'(P)\,
\mu\bb{
\prod_{P\in\mathcal E_t}
C'\bb{\eta\bb{i(P)},P}
}.
\end{equation}
\end{enumerate}
\end{proposition}

\begin{proof}[Proof of Proposition~\ref{prop:death-removal}]
Let~$H((A,\sigma),\eta)=\sigma\chi_A(\eta)$ be the duality function from
Section~\ref{sec:dual}. The pure-death term in
\eqref{eq:pure-death-decomposition} acts diagonally on~$H$:
\begin{equation*}
\sum_{i\in\Lambda}d_i\bb{E_{i,0}-I}H((A,\sigma),\eta)
=
-\bb{\sum_{i\in A}d_i}H((A,\sigma),\eta).
\end{equation*}
Hence, if
\begin{equation*}
\Lcal H
=
\Dcal H+V(A)H,
\end{equation*}
then
\begin{equation*}
\Lcal'H
=
\Dcal H
+
\bb{V(A)+\sum_{i\in A}d_i}H.
\end{equation*}
Thus the two systems have the same signed dual generator and therefore the same
dual process, successful-interaction skeleton, patch families, and
patch laws. Their potentials satisfy
\begin{equation}
\label{eq:pure-death-potential-comparison}
V_i'=V_i+d_i.
\end{equation}

We next compare the patch contributions between the two systems. Writing~$\psi_i'$ for the equivalent of~$\psi$ associated with~$\Lcal'$, we
have, by~\eqref{eq:dual-to-spin-contribution-data},
\begin{align*}
\psi_i(\Delta,z)
&=
\delta_i(\varnothing)
\int_0^\Delta
e^{-\bb{\lambda_i'(\mb 0)+d_i}u}\,du
+
z e^{-\bb{\lambda_i'(\mb 0)+d_i}\Delta},
\\
\psi_i'(\Delta,z)
&=
\delta_i(\varnothing)
\int_0^\Delta
e^{-\lambda_i'(\mb 0)u}\,du
+
z e^{-\lambda_i'(\mb 0)\Delta}.
\end{align*}
Hence
\begin{equation}
\label{eq:pure-death-psi-comparison}
0\leq\psi_i(\Delta,z)\leq\psi_i'(\Delta,z)
\qquad
\text{for }\Delta\in[0,\infty]\text{ and }z\in[0,1].
\end{equation}

First consider a patch with incoming initial interaction, i.e. of type~$\mathsf{II}$,~$\mathsf{IE}$ or~$\mathsf{IO}$. The first two
rows of~\eqref{eq:full-contribution-formulas},
together with~\eqref{eq:pure-death-potential-comparison} and
\eqref{eq:pure-death-psi-comparison}, give
\begin{equation*}
0\leq C(P)\leq C'(P).
\end{equation*}
If a patch starts with an outgoing interaction, i.e. is of type~$\mathsf{OI}$,~$\mathsf{OE}$ or~$\mathsf{OO}$, and the interaction has target~$S$, then patch positivity gives
\begin{equation*}
\beta_i(S)\sigma_i^\beta(S)
=
-\widehat{\lambda_i}(S)
\geq0.
\end{equation*}
The remaining rows of~\eqref{eq:full-contribution-formulas} therefore give the
same inequality. For the infinite patches, of types~$\mathsf{IE}$
and~$\mathsf{OE}$, take~$\Delta\to\infty$ in the same formulas. The
denominators are the same for both systems and positive almost surely on
realized infinite patches, since a zero consistency probability gives zero
probability of an infinite continuation at that site, and the coefficient
of~$\psi_i$ is nonnegative; so the comparison holds for them as well. We have proved
\begin{equation}
\label{eq:pure-death-bulk-patch-comparison}
0\leq C(P)\leq C'(P)
\qquad
\text{for every }P\in\mathcal P.
\end{equation}
Moreover, if~$P$ ends with an outgoing interaction, i.e. is of type~$\mathsf{IO}$ or~$\mathsf{OO}$, then
\eqref{eq:full-contribution-formulas} gives $C(P)=e^{V_i\Delta}$ for
type~$\mathsf{IO}$ and $C(P)=\sigma_i^\beta(S)e^{V_i\Delta}$ for
type~$\mathsf{OO}$, where $i=i(P)$ and $\Delta=e(P)-s(P)$, and the sign is the
same for both systems. Thus
\begin{equation}
\label{eq:pure-death-outgoing-terminal-factor}
C(P)
=
e^{-d_{i(P)}\bb{e(P)-s(P)}}C'(P).
\end{equation}

We next prove~(ii). Since~$C'(P)\geq C(P)\geq0$ for every~$P\in\mathcal P$, we conclude
from Proposition~\ref{prop:patch-contribution-positivity}(i)
that~$\Lcal'$ is patch-positive. The pure-death component changes only the constant
coefficient of~$\lambda_i$: for~$S\neq\vn$,
$\widehat{\lambda_i'}(S)=\widehat{\lambda_i}(S)$,
while~$\widehat{\lambda_i'\rho_i'}=\widehat{\lambda_i\rho_i}$. The
definition~\eqref{eq:threshold-profile} of~$\mb p^\star$ involves only these
coefficients, so the two threshold profiles agree.

Now we compare the end patch contributions. For an end patch~$P$ based at
site~$i$, set
\begin{equation*}
\kappa(P)=\partial_zC(z,P),
\qquad
\kappa'(P)=\partial_zC'(z,P).
\end{equation*}
Since~$\varphi_i$ is unchanged, the formulas in
\eqref{eq:end-contribution-formulas} give
\begin{equation}
\label{eq:pure-death-end-patch-comparison}
0\leq\kappa(P)\leq\kappa'(P),
\qquad
0\leq C\bb{p_i^\star,P}
\leq
C'\bb{p_i^\star,P}.
\end{equation}
The second inequality also uses Proposition~\ref{prop:patch-contribution-positivity}
for the two patch-positive systems.

Since~$C(z,P)$ is affine in~$z$,
\begin{equation*}
C\bb{\eta(i),P}
=
C\bb{p_i^\star,P}
+
\kappa(P)\bb{\eta(i)-p_i^\star}.
\end{equation*}
The expansion~\eqref{eq:centered-end-expansion} holds for both systems.
If~$\mu\in\mathcal M_\star$, every centered moment in
\eqref{eq:centered-end-expansion} is nonnegative, so
\eqref{eq:pure-death-end-patch-comparison} yields
\begin{equation}
\label{eq:integrated-end-patch-comparison}
0
\leq
\mu\bb{
\prod_{P\in\mathcal E_t}
C\bb{\eta\bb{i(P)},P}
}
\leq
\mu\bb{
\prod_{P\in\mathcal E_t}
C'\bb{\eta\bb{i(P)},P}
}.
\end{equation}
Combining this with~\eqref{eq:pure-death-bulk-patch-comparison} proves
\eqref{eq:death-removal-weight-comparison}.
\end{proof}

Now consider the special case of Proposition~\ref{prop:death-removal}, i.e. system~$\mathcal{L}$ with a pure-death component of rate~$\varepsilon$. Remove
it from~$\Lcal$ and let
\begin{equation}
\label{eq:uniform-death-decomposition}
\Lcal
=
\Lcal^\varepsilon
+
\varepsilon\Ncal^{\mb 0}.
\end{equation}
Let~$C^\varepsilon$ denote the patch contributions
of~$\Lcal^\varepsilon$, and define
\begin{align*}
W_t^{\varepsilon,\mu}
&=
\prod_{P\in\mathcal B_t}C^\varepsilon(P)\,
\mu\bb{
\prod_{P\in\mathcal E_t}
C^\varepsilon\bb{\eta\bb{i(P)},P}
},
\\
W^\varepsilon
&=
\prod_{P\in\mathcal P}C^\varepsilon(P)
\ind\bb{\abs{\mathcal P}<\infty}.
\end{align*}
We claim that these dominate~$W_t^\mu$ and~$W$, but they carry total mass at
most one.

\begin{corollary}[Death-removed contributions]
\label{cor:death-removed-weights}
Let~$A\Subset\Lambda$,~$\mu\in\mathcal M_\star$, and~$t\geq0$. Then,
almost surely,
\begin{equation}
\label{eq:death-removed-weight-comparison}
0\leq W_t^\mu\leq W_t^{\varepsilon,\mu},
\qquad
\E_A\Cb{W_t^{\varepsilon,\mu}}\leq1,
\end{equation}
and
\begin{equation}
\label{eq:limiting-death-removed-mass}
0\leq W\leq W^\varepsilon,
\qquad
\E_A\Cb{W^\varepsilon}
\leq1.
\end{equation}
\end{corollary}

\begin{proof}
By Proposition~\ref{prop:death-removal}, applied with~$d_i=\varepsilon$, the
two systems have the same dual and patch family,~$\Lcal^\varepsilon$ is
patch-positive, and the inequalities~$0\leq W_t^\mu\leq W_t^{\varepsilon,\mu}$ and~$0\leq W\leq W^\varepsilon$ are an immediate consequence of~\ref{prop:death-removal}(iii) and~\ref{prop:death-removal}(iv).

Now we prove the~$L^1$ bounds. For every finite~$R\supseteq A$, the finite-volume patch representation
for the death-removed system gives
\begin{equation*}
\E_A\Cb{W_t^{\varepsilon,\mu}\ind\bb{E_t^R}} =\mu\bb{P_t^{\varepsilon,R,0}}\chi_A \leq1.
\end{equation*} By
Lemma~\ref{lem:dual-local-finiteness}, the events~$E_t^R$ increase to an
event of probability one as~$R\uparrow\Lambda$. Since~$W_t^{\varepsilon,\mu}\geq0$,
monotone convergence gives~$\E_A\Cb{W_t^{\varepsilon,\mu}}\leq1$.

For the infinite time-horizon variable~$W^\varepsilon$, consider the variable~$W^{\varepsilon, \mu_{\mb 1}}$. On~$\cb{\abs{\mathcal P}<\infty}$,
there are finitely many successful interactions, and after the last one the patches no longer change, except the end patches get longer. On each end patch~$P\in\mathcal{E}_{\infty}$ the
contribution~$C^\varepsilon\bb{1,P^{\downarrow t}}$ converges
to~$C^\varepsilon(P)$. Hence, if we let~$\tau$ be the time of last successful interaction, then on~$\cb{\abs{\mathcal{P}}<\infty}$ and for~$t\geq \tau$ we have
\begin{align*}
W_t^{\varepsilon,\mu_{\mb 1}}
=
\prod_{P\in\mathcal B_{\tau}}C^\varepsilon(P)
\prod_{P\in\mathcal E_{\tau}}
C^\varepsilon\bb{1,P^{\downarrow t}}
\xrightarrow{t\to\infty}
\prod_{P\in\mathcal B_{\tau}}C^\varepsilon(P)
\prod_{P\in\mathcal E_{\tau}}
C^\varepsilon(P)
=W^\varepsilon,
\end{align*}
while~$W^\varepsilon=0$ on~$\cb{\abs{\mathcal P}=\infty}$. Hence
\begin{equation*}
W^\varepsilon
\leq
\liminf_{t\to\infty}
W_t^{\varepsilon,\mu_{\mb 1}},
\end{equation*}
and Fatou's lemma together with~$\E_A\Cb{W_t^{\varepsilon,\mu_{\mb 1}}}\leq1$
gives~$\E_A\Cb{W^\varepsilon}\leq1$.
\end{proof}

The last ingredient is the ancestry of a late interaction, which we
illustrate in Figure~\ref{fig:backward-outgoing-chain}.

\begin{lemma}[Ancestry of a successful interaction]
\label{lem:ancestry-chain}
Almost surely, for every successful interaction at a time~$u>0$ there are
patches~$P_1,\ldots,P_n\in\mathcal P$, each ending with an outgoing
interaction, such that
\begin{equation*}
s(P_1)=0,
\qquad
e(P_n)=u,
\qquad
e(P_k)=s(P_{k+1})
\quad\text{for }1\leq k<n.
\end{equation*}
In particular,
\begin{equation}
\label{eq:outgoing-chain-total-length}
\sum_{k=1}^n\bb{e(P_k)-s(P_k)}
=
u.
\end{equation}
\end{lemma}

\begin{proof}
The source of the interaction at time~$u$ lies in a patch ending with an
outgoing interaction. If that patch starts at time~$0$, it starts at the
formal incoming interaction, and the construction stops. Otherwise it
starts at an earlier successful interaction; move to the patch that ends with
the outgoing side of that interaction, at its source, and repeat. The
interaction times decrease strictly along the way, and only finitely many
successful interactions occur before time~$u$ by
Lemma~\ref{lem:dual-local-finiteness}. Hence the construction stops after
finitely many steps at a patch starting at time~$0$. Listing the patches in
chronological order gives~$P_1,\ldots,P_n$, and the lifetimes
telescope to~\eqref{eq:outgoing-chain-total-length}.
\end{proof}

\begin{proof}[Proof of Lemma~\ref{lem:late-interactions}]
On~$L_{T,t}^c$, let~$u\in(T,t]$ be the time of the first successful
interaction after~$T$, and let~$P_1,\ldots,P_n$ be its ancestry chain from
Lemma~\ref{lem:ancestry-chain}. Every patch in the chain belongs
to~$\mathcal B_t$. Applying~\eqref{eq:pure-death-outgoing-terminal-factor}
with~$d_i=\varepsilon$ and using~\eqref{eq:outgoing-chain-total-length},
\begin{equation*}
\prod_{k=1}^n C(P_k)
=
e^{-\varepsilon u}
\prod_{k=1}^n C^\varepsilon(P_k)
\leq
e^{-\varepsilon T}
\prod_{k=1}^n C^\varepsilon(P_k).
\end{equation*}
The remaining bulk patches satisfy~$0\leq C(P)\leq C^\varepsilon(P)$
by~\eqref{eq:pure-death-bulk-patch-comparison}, and the end-patch factor
of~$W_t^\mu$ is nonnegative and at most the corresponding factor
of~$W_t^{\varepsilon,\mu}$ by~\eqref{eq:integrated-end-patch-comparison}
with~$d_i=\varepsilon$. Hence
\begin{equation*}
W_t^\mu\ind\bb{E_t^R\cap L_{T,t}^c}
\leq
e^{-\varepsilon T}W_t^{\varepsilon,\mu}\ind\bb{E_t^R}.
\end{equation*}
Taking expectations and using
\eqref{eq:death-removed-weight-comparison} proves
\eqref{eq:finite-late-interaction-bound}.

For the infinite-horizon estimate the proof is analogous: apply the same ancestry comparison to the first successful interaction
after~$T$. The remaining patches satisfy $0\leq C(P)\leq C^\varepsilon(P)$
by Proposition~\ref{prop:death-removal}(iii). This gives
\begin{equation*}
W\ind\bb{L_T^c}
\leq e^{-\varepsilon T}W^\varepsilon.
\end{equation*}
Taking
expectations and using~\eqref{eq:limiting-death-removed-mass} proves
\eqref{eq:limiting-late-interaction-bound}.
\end{proof}

\subsection{\texorpdfstring{{Proof of the no-late-interactions estimate}, Lemma~\ref{lem:no-late-relaxation}}{Proof of the no-late-interactions estimate}}
\label{subsec:terminal-relaxation}

If no successful interaction occurs in~$(T,t]$, then the bulk patches
present at time~$T$ are unchanged and every end patch present at time~$T$
extends to time~$t$ (Definition~\ref{def:bulk-end-patches}).
The dependence of a single end patch contribution on~$\mu$ decays exponentially with the patch length. Thus we expect~$W_t^\mu$ and~$W$ to be very close both on~$L_{T,t}$ and on~$L_T$. We also expect both of these events to be almost the same, since the most likely reason for the dual process having no successful interactions in a long interval~$(T,t]$ is that the process died before or shortly after~$T$, in which case it continues to not have successful interactions after time~$t$.

The trick to prove~Lemma~\ref{lem:no-late-relaxation} is that both~$\CE{W_t^\mu \ind\bb{L_{T,t}}}{\mathcal{G}_T}$ and~$\CE{W\ind\bb{L_T}}{\mathcal{G}_T}$ factorize over the same bulk patches and end patches, with factors exponentially close together. Once this factorization is shown a telescoping identity proves the result. To attain the factorization we define, for each end patch~$P\in\mathcal{E}_T$ the events
\begin{align*}
    L_{T,t}^P &= \cb{\text{no successful interaction occurs at site }i(P)\text{ during }(T,t]},\\
    L_{T}^P &= \cb{\text{no successful interaction occurs at site }i(P)\text{ during }(T,\infty)}.
\end{align*}

\begin{lemma}
\label{prop:no-late-continuation-identity}
Fix~$0\leq T<t$.
\begin{enumerate}[label=\textup{(\roman*)}]
\item The events of no further successful interaction factor over the end patches:
\begin{equation}
\label{eq:local-no-late-factorization}
\begin{alignedat}{2}
L_{T,t}
  &= \bigcap_{P\in\mathcal E_T} L_{T,t}^P,
&\qquad L_T
  &= \bigcap_{P\in\mathcal E_T} L_T^P,\\
\CP{L_{T,t}}{\mathcal G_T}
  &= \prod_{P\in\mathcal E_T}\CP{L_{T,t}^P}{\mathcal G_T},
&\qquad \CP{L_T}{\mathcal G_T}
  &= \prod_{P\in\mathcal E_T}\CP{L_T^P}{\mathcal G_T}.
\end{alignedat}
\end{equation}
\item On~$L_{T,t}$,
\begin{equation*}
\mathcal B_t=\mathcal B_T,
\qquad
\mathcal E_t
=
\cb{P^{\uparrow t}:P\in\mathcal E_T}.
\end{equation*}
\item For every~$\mu\in\mathcal M_\star$,
\begin{align}
\CE{W_t^\mu\ind\bb{L_{T,t}}}{\mathcal G_T}
&=
\prod_{P\in\mathcal B_T}C(P)\,
\CP{L_{T,t}}{\mathcal G_T}
\mu\bb{
 \prod_{P\in\mathcal E_T}
 C\bb{\eta\bb{i(P)},P^{\uparrow t}}
}
\notag
\\
&=
\prod_{P\in\mathcal B_T}C(P)\,
\mu\bb{
 \prod_{P\in\mathcal E_T}
 C\bb{\psi_{i(P)}\bb{t-T,\eta\bb{i(P)}},P}
}.
\label{eq:appendix-finite-continuation-identity}
\end{align}
\item Furthermore,
\begin{equation}
\label{eq:appendix-infinite-continuation-identity}
\CE{W\ind\bb{L_T}}{\mathcal G_T}
=
\prod_{P\in\mathcal B_T}C(P)
\prod_{P\in\mathcal E_T}
C\bb{\rho_{i(P)}(\mb 0),P}.
\end{equation}
\end{enumerate}
The first equality in~\eqref{eq:appendix-finite-continuation-identity} also
holds for~$W_t^{\varepsilon,\mu}$ with~$C$ replaced by~$C^\varepsilon$.
\end{lemma}

\begin{proof}[Proof of Lemma~\ref{prop:no-late-continuation-identity}]
We first prove~(i). Before the first successful interaction after~$T$, the
only changes in activity are dual deaths. Every site active at time~$T$
belongs to an end patch, so the first interaction at an active site is the
same in the full process and in the one-site continuations. This proves the
identity for~$L_{T,t}$; intersecting over~$t=T+n$,~$n\geq1$, gives the
identity for~$L_T$. Conditional on~$\mathcal G_T$, the patch variables have
independent consistent laws, and their future Poisson marks are independent
at distinct sites. This proves the conditional probability identities.

For~(ii), note that on~$L_{T,t}$ no successful interaction starts a new
patch after time~$T$. Thus the bulk patches are unchanged, while
each~$P\in\mathcal E_T$ becomes its geometric extension~$P^{\uparrow t}$ from
\eqref{eq:geometric-patch-extension}.

We now prove~(iii). The first equality is just~$\mathcal{G}_T$-measurability of~$\mathcal{B}_T$ and~$\mathcal{E}_T$ and the fact that on event~$L_{T,t}$ the bulk patches don't change between time~$T$ and~$t$ and the end patches are just elongated. The system with the pure-death component removed has the same signed dual,
successful-interaction skeleton, and patch laws, by
Proposition~\ref{prop:death-removal}(i). Hence the event~$L_{T,t}$ and its
conditional probability are unchanged, which proves the corresponding first
equality for~$W_t^{\varepsilon,\mu}$.

We prove the second equality. Fix~$P\in\mathcal E_T$ and write
\begin{equation*}
i=i(P),
\qquad
S=S(P),
\qquad
\Delta=T-s(P).
\end{equation*}
The patch~$P^{\uparrow t}$ has length~$\Delta+t-T$. Therefore
\eqref{eq:end-contribution-formulas} gives
\begin{equation}
\label{eq:end-patch-continuation}
C\bb{z,P^{\uparrow t}}
=
\pw{
\displaystyle
\frac{\psi_i(\Delta+t-T,z)}{\varphi_i(\Delta+t-T)}
}{\mathsf X(P)=\mathsf I}
{
\displaystyle
\frac{
\delta_i(S)\sigma_i^\delta(S)
+
\beta_i(S)\sigma_i^\beta(S)\psi_i(\Delta+t-T,z)
}{
\delta_i(S)+\beta_i(S)\varphi_i(\Delta+t-T)
}
}{\mathsf X(P)=\mathsf O}.
\end{equation}
Here we apply the end-patch formula from Lemma~\ref{lem:contribution-formulas}
to the longer patch.

We combine the consistency probabilities from
Lemma~\ref{lem:contribution-formulas} with the patch factorization of
Lemma~\ref{thm:patch-factorization} to compute the probability that the
patch remains consistent during the added interval. Conditional on~$\mathcal G_T$,
\begin{equation}
\label{eq:local-continuation-probability}
\CP{L_{T,t}^P}{\mathcal G_T}
=
\pw{
\displaystyle
\frac{\varphi_i(\Delta+t-T)}{\varphi_i(\Delta)}
}{\mathsf X(P)=\mathsf I}
{
\displaystyle
\frac{
\delta_i(S)+\beta_i(S)\varphi_i(\Delta+t-T)
}{
\delta_i(S)+\beta_i(S)\varphi_i(\Delta)
}
}{\mathsf X(P)=\mathsf O}.
\end{equation}
Multiplying~\eqref{eq:end-patch-continuation} by
\eqref{eq:local-continuation-probability} and using
\begin{equation*}
\psi_i(\Delta+t-T,z)
=
\psi_i\bb{\Delta,\psi_i(t-T,z)}
\end{equation*}
from~\eqref{eq:all-zero-neighbourhood-relaxation}, we obtain
\begin{equation}
\label{eq:local-probability-weighted-continuation}
\CP{L_{T,t}^P}{\mathcal G_T}
C\bb{z,P^{\uparrow t}}
=
C\bb{\psi_i(t-T,z),P}.
\end{equation}
This proves~(iii). 

For~(iv), it remains to pass from a finite continuation to~$L_T$. By the pure-death
assumption,~$\lambda_i(\mb 0)\geq\varepsilon$, and therefore
\eqref{eq:all-zero-neighbourhood-relaxation} gives
\begin{equation*}
\psi_i(t-T,1)\longrightarrow\rho_i(\mb 0)
\qquad\text{as }t\to\infty.
\end{equation*}
For each~$P\in\mathcal E_T$,
\eqref{eq:local-probability-weighted-continuation} with~$z=1$ yields
\begin{equation}
\label{eq:local-infinite-continuation-limit}
\lim_{t\to\infty}
\CP{L_{T,t}^P}{\mathcal G_T}
C\bb{1,P^{\uparrow t}}
=
C\bb{\rho_i(\mb 0),P}.
\end{equation}
Finally, we have
\begin{align*}
    \CE{W\ind\bb{L_T}}{\mathcal{G}_T} &= \prod_{P\in\mathcal{B}_T}C(P)\prod_{P\in\mathcal{E}_T} \CP{L_{T}^P}{\mathcal G_T}C\bb{1,P^{\uparrow \infty}}\\
    &=\prod_{P\in\mathcal{B}_T}C(P)\prod_{P\in\mathcal{E}_T}\lim_{t\to\infty}
\CP{L_{T,t}^P}{\mathcal G_T}
C\bb{1,P^{\uparrow t}}\\
    &= \prod_{P\in\mathcal{B}_T}C(P)\prod_{P\in\mathcal{E}_T}C\bb{\rho_i(\mb 0),P}.
\end{align*}
The first equality follows by the same~$\mathcal{G}_T$-measurability argument as in~(iii) and the last equality in~(i). The second equality holds because~$L_T^P = \bigcap_{t\geq T}L_{T,t}^P$ and by continuity of measure and~$C$ with respect to patch length. The last equality follows by~\eqref{eq:local-infinite-continuation-limit} and the fact that~$\mathcal{E}_T$ is almost surely finite.
\end{proof}

In Lemma~\ref{prop:no-late-continuation-identity}, we showed that~$\CE{W_t^\mu \ind\bb{L_{T,t}}}{\mathcal{G}_T}$ and~$\CE{W \ind\bb{L_{T}}}{\mathcal{G}_T}$ factorize with the same factors over~$\mathcal{B}_T$ and with convergent factors~$\mathcal{E}_T$, which already gives us almost sure convergence
\begin{equation*}
\CE{W_t^\mu\ind\bb{L_{T,t}}}{\mathcal G_T}
\longrightarrow
\CE{W\ind\bb{L_T}}{\mathcal G_T}.
\end{equation*}
In the next lemma we upgrade this to an explicit upper bound, exponentially decreasing in~$t-T$.

\begin{lemma}[End-factor relaxation]
\label{lem:end-factor-relaxation}
Let~$0\leq T<t$. Almost surely, for every~$\mu\in\mathcal M_\star$, we have
\begin{align*}
&\abs{
 \CP{L_{T,t}}{\mathcal G_T}
 \mu\bb{
  \prod_{P\in\mathcal E_T}
  C\bb{\eta\bb{i(P)},P^{\uparrow t}}
 }
 -
 \prod_{P\in\mathcal E_T}
 C\bb{\rho_{i(P)}(\mb 0),P}
}
\\
&\qquad\leq
e^{-\varepsilon(t-T)}
\bb{
 \CP{L_{T,t}}{\mathcal G_T}
 \mu\bb{
  \prod_{P\in\mathcal E_T}
  C^\varepsilon\bb{\eta\bb{i(P)},P^{\uparrow t}}
 }
 +
 \abs{\mathcal E_T}
 \prod_{P\in\mathcal E_T}
 C\bb{\rho_{i(P)}(\mb 0),P}
}.
\end{align*}
\end{lemma}

\begin{proof}[Proof of Lemma~\ref{lem:end-factor-relaxation}]
For an end patch~$P$ based at~$i$, define
\begin{align*}
a_t(P)
&=
\CP{L_{T,t}^P}{\mathcal G_T}
C\bb{p_i^\star,P^{\uparrow t}},
&
k_t(P)
&=
\CP{L_{T,t}^P}{\mathcal G_T}
\partial_z C\bb{z,P^{\uparrow t}},
\\
a_t^\varepsilon(P)
&=
\CP{L_{T,t}^P}{\mathcal G_T}
C^\varepsilon\bb{p_i^\star,P^{\uparrow t}},
&
k_t^\varepsilon(P)
&=
\CP{L_{T,t}^P}{\mathcal G_T}
\partial_z C^\varepsilon\bb{z,P^{\uparrow t}},
\\
b(P)
&=
C\bb{\rho_i(\mb 0),P}.
\end{align*}
Also let
\begin{equation*}
\kappa(P)=\partial_zC(z,P),
\qquad
\kappa^\varepsilon(P)=\partial_zC^\varepsilon(z,P).
\end{equation*}
Since the end-patch contributions are affine in~$z$, neither~$\kappa(P)$ nor~$\kappa^\varepsilon(P)$
depends on~$z$.

By~\eqref{eq:local-probability-weighted-continuation},
\begin{equation}
\label{eq:continued-end-factor-data}
a_t(P)
=
C\bb{\psi_i(t-T,p_i^\star),P},
\qquad
k_t(P)
=
e^{-\lambda_i(\mb 0)(t-T)}\kappa(P).
\end{equation}
For the system with the pure-death component removed, the update rate at the
all-zero neighbourhood is~$\lambda_i(\mb 0)-\varepsilon$. The same calculation
therefore gives
\begin{equation*}
k_t^\varepsilon(P)
=
e^{-\bb{\lambda_i(\mb 0)-\varepsilon}(t-T)}
\kappa^\varepsilon(P).
\end{equation*}
We apply the end-patch comparison~\eqref{eq:pure-death-end-patch-comparison}
from Proposition~\ref{prop:death-removal} to~$P$ and its
extension~$P^{\uparrow t}$ and obtain
\begin{equation*}
0\leq a_t(P)\leq a_t^\varepsilon(P),
\qquad
0\leq\kappa(P)\leq\kappa^\varepsilon(P).
\end{equation*}
Consequently,
\begin{equation}
\label{eq:appendix-continued-factor-comparisons}
0\leq a_t(P)\leq a_t^\varepsilon(P),
\qquad
0\leq k_t(P)\leq e^{-\varepsilon(t-T)}k_t^\varepsilon(P).
\end{equation}

We also compare~$a_t(P)$ with its limiting value~$b(P)$. Since~$p_i^\star\leq\rho_i(\mb 0)$ and
\begin{equation*}
\psi_i(t-T,p_i^\star)
=
\rho_i(\mb 0)
-
\bb{\rho_i(\mb 0)-p_i^\star}
e^{-\lambda_i(\mb 0)(t-T)},
\end{equation*}
the affine dependence of~$C(z,P)$ on~$z$ gives
\begin{equation*}
b(P)-a_t(P)
=
\kappa(P)\bb{\rho_i(\mb 0)-p_i^\star}
e^{-\lambda_i(\mb 0)(t-T)}.
\end{equation*}
Moreover,
\begin{equation*}
b(P)
=
C\bb{p_i^\star,P}
+
\kappa(P)\bb{\rho_i(\mb 0)-p_i^\star}
\geq
\kappa(P)\bb{\rho_i(\mb 0)-p_i^\star},
\end{equation*}
where~$C(p_i^\star,P)\geq0$ by patch positivity. Since~$\lambda_i(\mb 0)\geq\varepsilon$,
\begin{equation}
\label{eq:continued-end-factor-limit-comparison}
0
\leq
b(P)-a_t(P)
\leq
e^{-\varepsilon(t-T)}b(P).
\end{equation}

We first control the terms that depend on centered moments of~$\mu$. For~$\mathcal Q\subseteq\mathcal E_T$, set
\begin{equation*}
I(\mathcal Q)
=
\cb{i(P):P\in\mathcal Q}.
\end{equation*}
The end-patch sites are distinct. Using
\eqref{eq:local-no-late-factorization} and expanding each factor around~$p_i^\star$ gives
\begin{equation}
\CP{L_{T,t}}{\mathcal G_T}
\mu\bb{
 \prod_{P\in\mathcal E_T}
 C\bb{\eta\bb{i(P)},P^{\uparrow t}}
}
=
\sum_{\mathcal Q\subseteq\mathcal E_T}
\mu\bb{\chi_{I(\mathcal Q)}^\star}
\prod_{P\in\mathcal Q}k_t(P)
\prod_{P\in\mathcal E_T\setminus\mathcal Q}a_t(P).
\label{eq:continued-end-factor-centered-expansion}
\end{equation}
Every term is nonnegative because~$\mu\in\mathcal M_\star$. The term with~$\mathcal Q=\varnothing$ is~$\prod_Pa_t(P)$. If~$\mathcal Q\neq\varnothing$, the inequalities in
\eqref{eq:appendix-continued-factor-comparisons} give
\begin{equation*}
\mu\bb{\chi_{I(\mathcal Q)}^\star}
\prod_{P\in\mathcal Q}k_t(P)
\prod_{P\in\mathcal E_T\setminus\mathcal Q}a_t(P)
\leq
e^{-\varepsilon(t-T)}
\mu\bb{\chi_{I(\mathcal Q)}^\star}
\prod_{P\in\mathcal Q}k_t^\varepsilon(P)
\prod_{P\in\mathcal E_T\setminus\mathcal Q}a_t^\varepsilon(P).
\end{equation*}
Summing over nonempty~$\mathcal Q$ and then adding the nonnegative empty-set
term on the right gives
\begin{align}
0
&\leq
\CP{L_{T,t}}{\mathcal G_T}
\mu\bb{
 \prod_{P\in\mathcal E_T}
 C\bb{\eta\bb{i(P)},P^{\uparrow t}}
}
-
\prod_{P\in\mathcal E_T}a_t(P)
\notag\\
&\leq
e^{-\varepsilon(t-T)}
\CP{L_{T,t}}{\mathcal G_T}
\mu\bb{
 \prod_{P\in\mathcal E_T}
 C^\varepsilon\bb{\eta\bb{i(P)},P^{\uparrow t}}
}.
\label{eq:continued-end-factor-centered-error}
\end{align}

It remains to replace $\prod_Pa_t(P)$ by $\prod_Pb(P)$. Enumerate
$\mathcal E_T=\cb{P_1,\ldots,P_n}$ and write $a_j=a_t(P_j)$ and $b_j=b(P_j)$.
Since $0\leq a_j\leq b_j$,
\begin{equation*}
\prod_{j=1}^n b_j-\prod_{j=1}^n a_j
=
\sum_{j=1}^n
\bb{b_j-a_j}
\prod_{k<j}a_k
\prod_{k>j}b_k.
\end{equation*}
Using~\eqref{eq:continued-end-factor-limit-comparison},
\begin{equation}
\label{eq:continued-end-factor-constant-error}
0
\leq
\prod_{P\in\mathcal E_T}b(P)
-
\prod_{P\in\mathcal E_T}a_t(P)
\leq
\abs{\mathcal E_T}e^{-\varepsilon(t-T)}
\prod_{P\in\mathcal E_T}b(P).
\end{equation}
Both quantities compared in the lemma are at least~$\prod_{P\in\mathcal E_T}a_t(P)$. Combining
\eqref{eq:continued-end-factor-centered-error} and
\eqref{eq:continued-end-factor-constant-error} proves the lemma.
\end{proof}

\begin{proof}[Proof of Lemma~\ref{lem:no-late-relaxation}]
The first inequality in~\eqref{eq:no-late-relaxation-bound} is proved as follows:
\begin{align*}
&\abs{
  \E_A\Cb{W_t^\mu\ind\bb{E_t^R\cap L_{T,t}}}
  -\E_A\Cb{W\ind\bb{E_T^R\cap L_T}}
}\\
&\quad\leq \E_A\Biggl[
  \prod_{P\in\mathcal B_T} C(P)\,\ind\bb{E_T^R}\\
&\qquad{}\times
  \CE{\abs{
    \CP{L_{T,t}}{\mathcal G_T}\,
    \mu\bb{\prod_{P\in\mathcal E_T} C\bb{\eta(i(P)),P^{\uparrow t}}}
    -\prod_{P\in\mathcal E_T} C\bb{\rho_{i(P)}(\mb 0),P}
  }}{\mathcal G_T}
\Biggr]\\
&\quad\leq e^{-\varepsilon(t-T)}\E_A\Biggl[
  \prod_{P\in\mathcal B_T} C(P)\,\ind\bb{E_T^R}\\
&\qquad{}\times\biggl(
  \CP{L_{T,t}}{\mathcal G_T}\,
  \mu\bb{\prod_{P\in\mathcal E_T} C^\varepsilon\bb{\eta(i(P)),P^{\uparrow t}}}
  +\abs{\mathcal E_T}\prod_{P\in\mathcal E_T} C\bb{\rho_{i(P)}(\mb 0),P}
  \biggr)
\Biggr]\\
&\quad\leq e^{-\varepsilon(t-T)}\E_A\Biggl[
  \ind\bb{E_T^R}\biggl(
  \prod_{P\in\mathcal B_T} C^\varepsilon(P)\,
  \CP{L_{T,t}}{\mathcal G_T}\,
  \mu\bb{\prod_{P\in\mathcal E_T} C^\varepsilon\bb{\eta(i(P)),P^{\uparrow t}}}\\
&\qquad\qquad{}
  +\abs R\prod_{P\in\mathcal B_T} C(P)
  \prod_{P\in\mathcal E_T} C\bb{\rho_{i(P)}(\mb 0),P}
  \biggr)
\Biggr]\\
&\quad=
  e^{-\varepsilon(t-T)}
  \bb{
    \E_A\Cb{W_t^{\varepsilon,\mu}\ind\bb{E_t^R\cap L_{T,t}}}
    +\abs R\,\E_A\Cb{W\ind\bb{E_T^R\cap L_T}}
  }\\
&\quad\leq (1+\abs R)e^{-\varepsilon(t-T)}.
\end{align*}
For the first inequality, we use Lemma~\ref{prop:no-late-continuation-identity}(iii)-(iv) and the triangle inequality. Note that~$E_t^R\cap L_{T,t}=E_T^R\cap L_{T,t}$ and~$E_T^R$ is~$\mathcal G_T$-measurable, as is the bulk product.
For the second inequality, we apply Lemma~\ref{lem:end-factor-relaxation}
inside the conditional expectation and use the fact that the upper bound
is~$\mathcal G_T$-measurable. The two end-patch terms are nonnegative,
as shown in the proof of that lemma. In the third inequality,
we use~$C(P)\leq C^\varepsilon(P)$ for the first term and~$\abs{\mathcal E_T}\leq\abs{R}$ on~$E_T^R$ for the second.

For the equality, for the first term we use the first identity in
Lemma~\ref{prop:no-late-continuation-identity}(iii) with~$C^\varepsilon$, and for the second term we use part~(iv). The last inequality uses $\E_A\Cb{W_t^{\varepsilon,\mu}}\leq1$
and $\E_A\Cb{W}\leq\E_A\Cb{W^\varepsilon}\leq1$, see
Corollary~\ref{cor:death-removed-weights}. This proves the first inequality in~\eqref{eq:no-late-relaxation-bound}.

For the remaining bound, the identity
\eqref{eq:appendix-infinite-continuation-identity} gives
$\CE{W\ind\bb{L_T}}{\mathcal G_T}\leq W_T^{\mu_{\mb 1}}$
because $0\leq C(\rho_{i(P)}(\mb 0),P)\leq C(1,P)$ for every end patch.
For finite~$S\supseteq R$, the event~$E_T^S\setminus E_T^R$
belongs to~$\mathcal G_T$. Conditioning on~$\mathcal G_T$ and using the
finite-volume representation~\eqref{eq:patch-representation} twice gives
\begin{align*}
\E_A\Cb{W\ind\bb{L_T\cap(E_T^S\setminus E_T^R)}}
&\leq\E_A\Cb{W_T^{\mu_{\mb 1}}\ind\bb{E_T^S\setminus E_T^R}}\\
&=\mu_{\mb 1}\bb{\bb{P_T^{S,0}-P_T^{R,0}}\chi_A}.
\end{align*}
As~$S\uparrow\Lambda$, local finiteness and monotone convergence give
$\E_A\Cb{W\ind\bb{(E_T^R)^c\cap L_T}}$ on the left, and
Lemma~\ref{lem:finite-propagation}(b) gives
$\mu_{\mb 1}\bb{\bb{P_T-P_T^{R,0}}\chi_A}$ on the right.
We bound this by~$\norm{\bb{P_T-P_T^{R,0}}\chi_A}_\infty$, which proves~\eqref{eq:no-late-relaxation-bound} and thus Lemma~\ref{lem:no-late-relaxation}.
\end{proof}

\Needspace{9\baselineskip}
\subsection{\texorpdfstring
  {Proof of Corollary~\ref*{cor:uniform-exponential-ergodicity} and Corollary~\ref*{cor:oriented-ergodicity}}
  {Proof of Corollary \getrefnumber{cor:uniform-exponential-ergodicity} and Corollary \getrefnumber{cor:oriented-ergodicity}}}
\label{subsec:ergodicity-corollaries}
We prove Corollary~\ref{cor:uniform-exponential-ergodicity} and Corollary~\ref{cor:oriented-ergodicity}. The former is a special case of Theorem~\ref{thm:common-invariant-limit}. The latter we prove via a burn-in argument.

\subsubsection{Uniform exponential ergodicity when~$\mb p^\star\leq\mbs{\tfrac12}$}
\label{subsec:proof-uniform-exponential}
For~$\mb p^\star\leq\mbs{\tfrac12}$ the family of measures~$\mathcal{M}_{-,1}$
contains all probability measures, so we obtain uniform exponential
ergodicity by applying Theorem~\ref{thm:common-invariant-limit}.
\begin{proof}[Proof of Corollary~\ref{cor:uniform-exponential-ergodicity}]
Suppose~$\mb p^\star\leq\mbs{\tfrac12}$. For every configuration~$\eta$
and every finite~$A$,
\begin{equation*}
\chi_A^\star(\eta)
\geq
-\prod_{i\in A}\bb{1-p_i^\star}
=
-\mu_{\mb 1}(\chi_A^\star),
\end{equation*}
because~$p_i^\star\leq1-p_i^\star$ at every site. Hence, for every
probability measure~$\mu$,
\begin{equation*}
\frac{\mu+\mu_{\mb 1}}2
\bb{\chi_A^\star}
\geq0
\qquad
\text{for every }A\Subset\Lambda.
\end{equation*}
Thus every probability measure belongs to~$\mathcal M_{-,1}$. We apply
Theorem~\ref{thm:common-invariant-limit} with~$K=1$ and obtain
\eqref{eq:uniform-pointwise-convergence} after replacing~$K_f$ by~$3K_f$. If~$\nu$ is an invariant probability measure, integrating
\eqref{eq:uniform-pointwise-convergence} against~$\nu$ gives
$\abs{\nu(f)-\pi(f)}\leq3K_f(1+t)^De^{-\varepsilon t/2}$ for every~$t$, so
$\nu=\pi$.
\end{proof}

\subsubsection{Uniform ergodicity under orientability}
\label{subsec:proof-orientability}

By Theorem~\ref{thm:common-invariant-limit}, every initial law
in~$\mathcal M_\star$ converges to~$\pi$. To prove
Corollary~\ref{cor:oriented-ergodicity}, we show that the evolution of any
initial measure approaches~$\mathcal M_\star$. This is immediate for finite~$\Lambda$, because
then~$\mathcal M_-$ contains all probability measures. The point of
orientability is that, for any finite set~$A$, we can condition on the
spin-system trajectory on the set~$U(A)$ from
Definition~\ref{def:orientability} and get a finite, patch-positive spin system on~$\bar A$. This reduced system has time-inhomogeneous boundary conditions,
so it may not converge to an invariant measure, but it still preserves
nonnegative centered moments on~$\bar A$, with the original centering
profile. This is the content of the next lemma, and it is what allows a
burn-in argument to go through. More precisely, we prove that
\begin{equation*}
\liminf_{t\to\infty}
\inf_{\eta\in\Omega}
P_t\chi_A^\star(\eta)
\geq0
\qquad
\text{for every }A\Subset\Lambda,
\end{equation*}
and then combine this with compactness and
Theorem~\ref{thm:common-invariant-limit}.

\begin{lemma}[{Preservation of~$\mathcal{M}_*$ under boundary conditions}]
\label{lem:boundary-moment-preservation}
{Suppose that the spin system is patch-positive.} Let~$I\Subset\Lambda$ and~$U\subseteq\Lambda\setminus I$ satisfy
$N(i)\subseteq I\cup U$ for every~$i\in I$. For~$\xi\in\cb{0,1}^U$
and~$\eta\in\cb{0,1}^I$, let~$\eta^{U,\xi}\in\Omega$ be any configuration
that agrees with~$\xi$ on~$U$ and with~$\eta$ on~$I$; the rates at the sites
of~$I$ do not depend on the choice. Fixing the spins on~$U$ at~$\xi$ defines
a spin system on~$\cb{0,1}^I$ with rates
\begin{equation*}
\lambda_i^{I,\xi}(\eta)=\lambda_i\bb{\eta^{U,\xi}},
\qquad
\rho_i^{I,\xi}(\eta)=\rho_i\bb{\eta^{U,\xi}}
\qquad(i\in I);
\end{equation*}
let~$P_t^{I,\xi}$ be its semigroup. With the original centering
profile~$(p_i^\star)_{i\in I}$, let
\begin{equation*}
\mathcal M_\star^I
=
\cb{
\mu\in\mathcal P\bb{\cb{0,1}^I}:
\mu(\chi_B^\star)\geq0
\text{ for every }B\subseteq I
}.
\end{equation*}
Then~$\mathcal M_\star^I P_t^{I,\xi}\subseteq\mathcal M_\star^I$ for
every~$t\geq0$.

\end{lemma}

\begin{proof}
We first check that
\begin{equation}
\label{eq:centered-generator-positive-coefficients}
\widehat{\lambda_i\rho_i}(S)-p_i^\star\widehat{\lambda_i}(S)
\geq0
\qquad
\text{for every }S\subseteq N(i).
\end{equation}
For~$S=\vn$, the difference is~$\lambda_i(\mb 0)\bb{\rho_i(\mb 0)-p_i^\star}\geq0$
by~\eqref{eq:threshold-below-zero-neighbour-refresh}. If $S\neq\vn$ and
$\widehat{\lambda_i}(S)<0$, the inequality follows from
\eqref{eq:threshold-profile}. If $S\neq\vn$ and $\widehat{\lambda_i}(S)=0$,
the difference equals
\begin{equation*}
\widehat{\lambda_i\bb{\rho_i-\rho_i(\mb 0)}}(S),
\end{equation*}
which is nonnegative by patch positivity.

For~$S\subseteq N(i)\cap I$,
substituting the fixed spins~$\xi$ into the monomial expansions gives
\begin{align}
\widehat{\lambda_i^{I,\xi}}(S)
&=
\sum_{T\subseteq N(i)\cap U}
\widehat{\lambda_i}(S\cup T)
\prod_{j\in T}\xi(j),
\label{eq:restricted-lambda-coefficients}
\\
\widehat{
\lambda_i^{I,\xi}
\bb{\rho_i^{I,\xi}-p_i^\star}
}(S)
&=
\sum_{T\subseteq N(i)\cap U}
\bb{
\widehat{\lambda_i\rho_i}(S\cup T)
-
p_i^\star\widehat{\lambda_i}(S\cup T)
}
\prod_{j\in T}\xi(j).
\label{eq:restricted-centered-coefficients}
\end{align}
If~$S\neq\varnothing$, then every~$S\cup T$ appearing in
\eqref{eq:restricted-lambda-coefficients} is nonempty, so patch positivity
gives
\begin{equation}
\label{eq:restricted-lambda-sign}
\widehat{\lambda_i^{I,\xi}}(S)\leq0.
\end{equation}
By~\eqref{eq:centered-generator-positive-coefficients}, every term in
\eqref{eq:restricted-centered-coefficients} is nonnegative. Hence
\begin{equation}
\label{eq:restricted-centered-sign}
\widehat{
\lambda_i^{I,\xi}
\bb{\rho_i^{I,\xi}-p_i^\star}
}(S)
\geq0
\qquad
\text{for every }S\subseteq N(i)\cap I.
\end{equation}

For~$B\subseteq I$, an update at~$i\in B$ replaces the centered factor
at~$i$ by~$\rho_i^{I,\xi}-p_i^\star$ in expectation. Expanding the rates
in monomials gives
\begin{align*}
\Lcal^{I,\xi}\chi_B^\star
={}&
-\sum_{i\in B}
\lambda_i^{I,\xi}(\mb 0)\chi_B^\star
\\
&+
\sum_{i\in B}
\sum_{S\subseteq N(i)\cap I}
\widehat{
\lambda_i^{I,\xi}
\bb{\rho_i^{I,\xi}-p_i^\star}
}(S)
\chi_S\chi_{B\setminus\cb{i}}^\star
\\
&+
\sum_{i\in B}
\sum_{\substack{S\subseteq N(i)\cap I\\S\neq\varnothing}}
\bb{-\widehat{\lambda_i^{I,\xi}}(S)}
\chi_S\chi_B^\star.
\end{align*}
Choose~$r>0$ so that
\begin{equation*}
r\sum_{i\in I}\lambda_i^{I,\xi}(\mb 0)\leq1.
\end{equation*}
By~\eqref{eq:restricted-lambda-sign} and
\eqref{eq:restricted-centered-sign}, all coefficients in
$(I+r\Lcal^{I,\xi})\chi_B^\star$ are nonnegative before we expand the
products~$\chi_S\chi_D^\star$. For~$j\in S\cap D$, we use
\begin{equation*}
\eta(j)\bb{\eta(j)-p_j^\star}=\bb{1-p_j^\star}\eta(j).
\end{equation*}
We then write each uncentered factor as~$\eta(j)=p_j^\star+\bb{\eta(j)-p_j^\star}$.
Since~$0\leq p_j^\star\leq1$, the resulting coefficients are nonnegative.
Thus~$(I+r\Lcal^{I,\xi})\chi_B^\star$ is a nonnegative linear combination
of~$\cb{\chi_D^\star:D\subseteq I}$.

Since~$I$ is finite, for every~$t\geq0$,
\begin{equation*}
P_t^{I,\xi}
=
\lim_{n\to\infty}
\bb{
I+\frac{t}{n}\Lcal^{I,\xi}
}^n
\end{equation*}
as operators on functions on~$\cb{0,1}^I$. For all sufficiently large~$n$,
the preceding argument applies with~$r=t/n$. Thus, if~$\mu\in\mathcal M_\star^I$, then for every~$B\subseteq I$,
\begin{equation*}
(\mu P_t^{I,\xi})(\chi_B^\star)
=
\lim_{n\to\infty}
\mu\bb{
\bb{
I+\frac{t}{n}\Lcal^{I,\xi}
}^n
\chi_B^\star
}
\geq0.
\end{equation*}
Hence~$\mu P_t^{I,\xi}\in\mathcal M_\star^I$.
\end{proof}

\begin{proof}[Proof of Corollary~\ref{cor:oriented-ergodicity}]

Fix a nonempty finite set~$A$ and take~$I=\bar A$ and~$U=U(A)$
from Definition~\ref{def:orientability}. Since~$I\cup U$ is autonomous,~$N(i)\subseteq I\cup U$ for every~$i\in I$.
We can therefore apply Lemma~\ref{lem:boundary-moment-preservation} to every
boundary condition~$\xi\in\cb{0,1}^U$, with~$P_t^{I,\xi}$
and~$\mathcal M_\star^I$ as in the lemma.
Fix an initial configuration~$\eta$ and define the~$\sigma$-algebra
\begin{equation*}
\mathcal G_U:=\sigma\cb{\eta_s(j):j\in U,\ s\geq0}.
\end{equation*}
Since both~$U$ and~$I\cup U$ are autonomous,
conditional on~$\mathcal G_U$, the process on~$I$ evolves with
the revealed trajectory on $U$ as a boundary condition. Only the states on the boundary $\partial I:=N_*(I)\setminus I\subseteq U$ matter.
This set is finite and the update rates are bounded, so the boundary conditions change finitely many times before time~$t$.
If the successive boundary configurations are~$\xi_1,\ldots,\xi_n$,
extended arbitrarily from~$\partial I$ to~$U$, and their durations
are~$t_1,\ldots,t_n$ with~$\sum_{k=1}^n t_k=t$, then, for~$f$ supported on~$I$,
\begin{equation}
\label{eq:revealed-boundary-product}
\CE{f(\eta_t)}{\mathcal G_U}
=P_{t_1}^{I,\xi_1}\cdots P_{t_n}^{I,\xi_n}(f)(\eta).
\end{equation}

\enlargethispage{-0.5\baselineskip}
We next show that, for every initial configuration~$\eta$ and boundary conditions, after a unit of time a fixed positive part of the evolution law enters~$\mathcal{M}_\star^I$. To be more precise we obtain a decomposition
\begin{equation}\label{eq:mass_leak_decomposition}
    \delta_\eta P_{s_1}^{I,\xi_1}\cdots P_{s_n}^{I,\xi_n} = \alpha_I \delta_{\zeta_I} + \bb{1-\alpha_I}\nu,
\end{equation}
where~$\zeta_I(i) = \ind\bb{p_i^\star > 0}$ for~$i\in I$,~$\nu$ is a probability measure,~$\alpha_I>0$ and~$\xi_1,\dots,\xi_n$ and~$s_1+\dots +s_n\geq 1$ are arbitrary boundary conditions and their durations. The content of~\eqref{eq:mass_leak_decomposition} is that~$\alpha_I$ is uniform in~$\eta$,~$s_1,\dots,s_n$ and the boundary conditions~$\xi_1,\dots,\xi_n$. This is useful because~$\delta_{\zeta_I}\in\mathcal{M}_\star^I$.

We verify that the restricted process has positive probability of
entering~$\zeta_I$ after any positive time, uniformly in the boundary path.
By~\eqref{eq:centered-generator-positive-coefficients}, the rate at which a
site turns to state~$1$ is~$\lambda_i\rho_i\geq p_i^\star\lambda_i\geq\varepsilon p_i^\star$,
and the pure-death component makes the rate at which it turns to state~$0$
at least~$\varepsilon$. Set~$r_i=\varepsilon p_i^\star$ if~$\zeta_I(i)=1$
and~$r_i=\varepsilon$ if~$\zeta_I(i)=0$; then~$r_i>0$ is a lower bound,
uniform in the neighbourhood configuration, for the rate at which site~$i$
takes the state~$\zeta_I(i)$. Set~$\overline\lambda=\sup_{i,\eta}\lambda_i(\eta)$.
Construct the restricted process from independent rate-$\overline\lambda$
proposal clocks at the sites of~$I$, as in the proof of
Lemma~\ref{lem:finite-propagation}. In the last unit interval,
require exactly one proposal at each site of~$I$. Conditional on the proposal
times and on the boundary path, read the proposals in time order: at each
proposal the site is in state~$\zeta_I(i)$ immediately afterwards with
conditional probability at least~$r_i/\overline\lambda$, and no later
proposal occurs at that site. Thus the configuration on~$I$ at the end of
the interval equals~$\zeta_I$ with probability at least
\begin{equation*}
\bb{\overline\lambda e^{-\overline\lambda}}^{\abs I}
\prod_{i\in I}\frac{r_i}{\overline\lambda}
=
e^{-\abs I\overline\lambda}\prod_{i\in I}r_i
=:\alpha_I>0,
\end{equation*}
uniformly in the initial configuration and the boundary path. Applying this
to the last unit interval gives~\eqref{eq:mass_leak_decomposition}.

Now we show that the process converges towards~$\mathcal{M}_\star$ from any initial configuration~$\eta$.
For~$t\geq1$ split the evolution in~\eqref{eq:revealed-boundary-product} into~${\lfloor t\rfloor}$
blocks: the first~${\lfloor t\rfloor}-1$ have length~$1$, and the last has length~$t-{\lfloor t\rfloor}+1\geq 1$.
Splitting the boundary intervals at the block endpoints gives
\begin{equation*}
\CE{f(\eta_t)}{\mathcal G_U}
=\bb{\prod_{j=1}^{{\lfloor t\rfloor}}
 \bb{P_{s_1^j}^{I,\xi_1^j}\cdots
     P_{s_{n_j}^j}^{I,\xi_{n_j}^j}}}(f)(\eta),
\end{equation*}
where~$\sum_l s_l^j=1$ for~$j<{\lfloor t\rfloor}$ and~$\sum_l s_l^{\lfloor t\rfloor}=t-{\lfloor t\rfloor}+1\geq1$.

By~\eqref{eq:mass_leak_decomposition}, each block allows us to extract
mass~$\alpha_I$ at~$\zeta_I$ from the remaining measure.
Repeating the extraction gives
\begin{align*}
\CE{\chi_A^\star(\eta_t)}{\mathcal G_U}
&=\sum_{k=1}^{{\lfloor t\rfloor}} \alpha_I(1-\alpha_I)^{k-1}
 \bb{\prod_{j=k+1}^{{\lfloor t\rfloor}}
  \bb{P_{s_1^j}^{I,\xi_1^j}\cdots
      P_{s_{n_j}^j}^{I,\xi_{n_j}^j}}}(\chi_A^\star)(\zeta_I)\\
&\quad +(1-\alpha_I)^{\lfloor t\rfloor}\nu_t(\chi_A^\star)\\
&\geq-\norm{\chi_A^\star}_\infty(1-\alpha_I)^{\lfloor t\rfloor},
\end{align*}
where~$\nu_t$ is a probability measure.
The term indexed by~$k$ is the mass first extracted at the end of
block~$k$, evolved through the remaining blocks.
The sum is non-negative because~$\delta_{\zeta_I}\in\mathcal M_\star^I$
and each~$P_s^{I,\xi}$ preserves~$\mathcal M_\star^I$.
The lower bound side tends to zero exponentially,
and hence, for every~$A\Subset\Lambda$,
\begin{equation}\label{eq:convergence-towards-m-star}
\liminf_{t\to\infty}\inf_\eta P_t(\chi_A^\star)(\eta)\geq0.
\end{equation}

To obtain convergence to~$\pi$, we combine~\eqref{eq:convergence-towards-m-star}
with Theorem~\ref{thm:common-invariant-limit}.
Fix a local function~$f$ and~$s>0$. Choose~$t_n\to\infty$ and
configurations~$\eta_n$ such that
$\abs{P_{t_n}f(\eta_n)-\pi(f)}$ converges to
$\limsup_{t\to\infty}\sup_\eta\abs{P_tf(\eta)-\pi(f)}$ and
$\delta_{\eta_n}P_{t_n-s}\Rightarrow\mu$ for some probability measure~$\mu$. This is possible because the space of probability measures on~$\Omega$ is compact.
For every finite~$B\subseteq\Lambda$,~\eqref{eq:convergence-towards-m-star} implies
\begin{equation*}
\mu(\chi_B^\star)
=\lim_{n\to\infty}P_{t_n-s}(\chi_B^\star)(\eta_n)\geq0,
\end{equation*}
so~$\mu\in\mathcal M_\star$. By the Feller property,~$P_sf$ is continuous,
and hence
\begin{equation*}
P_{t_n}f(\eta_n)
=(\delta_{\eta_n}P_{t_n-s})(P_sf)
\longrightarrow(\mu P_s)(f).
\end{equation*}
By Theorem~\ref{thm:common-invariant-limit} we obtain
\begin{equation*}
\limsup_{t\to\infty}\sup_\eta\abs{P_tf(\eta)-\pi(f)}
=\abs{(\mu P_s)(f)-\pi(f)}
\leq K_f(1+s)^D e^{-\varepsilon s/2}.
\end{equation*}
Letting~$s\to\infty$ proves convergence uniformly over initial
configurations for every local~$f$ and thus uniform ergodicity. Integrating
this uniform convergence against an invariant probability measure~$\nu$
gives~$\nu(f)=\pi(f)$ for every local~$f$, so~$\pi$ is the unique invariant
measure.
\end{proof}

\appendix

\section*{Appendix}
\addcontentsline{toc}{section}{Appendix}
\addtocontents{toc}{\protect\setcounter{tocdepth}{0}}

\section{Finite-volume approximation}
\label{app:passage-to-infinite-volume}

We first prove the finite-propagation estimate for spin systems with bounded
rates. {For qualitative convergence of finite-volume semigroups, see
also~\cite[Chapter~I, Corollary~3.14]{Liggett1985}.}

\begin{proof}[Proof of Lemma~\ref{lem:finite-propagation}]
Set
\begin{equation*}
\overline\lambda
=
\sup_{i\in\Lambda,\,\eta\in\Omega}\lambda_i(\eta),
\qquad
m
=
\sup_{i\in\Lambda}\abs{N_*(i)}.
\end{equation*}
Construct the process using rate-$\overline\lambda$ Poisson clocks at each site,
together with independent randomness deciding whether a proposed update is
accepted and, if so, the refreshed spin. Trace the dependence of the spins in~$A$ at time~$t$ backwards through this construction. Whenever a clock at~$i$ is encountered, the earlier states at the sites in~$N(i)\cup\cb{i}$ may be needed, so these sites are added to the dependence
set.

Couple the infinite-volume process and the process with zero boundary outside~$R$ using the same graphical construction in~$R$. If the backward dependence
set stays inside~$R$, the two processes agree on~$A$ at time~$t$.

To leave~$B(A,n)$, the backward dependence set must contain a time-ordered path
of at least~$n$ clock rings~\cite[proof of Proposition~3.12]{HartarskyToninelli2025}.
For a path of length~$k$, there are at most~$\abs A\,m^k$ possible spatial paths, while the expected number of
chronologically ordered~$k$-tuples of clock rings along any fixed path is~$(\overline\lambda t)^k/k!$. Hence, for every~$\theta>0$,
\begin{align*}
\pP\bb{\text{the backward dependence set leaves }B(A,n)}
&\leq
\abs A
\sum_{k\geq n}
\frac{\bb{m\overline\lambda t}^k}{k!}
\\
&\leq
\abs A
\exp\bb{
m\overline\lambda t e^\theta-\theta n
}.
\end{align*}
Choose~$v$ so that
\begin{equation*}
\theta v-m\overline\lambda e^\theta>a.
\end{equation*}
With~$n=\lceil vt\rceil$, the preceding probability is bounded by~$\abs A e^{-at}$. Therefore, for
every~$\eta\in\Omega$,
\begin{align*}
\abs{P_tf(\eta)-P_t^{R,0}f(\eta)}
&\leq
2\norm f_\infty
\pP\bb{\text{the backward dependence set leaves }R}
\\
&\leq
2\abs A e^{-at}\norm f_\infty.
\end{align*}
This proves~\eqref{eq:finite-propagation}. For~(b), if~$R\supseteq B(A,n)$,
the same bound gives
\begin{equation*}
\norm{\bb{P_t-P_t^{R,0}}f}_\infty
\leq
2\abs A\norm f_\infty
\sum_{k\geq n}
\frac{\bb{m\overline\lambda t}^k}{k!},
\end{equation*}
which tends to~$0$ as~$n\to\infty$.
\end{proof}

We now extend the patch representation to infinite volume without assuming
integrability of the signed Feynman--Kac weight. We use the weight~$W_t^\mu$ from Theorem~\ref{thm:patch-representation},
which we define from the original dual and its patch family.

\begin{remark}[Infinite volume]
\label{rem:fk-infinite-volume}
On $E_t^R$, the integrand in~\eqref{eq:monomial-feynman-kac-duality} is bounded
by
\begin{equation*}
\exp\bb{t\abs R\sup_i\abs{V_i}},
\end{equation*}
so no integrability condition is needed.
If~$\E_A\Cb{\exp\bb{\int_0^tV(A_s)\,ds}}<\infty$ for all~$A$ and~$t$, we
use Lemma~\ref{lem:finite-propagation} and dominated convergence to pass to
the limit~$R\uparrow\Lambda$ and obtain the infinite-volume
duality~\eqref{eq:intro-signed-duality}.
We do not use this condition: for patch-positive systems, nonnegativity
of~$W_t^\mu$ for~$\mu\in\mathcal M_\star$ takes its place
(Corollary~\ref{cor:infinite-volume-representation}).
\end{remark}

\begin{corollary}[Patch representation in infinite volume]
\label{cor:infinite-volume-representation}
Suppose that the spin system is patch positive, and let~$\mu\in\mathcal M_\star$.
Then, for every~$A\Subset\Lambda$ and~$t\geq0$, almost surely~$W_t^\mu\geq0$, and
\begin{equation}
\label{eq:infinite-volume-representation}
(\mu P_t)(\chi_A)
=
\E_A\Cb{W_t^\mu}.
\end{equation}
\end{corollary}

\begin{proof}{
By Proposition~\ref{prop:patch-contribution-positivity},~$W_t^\mu\geq0$.
}Let~$R_1\subseteq R_2\subseteq\cdots$ be finite sets containing~$A$ with union~$\Lambda$.
By Theorem~\ref{thm:patch-representation},
\begin{equation*}
(\mu P_t^{R_n,0})(\chi_A)
=
\E_A\Cb{W_t^\mu\ind\bb{E_t^{R_n}}}.
\end{equation*}
By Lemma~\ref{lem:finite-propagation}(b), the left-hand side converges
to~$(\mu P_t)(\chi_A)$. The events~$E_t^{R_n}$ increase in~$n$, and their union
is the event that~$\bigcup_{s\leq t}A_s$ is finite. By
Lemma~\ref{lem:dual-local-finiteness}, only finitely many successful
interactions occur up to time~$t$, each adding finitely many sites, so this
event has probability one. Since~$W_t^\mu\geq0$, monotone convergence gives
$\E_A\Cb{W_t^\mu\ind\bb{E_t^{R_n}}}\to\E_A\Cb{W_t^\mu}$.
\end{proof}

\section{Conditional factorization}
\label{app:conditional-factorization}

We give the proof of Lemma~\ref{thm:patch-factorization}. The patch laws and consistency events are defined in
Section~\ref{subsec:patch-laws}.

\begin{proof}[Proof of Lemma~\ref{thm:patch-factorization}]
Since the initial record~$(\infty,0,A)$ is deterministic under~$\pP_A$, let~$G_T$ denote the positive-time records in~$\mathcal I_T$. Then~$\mathcal G_T=\sigma(G_T)$.
We compute the joint law of~$G_T$ and the patch variables relative to a
reference measure on possible successful-interaction records. Setting all patch
test functions equal to~$1$ then gives the marginal law of~$G_T$, and the
quotient of the two densities gives the conditional law.

By Lemma~\ref{lem:dual-local-finiteness}, the records in~$G_T$ almost surely
form a finite chronological list
\begin{equation*}
g=\bb{(i_k,t_k,S_k)}_{k=1}^n,
\qquad
0<t_1<\cdots<t_n\leq T,
\end{equation*}
where~$\varnothing\neq S_k\subseteq N(i_k)$ and
$i_k\in A\cup\bigcup_{\ell<k}S_\ell$. We call the discrete sequence
$((i_k,S_k))_{k=1}^n$ admissible when these conditions hold. On the disjoint
union of all admissible record lists, define
\begin{equation}
\label{eq:patch-record-reference-measure}
\begin{aligned}
m_T(dg)
&=
\sum_{n\geq0}
\sum_{\substack{((i_k,S_k))_{k=1}^n\\\text{admissible}}}
\ind\cb{0<t_1<\cdots<t_n\leq T}\\
&\quad\times
\prod_{k=1}^n
\bb{\delta_{i_k}(S_k)+\beta_{i_k}(S_k)}
\,dt_1\cdots dt_n.
\end{aligned}
\end{equation}
Thus~$m_T$ is the reference measure on ordered admissible record lists; its
density is the product of the intensities of the proposed nonempty-target
marks, with the split and birth kinds suppressed.

Fix such a record list~$g$, and let~$\mathcal P_T(g)$ be the patch family it
defines. The interactions in~$g$ are exactly the successful interactions up to
time~$T$ if and only if every~$P\in\mathcal P_T(g)$ satisfies
$\operatorname{Con}(P)$. Indeed, suppose first that every patch is consistent.
Then there are no successful interactions that originate from the interior of a
patch, and every interaction with record in~$g$ is successful because the
consistency event for its source patch ensures it happens at an active site.
Thus~$g$ is exactly the list of records of successful interactions. Conversely,
if~$g$ is the successful-interaction record, then each~$X^P$ is the restriction
of the global active indicator to its patch, so every patch is consistent.

We can now compute the joint law. For nonnegative measurable~$h$ and~$f_P$, the
multivariate Mecke formula~\cite[Theorem~4.4]{LastPenrose2017}, applied to the
ordered nonempty-target marks, gives
\begin{equation}
\label{eq:patch-factorization-mecke}
\E_A\Cb{
h(G_T)
\prod_{P\in\mathcal P_T}f_P\bb{\Sigma_P}
}
=
\int h(g)
\prod_{P\in\mathcal P_T(g)}
\E_P\Cb{
f_P\bb{\Sigma_P}\ind_{\operatorname{Con}(P)}
}
\,m_T(dg).
\end{equation}
Here a record~$(i,t,S)$ has combined intensity
$\bb{\delta_i(S)+\beta_i(S)}\,dt$, while conditional on this record its
hidden kind has the law in~\eqref{eq:patch-initial-kind-law}. After these
records are fixed, the remaining marks on distinct patch intervals are
independent with the reference laws~$\pP_P$. The consistency indicators impose
exactly the condition that~$g$ is the successful-interaction
record.

Setting~$f_P=1$ in~\eqref{eq:patch-factorization-mecke} gives the marginal law
of the skeleton record,
\begin{equation*}
\pP_A\bb{G_T\in dg}
=
\prod_{P\in\mathcal P_T(g)}
\pP_P\bb{\operatorname{Con}(P)}
\,m_T(dg).
\end{equation*}
In particular, for~$\pP_A\circ G_T^{-1}$-almost every~$g$, every
factor~$\pP_P\bb{\operatorname{Con}(P)}$ with~$P\in\mathcal P_T(g)$ is
positive. Therefore, for~$\pP_A\circ G_T^{-1}$-almost every~$g$, the Radon--Nikodym
quotient of the joint density in~\eqref{eq:patch-factorization-mecke} by this
marginal density is
\begin{align*}
\E_A\Cb{
\prod_{P\in\mathcal P_T}f_P\bb{\Sigma_P}
\,\middle|\,G_T=g
}
&=
\prod_{P\in\mathcal P_T(g)}
\frac{
\E_P\Cb{f_P\bb{\Sigma_P}\ind_{\operatorname{Con}(P)}}
}{
\pP_P\bb{\operatorname{Con}(P)}
}
\\
&=
\prod_{P\in\mathcal P_T(g)}
\E_P^{\mathrm{con}}\Cb{f_P\bb{\Sigma_P}}.
\end{align*}
This proves~\eqref{eq:patch-factorization} for nonnegative factors and
identifies the conditional law. The integrable signed case follows from this
conditional-law statement.
\end{proof}

\section*{Statement on AI use}
We developed the mathematics underlying the patch representation without AI
assistance. We used Claude and ChatGPT to assist with planning the exposition and formatting, prepare preliminary
drafts, check the proofs for gaps and errors, search and review the literature, and produce the figures. The final manuscript is almost entirely human-written
and we take full responsibility for the mathematical content and the exposition.

\addtocontents{toc}{\protect\setcounter{tocdepth}{1}}
\bibliographystyle{amsplain-fullnames}
\bibliography{references}

\begin{acks}
MG has been supported by NSF award DMS-2348113.
We thank Marek Biskup for discussions.
\end{acks}

\end{document}